\documentclass[12pt,reqno]{amsart}
\usepackage{amsmath}
\usepackage{amssymb}
\usepackage{verbatim}
\usepackage{mathrsfs}
\usepackage{graphicx} 
\usepackage{caption}
\usepackage{subcaption}
\usepackage{hyperref}
\usepackage{color}
\usepackage[noadjust]{cite}

\usepackage[top=0.5in,bottom=0.5in,left=0.7in,right=0.7in]{geometry}

\newtheorem{lemma}{Lemma}
\newtheorem{prop}[lemma]{Proposition}

\newtheorem{theorem}[lemma]{Theorem}

\newtheorem{rem}[lemma]{Remark}
\newtheorem{exmp}{Example}

\newtheorem{algorithm}{Algorithm}

\numberwithin{equation}{section}
\numberwithin{lemma}{section}

\newcommand{\C}{\mathbb{C}}    
\newcommand{\N}{\mathbb{N}}    
\newcommand{\R}{\mathbb{R}}    
\newcommand{\T}{\mathbb{T}}    
\newcommand{\Z}{\mathbb{Z}}    

\newcommand{\wh}{\widehat}

\renewcommand{\le}{\leqslant}
\renewcommand{\ge}{\geqslant}
\newcommand{\bs}{\backslash}
\newcommand{\ol}{\overline}
\newcommand{\la}{\langle}
\newcommand{\ra}{\rangle}

\newcommand{\eps}{\epsilon}

\newcommand{\pa}{\mathsf{a}}
\newcommand{\pb}{\mathsf{b}}

\newcommand{\pd}{\mathsf{d}}
\newcommand{\pe}{\mathsf{e}}
\newcommand{\pu}{\mathsf{u}}
\newcommand{\pv}{\mathsf{v}}
\newcommand{\pw}{\mathsf{w}}
\newcommand{\pp}{\mathsf{p}}
\newcommand{\pq}{\mathsf{q}}
\newcommand{\pr}{\mathsf{r}}
\newcommand{\ps}{\mathsf{s}}
\newcommand{\pt}{\mathsf{t}}
\newcommand{\pA}{\mathsf{A}}
\newcommand{\pB}{\mathsf{B}}

\newcommand{\pD}{\mathsf{D}}
\newcommand{\pE}{\mathsf{E}}
\newcommand{\pF}{\mathsf{F}}
\newcommand{\pI}{\mathsf{I}}
\newcommand{\pN}{\mathsf{N}}
\newcommand{\pU}{\mathsf{U}}
\newcommand{\pV}{\mathsf{V}}
\newcommand{\pW}{\mathsf{W}}
\newcommand{\pP}{\mathsf{P}}
\newcommand{\pzr}{\mathbf{0}}

\newcommand{\ptA}{\tilde{\mathsf{A}}}
\newcommand{\ptP}{\tilde{\mathsf{P}}}
\newcommand{\ptU}{\tilde{\mathsf{U}}}
\newcommand{\pQ}{\mathsf{Q}}
\newcommand{\pR}{\mathsf{R}}

\newcommand{\pth}{\pmb{\theta}}

\newcommand{\pTh}{\mathsf{\Theta}}

\newcommand{\ptb}{\tilde{\mathsf{b}}}
\newcommand{\ptd}{\tilde{\mathsf{d}}}

\newcommand{\ptp}{\tilde{\mathsf{p}}}
\newcommand{\ptq}{\tilde{\mathsf{q}}}
\newcommand{\brpA}{\breve{\mathsf{A}}}

\newcommand{\brpU}{\breve{\mathsf{U}}}

\newcommand{\mrb}{{\mathring{b}}}

\newcommand{\mrpq}{\mathring{\mathsf{q}}}
\newcommand{\mrpA}{\mathring{\mathsf{A}}}

\newcommand{\mrpU}{\mathring{\mathsf{U}}}

\newcommand{\mrpb}{\mathring{\mathsf{b}}}
\newcommand{\mrpd}{\mathring{\mathsf{d}}}

\newcommand{\DG}{{\text{Diag}}}
\newcommand{\mz}{\mathsf{Z}} 

\newcommand{\er}{\eqref}

\newcommand{\bp}{ \begin{proof} }
	\newcommand{\ep}{\end{proof} }
\newcommand{\be}{ \begin{equation} }
\newcommand{\ee}{ \end{equation} }
\newcommand{\tp}{\mathsf{T}}

\newcommand{\sr}{\operatorname{sr}}  
\newcommand{\sm}{\operatorname{sm}}  
\newcommand{\vmo}{\operatorname{vm}}

\newcommand{\lp}[1]{l_{#1}(\mathbb{Z})}

\newcommand{\lrs}[3]{(l_{#1}(\mathbb{Z}))^{#2\times #3}}

\newcommand{\setsp}{\;:\;}     
\newcommand{\td}{\boldsymbol{\delta}}  
\newcommand{\Lp}[1]{L_{#1}(\mathbb{R})}
\newcommand{\cM}{\mathcal{M}}

\newcommand{\cN}{\mathcal{N}}
\newcommand{\JJ}{\begin{bmatrix}1 & \\ & -1\end{bmatrix}}
\newcommand{\adj}{\operatorname{adj}}

\newcommand{\sym}{\mathsf{S}}   

\newcommand{\re}{\operatorname{Re}}
\newcommand{\im}{\operatorname{Im}}
\newcommand{\odd}{\operatorname{odd}}
\newcommand{\ldeg}{\operatorname{ldeg}}
\newcommand{\len}{\operatorname{len}}
\newcommand{\fsupp}{\operatorname{fsupp}}

\begin{document}
\title[Quasi-Tight Framelets with Symmetry]
{Generalized Matrix Spectral Factorization with Symmetry and Applications to Symmetric Quasi-Tight Framelets}

\author{Chenzhe Diao}
\author{Bin Han}
\address{Department of Mathematical and Statistical Sciences,
	University of Alberta, Edmonton,\quad Alberta, Canada T6G 2G1.
	\quad {\tt diao@ualberta.ca},\quad  {\tt bhan@ualberta.ca} }

\author{Ran Lu}
\address{College of Science, Hohai University, Nanjing, China 211100.
\quad  {\tt rlu3@hhu.edu.cn}}

\thanks{The research of the first and the second authors was supported in part by the Natural Sciences and Engineering Research Council of Canada (NSERC) under Grant RGPIN-2019-04276.
}	

\thanks{The research of the third author was supported by the Hohai University start-up funds.
}

\makeatletter \@addtoreset{equation}{section} \makeatother
\begin{abstract}

Factorization of matrices of Laurent polynomials plays an important role
in mathematics and engineering such as wavelet frame construction and filter bank design.
Wavelet frames (a.k.a. framelets) are useful in applications such as signal and image processing.
Motivated by the recent development of quasi-tight framelets,
we study and characterize generalized spectral factorizations with symmetry for $2\times 2$ matrices of Laurent polynomials.
Applying our result on generalized matrix spectral factorization, we establish a necessary and sufficient condition for the existence of symmetric quasi-tight framelets with two generators. The proofs of all our main results are constructive and therefore, one can use them as construction algorithms.
We provide several examples to illustrate our theoretical results on generalized matrix spectral factorization and quasi-tight framelets with symmetry.

\end{abstract}


 \maketitle

\pagenumbering{arabic}


\section{Introduction, Motivations and Main Results}\label{sec:introduction}

Wavelet frames (also called as framelets) are often constructed
through their associated filter banks
from refinable functions, which further boil down to the fundamental problem of factorizing matrices of Laurent/trigonometric polynomials with various properties. We refer the readers to \cite{cs08,ch00,chs02,daubook,dh04,dhrs03,dhacha,dh18pp,ds07,eh08,
han97,han09,han13,han14acha,hanbook,hl19pp,hl20pp,hm04,hm05,hr08,jiang03,js15,kps16book,ml11,rs97,sz18,sel01,sa04,sl11} for more details on framelets. As a consequence, appropriately factorizing matrices of Laurent polynomials plays a crucial role in the construction of framelets with desired properties such as symmetry and vanishing moments.
Both symmetry and vanishing moments are highly desirable among many key properties of framelets for sparse representations and boundary treatment.
This motivates us to study generalized matrix spectral factorization with symmetry.

\subsection{Main result on generalized matrix spectral factorization}

Let $\pu(z)=\sum_{k\in\Z}u(k)z^k, z\in \C\bs\{0\}$ be a \emph{Laurent polynomial}, where its coefficient sequence $u=\{u(k)\}_{k\in \Z}$ is finitely supported with $u(k)\in \C$ for all $k\in \Z$ and $u$ is often called a filter in the literature of engineering. We say that $\pu$ has  \emph{symmetry} if its coefficient filter $u$ has symmetry:
\be \label{eq:RealSymTime}
\pu(z)=\epsilon z^c \pu(z^{-1})
\quad \mbox{or equivalently},\quad
u(k) = \epsilon u(c-k), \qquad \forall~ k\in \Z,
\ee
for some $ \epsilon \in \{-1, 1\} $ and $c\in \Z$.
For any nonzero Laurent polynomial $\pu$ having symmetry, we define the following symmetry operator $\sym$:
\be\label{eq:RealSymZ} \sym \pu (z) := \frac{\pu (z)}{\pu(z^{-1})},  \qquad z\in \C\bs\{0\}.
\ee
It is then easy to see that a Laurent polynomial $\pu$ has symmetry as in \er{eq:RealSymTime} if and only if $\sym\pu(z)=\eps z^c$. In this case, we say that $\pu$ \emph{has symmetry with type $\eps z^c$}. For the zero Laurent polynomial $\mathbf{0}$, we say that it \emph{has symmetry with any type}.

For any $t\times r$ matrix $\pP(z)=\sum_{k\in\Z}P(k)z^k$ of Laurent polynomials, where $P(k)\in\C^{t\times r}$ for all $k\in\Z$, we define the \emph{Hermitian conjugate} of $\pP$ via
$$
\pP^\star(z):=\ol{\pP(\ol{z^{-1}})}^\tp=\sum_{k\in\Z}\ol{P(k)}^\tp z^{-k},\quad z\in\T:=\{\zeta\in \C \setsp |\zeta|=1\}.
$$
Let $\pA(z)$ be a \emph{Hermitian} $2\times 2$ matrix of Laurent polynomials, i.e. $\pA(z)=\pA^\star(z)$. Moreover, suppose that all elements of $\pA$ have symmetry. The problem of generalized matrix spectral factorization with symmetry is to decompose the matrix $ \pA(z) $ into the following form:
\be \label{eq:spectfact:0}
\pA(z) = \pU(z)\DG(\epsilon_1,\epsilon_2) \pU^\star(z)
\quad \mbox{with}\quad
\epsilon_1,\epsilon_2\in \{-1,1\} \qquad
\mbox{for all}~ z\in \T,
\ee
where $\pU$ is some $2\times 2$ matrix of Laurent polynomials such that all the entries of $\pU$ have symmetry.
The matrix spectral factorization with symmetry corresponds to \eqref{eq:spectfact:0} with $\epsilon_1=\epsilon_2=1$ has been characterized in \cite[Theorem~2.3]{hm04} and \cite[Theorems~3.2]{han13}.
\eqref{eq:spectfact:0} without requiring the symmetry property of $\pU$ for $\epsilon_1=\epsilon_2=1$
is just the well-known standard matrix spectral factorization (also called matrix Fej\'er-Riesz lemma, e.g., see \cite{hjs04}) in the literature. However, the case $\epsilon_1\ne \epsilon_2$ in \eqref{eq:spectfact:0} is much more difficult and technical, because we lose the positive semidefinite property $\pA(z)\ge 0$ for all $z\in \T$ and we have to deal with a difference of Hermitian squares instead of sum of Hermitian squares as in \cite{han13,hm04}.
For a Laurent polynomial $\pu$, we say that $\pu$ has the \emph{difference of (Hermitian) squares (DOS) property with respect to symmetry type $\eps z^c$} with $\eps\in\{\pm1\}$ and $c\in\Z$, if there exist Laurent polynomials $\pu_1$ and $\pu_2$ having symmetry, such that
\be\label{dos}
\pu_1(z)\pu_1^\star(z)-\pu_2(z)\pu_2^\star(z)=\pu(z)\quad\mbox{and}\quad \frac{\sym\pu_1(z)}{\sym\pu_2(z)}=\eps z^c.
\ee
For $z_0\in \C\bs\{0\}$,
we denote by $\mz(\pu,z_0)$ the multiplicity of the zero $z_0$ of the Laurent polynomial $\pu$.

We now complete the picture on generalized matrix spectral factorization by
stating the following main theorem on the generalized spectral factorization of $2\times 2$ matrices of Laurent polynomials.

\begin{theorem} \label{thm:SymSpecFactSym}
Let $ \pA $ be a $ 2\times 2 $ Hermitian matrix of Laurent polynomials having symmetry and define $\alpha(z):=\sym \pA_{1,2}(z)$, where $\pA_{j,k}$ is the $(j,k)$-entry of $\pA$ for $j,k\in \{1,2\}$.
	Then there is a $2\times 2$ matrix $ \pU(z)$ of Laurent polynomials with all its entries having symmetry
	such that
	$\pA(z) = \pU(z)\DG(1,-1)
\pU^\star(z)$ holds (i.e., \eqref{eq:spectfact:0} holds with $\epsilon_1=1$ and $\epsilon_2=-1$) and
	the symmetry type of $\pU$ satisfies
\begin{equation} \label{eq:Sym1}	 \frac{\sym\pU_{1,1}(z)}{\sym\pU_{2,1}(z)} = \frac{\sym\pU_{1,2}(z)}{\sym\pU_{2,2}(z)} =
\alpha(z),
\end{equation}
%
	if and only if the following two conditions hold:
	\begin{enumerate}
		\item  $ \det(\pA(z)) = -\pd(z)\pd^\star(z) $ for some Laurent polynomial $ \pd$ with symmetry.
		
		\item Define $ \pp_0(z):=\gcd(\pA_{1,1}(z), \pA_{1,2}(z), \pA_{2,1}(z), \pA_{2,2}(z)) $ and
		\be\label{pp}\pp(z):= \frac{\pp_0(z)}{(z-1)^{\mz(\pp_0, 1)} (z+1)^{\mz(\pp_0, -1)}}.
\ee
Then $ \pp(z) $ satisfies the DOS (Difference of Squares) condition  with respect to type
$ \alpha(z) \sym\pd(z) $.
	\end{enumerate}
\end{theorem}


The technical long proof of Theorem~\ref{thm:SymSpecFactSym} will be given in Section~\ref{sec:fac:lau}, where we shall also explain the compatibility condition in \eqref{eq:Sym1} for the symmetry patterns between the matrices $\pA$ and $\pU$.


Writing a Laurent polynomial into sums and differences of Hermitian squares is inevitable for a generalized spectral factorization problem in Theorem~\ref{thm:SymSpecFactSym}. We now characterize the difference of (Hermitian) squares (DOS) property below whose proof is given in Subsection~\ref{subsec:dos}.

\begin{theorem}\label{thm:dos}
Suppose $\eps\in\{\pm 1\}$ and $c\in\Z$. A Laurent polynomial $\pu$ has the DOS property with respect to the symmetry type $\eps z^c$ if and only if
\begin{enumerate}
	\item[(1)] The Laurent polynomial $\pu$ has real coefficients and $\pu^\star=\pu$;
	
	\item[(2)] The Laurent polynomial $\pu$ satisfies exactly one of the following technical conditions:
	
	\begin{enumerate}
		\item[(i)] If $\eps=1$ and $c\in2\Z$, then there is no condition on $\pu$;
		
		\item[(ii)]If $\eps=1$ and $c\in2\Z+1$, then $\mz(\pu,x)\in2\Z$ for all $x\in(-1,0)$;
		
		\item[(iii)]If $\eps=-1$ and $c\in2\Z$, then $\mz(\pu,x)\in2\Z$ for all $x\in(-1,0)\cup (0,1)$;
		
		\item[(iv)]If $\eps=-1$ and $c\in2\Z+1$, then $\mz(\pu,x)\in2\Z$ for all $x\in(0,1)$.
	\end{enumerate}
\end{enumerate}
\end{theorem}

\subsection{Application of generalized matrix spectral factorization on constructing quasi-tight framelets}

Symmetric wavelets and framelets have been extensively studied and applied in applications, e.g., see \cite{ch00,dh04,ds07,eh08,han97,han09,han13,han14acha,hanbook,
hm04,hm05,hz10,jiang03,js15,kps16book,ml11,mz12,rs97,sz18,sel01,
sa04,sl11,zhuang12} and references therein. Construction of wavelets and framelets is intrinsically connected with factorizing matrices of Laurent polynomials. In this paper, we are interested in \emph{quasi-tight framelets} with two generators having the symmetry property. Our investigation is inspired by the recent development of quasi-tight framelets \cite{dhacha,dh18pp,hanbook,hl19pp,hl20pp} and the importance of (anti-)symmetric wavelet frames in many applications. Using Theorem~\ref{thm:SymSpecFactSym}, we can obtain a necessary and sufficient condition for the existence of quasi-tight framelets with symmetry.

To present our main result on quasi-tight framelets with symmetry, we need several notations. By $ \lrs{0}{r}{s} $ we denote the space of all finitely supported matrix-valued sequences/filters on $ \Z $, i.e., $u\in\lrs{0}{r}{s}$ if $u:\Z\to\C^{r\times s}$ has only finitely many nonzero terms. For $ u = \{u(k)\}_{k\in \Z}\in \lrs{0}{r}{s}$, its \emph{associated matrix Laurent polynomial} is $\pu(z):= \sum_{k\in \Z} u(k) z^k$ for $z\in \C \bs \{0\}$.
Define $u^\star:=\ol{u(-\cdot)}^{\tp}$ and then
$\pu^\star(z):=[\pu(z)]^{\star}=\sum_{k\in \Z} \ol{u(k)}^\tp z^{-k}$.
For $u\in\lrs{0}{s}{r}$ and $v\in\lrs{0}{t}{r}$, their \emph{convolution} is
$u*v:=\sum_{j\in\Z}u(\cdot-j)v(j)$ and the matrix Laurent polynomial associated with $u*v$ is $\pu(z)\pv(z)$.

One important feature for a framelet is the orders of vanishing moments of its generators. For any compactly supported $\psi\in(\Lp{2})^s$, we say that $\psi$ has \emph{order $m$ vanishing moments} if
$$
\int_{\R} x^j \psi(x) dx = 0 ,\qquad j = 0, \ldots, m-1.
$$
We define $ \vmo(\psi) := m$ with $m$ being the largest positive integer satisfying the above identities.
For a filter $ b \in \lp{0} $, we say that $ b $ has \emph{order $m$ vanishing moments} if $ (z-1)^m $ divides the Laurent polynomial $\pb(z)$, and we similarly define $ \vmo(b) := \vmo(\pb) := m$ with $m$ being the largest positive integer satisfying this.
The order of vanishing moments is  naturally linked to the concept of the sum rules of a filter. For a sequence $ a \in \lp{0} $, we say that $ a $ has $ n $ \emph{sum rules} if $ (z+1)^n $ divides the Laurent polynomial $ \pa(z) $. We define $ \sr(a):= \sr(\pa(z)):= n $, with $ n $ being the largest such integer. For an integer $k\in \Z$, we define
\be\label{odd}
\odd(k):=\frac{1-(-1)^k}{2}=
\begin{cases}0, &\text{if $k$ is even},\\
1, &\text{if $k$ is odd},\end{cases}
\qquad k\in \Z.
\ee

As an application of Theorem~\ref{thm:SymSpecFactSym},
our main result on symmetric quasi-tight framelets with two generators is the following result, whose proof is given in Section~\ref{sec:qtf}.


\begin{theorem} \label{thm:qtfsym}
	Let $ \Theta, a \in \lp{0} $ be given filters having symmetry $\sym \pTh(z) = 1$ and $ \sym \pa(z) = z^c $ for some $c\in \Z$, and $\Theta^\star=\Theta$.
Let $ n_b $ be any chosen integer satisfying
	\be\label{eq:nbRange}
	1\le n_b\le\min\{\sr(a),\tfrac{1}{2}\vmo(\pTh(z)-\pTh(z^2)\pa^\star(z)\pa(z))\}.
	\ee
	Define $2\times 2$ matrices $\cM_{\pa, \pTh|n_b}(z)$ and $\cN_{\pa,\pTh|n_b}$ of Laurent polynomials as follows:
	\be\label{eq:DefcMnb}\cM_{\pa, \pTh|n_b}(z):=\DG((1-z^{-1})^{-n_b},
(1+z^{-1})^{-n_b})\cM_{\pa,\pTh}(z)\DG((1-z)^{-n_b},(1+z)^{-n_b}),\ee
	\be\label{eq:DefcN}\cN_{\pa, \pTh|n_b}(z^2)=\begin{bmatrix}1 & z^{-1}\\
		1 &-z^{-1}\end{bmatrix}^\star\cM_{\pa,\pTh|n_b}(z)\begin{bmatrix}1 & z^{-1}\\
		1 &-z^{-1}\end{bmatrix},\ee
	where the $2\times 2$ matrix $\cM_{\pa,\pTh}$ of Laurent polynomials is defined to be
	\be
	\label{cond:oep:tf}
\cM_{\pa,\pTh}(z):=\begin{bmatrix}\pTh(z)-\pTh(z^2)\pa^\star(z)\pa(z) & -\pTh(z^2)\pa^\star(z)\pa(-z)\\
		-\pTh(z^2)\pa^\star(-z)\pa(z) &\pTh(-z)-\pTh(z^2)\pa^\star(-z)\pa(-z)\end{bmatrix}.
	\ee
Let $\eps_1=1$ and $\eps_2=-1$.
	Then there exist filters $b_1, b_2\in \lp{0}$ such that $b_1, b_2$ have symmetry with order $n_b$ vanishing moments and $\{a; b_1, b_2\}_{\Theta, (\eps_1, \eps_2)}$
forms a quasi-tight framelet filter bank, i.e.,
	\begin{align}
	&\pTh(z^2) \pa^\star(z)\pa(z) +\eps_1 \pb_1^\star(z)\pb_1(z)+\eps_2\pb_2^\star(z)\pb_2(z)=\pTh(z),
	\qquad z\in \C\bs\{0\}, \label{tffb:1}\\
	&\pTh(z^2) \pa^\star(z)\pa(-z) +\eps_1\pb_1^\star(z)\pb_1(-z)+\eps_2\pb_2^\star(z)\pb_2(-z)=0,
	\qquad z\in \C\bs\{0\},\label{tffb:0}
	\end{align}
if and only if the following two conditions hold:
	\begin{enumerate}
		\item $ \det(\cN_{\pa, \pTh|n_b}(z)) = -\pd_{n_b}(z)\pd_{n_b}^\star(z) $ for some Laurent polynomial $ \pd_{n_b}(z) $ having symmetry.
		\item
		$ \pp(z) $ has the difference of (Hermitian) squares (DOS) property with respect to symmetry type $ (-1)^{c+n_b}z^{\odd(c+n_b)-1}\sym \pd_{n_b}(z) $, where $ \pp(z) $ is defined through
		\begin{align}
		\pp_0(z):=&\gcd\left(\left[ \cN_{\pa, \pTh|n_b}(z) \right]_{1,1},
		\left[ \cN_{\pa, \pTh|n_b}(z) \right]_{1,2},
		\left[ \cN_{\pa, \pTh|n_b}(z) \right]_{2,1},
		\left[ \cN_{\pa, \pTh|n_b}(z) \right]_{2,2}
		\right),\label{pp0:z}\\
		\pp(z):=& \frac{\pp_0(z)}{(z-1)^{\mz(\pp_0, 1)}(z+1)^{\mz(\pp_0, -1)}}.\label{pp1:z}
		\end{align}
		Moreover, if $\pTh(z)=1$, then $ \pp_0(z) = \pp(z) = 1 $ and item (2) is automatically satisfied.
	\end{enumerate}


\noindent
If in addition $\pa(1)=\pTh(1)=1$ and
	$\phi\in \Lp{2}$, where
	$\wh{\phi}(\xi):=\prod_{j=1}^\infty \pa(e^{-i 2^{-j}\xi})$ for $\xi\in \R$, define
\be \label{eq:psi}
	\eta := \sum_{k\in \Z} \Theta(k) \phi(\cdot - k) \quad \mbox{and}\quad
	\psi^\ell:= 2\sum_{k\in \Z} b_\ell(k) \phi(2\cdot - k)
\ee
for $\ell=1,2$,
then all $\eta, \psi^1, \psi^2\in \Lp{2}$ have symmetry with $\psi^1, \psi^2$ having at least order $n_b$ vanishing moments, and $\{\eta,\phi;\psi^1, \psi^2\}_{(\eps_1,\eps_2)}$ forms a compactly supported quasi-tight framelet in $\Lp{2}$, i.e.,
	\be\label{qtfr}f =  \sum_{k\in \Z}  \la f, \eta(\cdot - k) \ra \phi(\cdot - k) +\sum_{\ell=1}^2
	\sum_{j = 0}^{\infty} \sum_{k\in \Z}\eps_\ell \la f, \psi^\ell_{2^j;k} \ra \psi^\ell_{2^j;k},
	\qquad \forall~ f\in \Lp{2}
	\ee
with the above series converging unconditionally in $\Lp{2}$, where $\psi^\ell_{2^j;k}:=2^{j/2}\psi^\ell(2^j\cdot-k)$.
\end{theorem}

\subsection{Paper structure}

The paper is organized as follows.
In Section~\ref{sec:ex}, we briefly review some necessary background on framelets. In particular, we shall explain more on our motivations of considering quasi-tight framelets with symmetry, and some illustrative examples will be provided.
In Section~\ref{sec:fac:lau}, we introduce some notations for Laurent polynomials with symmetry first. Next, we consider two special cases of Theorem~\ref{thm:SymSpecFactSym} before we prove Theorem~\ref{thm:SymSpecFactSym} on the factorization of a $ 2 \times 2 $ matrix of Laurent polynomials with symmetry.
Finally, in Section~\ref{sec:qtf}, we use Theorem~\ref{thm:SymSpecFactSym} to prove Theorem~\ref{thm:qtfsym}, which characterizes quasi-tight framelet filter banks $ \{a; b_1,b_2\}_{\Theta, (1, -1)} $ with symmetry, and their associated symmetric quasi-tight framelets $ \{\eta, \phi; \psi^1, \psi^2 \}_{(1, -1)} $ in $\Lp{2}$.

\section{Motivation and Examples of Quasi-tight Framelets with Symmetry}\label{sec:ex}

To further explain our motivations, in this section we shall first review some necessary background on framelets and then provide several illustrative examples of symmetric quasi-tight framelets as an application of Theorems~\ref{thm:SymSpecFactSym} and~\ref{thm:qtfsym}.

\subsection{Introduction to quasi-tight framelets}

For $\phi,\psi^1,\ldots,\psi^s\in\Lp{2}$, we say that $\{ \phi; \psi^1,\ldots,\psi^s\}$ is a \emph{framelet} (i.e., a wavelet frame) in $\Lp{2}$ if there exist positive constants $ C_1, C_2 > 0 $ such that
\[
C_1 \| f \|_{\Lp{2}}^2 \leqslant
\sum_{k\in \Z} | \la f, \phi(\cdot - k) \ra |^2 +
\sum_{\ell = 1}^{s}\sum_{j = 0}^{\infty}\sum_{k\in \Z} |\la f, \psi^\ell_{2^j;k}\ra|^2
\leqslant C_2 \| f \|_{\Lp{2}}^2 , \qquad
\forall ~f\in \Lp{2},
\]
where $\psi^\ell_{2^j;k}:=2^{j/2}\psi^\ell(2^j\cdot-k)$.
For $\phi,\eta,\psi^1,\dots,\psi^s\in\Lp{2}$ and $\eps_1,\dots,\eps_s\in\{\pm 1\}$,
we say that $\{\eta, \phi;\psi^1,\ldots,\psi^s\}_{(\eps_1,\dots,\eps_s)}$ is \emph{a quasi-tight framelet} in $\Lp{2}$ if $\{\eta;\psi^1,\ldots,\psi^s\}$ and $\{\phi;\psi^1,\ldots,\psi^s\}$ are framelets in $\Lp{2}$ and
\be\label{qtfr}f =  \sum_{k\in \Z}  \la f, \eta(\cdot - k) \ra \phi(\cdot - k) +
\sum_{\ell = 1}^{s} \sum_{j = 0}^{\infty}\sum_{k\in \Z}\eps_\ell \la f, \psi^\ell_{2^j;k} \ra \psi^\ell_{2^j;k},
\qquad \forall~ f\in \Lp{2},
\ee
with the above series converging unconditionally in $ \Lp{2} $. We say that $\{\phi;\psi^1,\dots,\psi^s\}_{(\eps_1,\dots,\eps_s)}$ is a quasi-tight framelet in $\Lp{2}$ if
$\{\phi, \phi;\psi^1,\dots,\psi^s\}_{(\eps_1,\dots,\eps_s)}$ is a quasi-tight framelet in $\Lp{2}$. By \cite[Propositions 4 and~5]{han12} or \cite[Proposition 4.3]{han17}, \er{qtfr} implies that $\{\psi^1,\dots,\psi^s\}_{(\eps_1,\dots,\eps_s)}$ is \emph{a homogeneous quasi-tight framelet} in $\Lp{2}$, that is,
\be\label{qtfr}f =
\sum_{\ell = 1}^{s} \sum_{j\in\Z}\sum_{k\in \Z}\eps_\ell \la f, \psi^\ell_{2^j;k} \ra \psi^\ell_{2^j;k},
\qquad \forall~ f\in \Lp{2},
\ee
with the above series converging unconditionally in $ \Lp{2} $. A quasi-tight framelet $\{\phi;\psi^1,\ldots,\psi^s\}_{(\eps_1,\ldots,\eps_s)}$ in $\Lp{2}$ with $\eps_1=\dots=\eps_s=1$
is often called \emph{a tight framelet} in the literature. In this case, $\{\psi^1,\dots,\psi^s\}$ is \emph{a homogeneous tight framelet} in $\Lp{2}$.

In practice, compactly supported wavelets and framelets are highly desired to reduce computational complexity. Quite often, we derive them from a compactly supported \emph{refinable function} $\phi$, i.e.,
%
\begin{equation} \label{eq:refinement}
\phi = 2 \sum_{k\in \Z} a(k) \phi(2\cdot - k)
\quad \mbox{or equivalently},\quad \wh{\phi}(2\xi)=\pa(e^{-i\xi}) \wh{\phi}(\xi),
\end{equation}
holds for some finitely supported sequence $ a\in \lp{0} $, where the Fourier transform is defined as $ \wh{f}(\xi) := \int_{\R} f(x)e^{-ix\xi} dx $
for $ f \in \Lp{1} $ and can be naturally extended to square integrable functions and tempered distributions.
For a filter $ a\in \lp{0} $ satisfying $ \pa(1) = 1 $ (often called
\emph{a refinement mask} or \emph{a low-pass filter} in the literature), it is easy to check that
\eqref{eq:refinement} holds for a compactly supported function/distribution $\phi$ by defining $\phi$ through its Fourier transform
\be \label{def:phi}
\wh{\phi}(\xi) := \prod_{j = 1}^{\infty} \pa(e^{-i 2^{-j}\xi}), \qquad \xi \in \R.
\ee
%

From a given refinable function $\phi\in\Lp{2}$ associated with a low-pass filter $a\in \lp{0}$, to construct a quasi-tight framelet, we first find an appropriate filter $\Theta\in\lp{0}$ with $\Theta^\star=\Theta$ and $\pTh(1)=1$, which is often called \emph{a moment correcting filter}. Define $\cM_{\pa,\pTh}$ as in \er{cond:oep:tf}. Next, we want to find $b=(b_1,\ldots,b_s)^\tp \in\lrs{0}{s}{1}$ and $\eps_1,\dots,\eps_s\in\{\pm1\}$ such that $\pb_1(1)=\cdots=\pb_s(1)=0$ and
\be\label{fac:1}\cM_{\pa,\pTh}(z)=\begin{bmatrix}\pb(z)& \pb(-z)\end{bmatrix}^{\star}\DG(\eps_1,\dots,\eps_s)\begin{bmatrix}\pb(z),\pb(-z)\end{bmatrix}.
\ee
$\{a;b_1,\ldots,b_s\}_{\Theta, (\eps_1,\dots,\eps_s)}$
satisfying \eqref{fac:1} is called
\emph{a quasi-tight framelet filter bank}. Define $\eta$ and $\psi^\ell$ as in \er{eq:psi} for $\ell=1,\ldots,s$. If $\phi\in \Lp{2}$,
then $\{\eta, \phi;\psi^1,\ldots,\psi^s\}_{(\eps_1,\dots,\eps_s)}$ is a quasi-tight framelet in $\Lp{2}$ satisfying $\wh{\psi^\ell}(0)=0$ for all $\ell=1,\ldots,s$. This approach of constructing quasi-tight framelets is known as a special case of \emph{the Oblique Extension Principle (OEP)} in the literature, and such constructed framelets are called \emph{OEP-based framelets}, see \cite{chs02,dh04,dhrs03,hanbook} and references therein for details.

To construct a symmetric quasi-tight framelet, we require that our low-pass filter $a\in\lp{0}$ should have symmetry. For the construction of an orthogonal wavelet from $\phi$, its associated refinement mask $a\in\lp{0}$ must satisfy (\cite{daubook})
\be\label{sym:wav}
|\pa(z)|^2+|\pa(-z)|^2=1,\qquad \forall z\in\T.
\ee
As pointed out by Daubechies in \cite{daubook}, there does not exist a compactly supported real-valued (dyadic) orthonormal wavelet with symmetry, except for the Haar wavelet which is discontinuous.
Due to the aforementioned restrictions \eqref{sym:wav} on compactly supported (dyadic) real-valued orthogonal wavelets, it is natural to consider tight framelets with symmetry, which only require the following less restrictive condition (\cite{dhrs03,rs97}):
\be\label{sym:tf}
|\pa(z)|^2 + |\pa(-z)|^2 \leqslant 1, \qquad \forall z \in \T.
\ee
The redundancy of tight framelets brings much more freedom for construction and robustness for their performance in applications. To achieve lower computational complexity in the framelet transform, researchers are trying to find tight framelets with small number of generators. Symmetric tight framelet with three generators have been obtained in \cite{ch00,han14acha,hm05}.
As for tight framelets with two generators and symmetry, the general characterization based on the OEP was given in \cite[Theorem~2.4]{hm04} and \cite[Theorem~4.2]{han13}. It is shown in \cite{hm04,han13} that the construction
(i.e., \er{fac:1} with $s=2$ and $\eps_1=\eps_2=1$)
is related to the problem of splitting a $ 2\times 2 $ Hermitian matrix $ \pA(z) $ of Laurent polynomials with symmetry into the form of
\begin{equation} \label{eq:pdFactor}
\pA(z) = \pU(z) \pU^\star(z), \qquad z \in \C \setminus \{0\},
\end{equation}
where $ \pU(z) $ is a $ 2\times 2 $ matrix of Laurent polynomials, and all of its elements have symmetry.
The necessary and sufficient condition for such a factorization to exist is given in \cite[Theorem~2.3]{hm04} and \cite[Theorem~3.2]{han13}.
Notice that the factorization in \eqref{eq:pdFactor} above is similar to the problem of solving the factorization in \eqref{eq:spectfact:0} we stated at the beginning of this paper. The only difference is that \eqref{eq:pdFactor} is a standard spectral factorization, wherein \eqref{eq:spectfact:0} we have a generalized spectral factorization of $\pA$. However, our generalized matrix spectral factorization problem has some fundamental differences and essential difficulties compared to the standard one. For example, \eqref{eq:pdFactor} implies that $\pA(z)$ is positive semidefinite for all $z\in \T$. We will further explain this issue in later sections.

On the other hand, the sparsity of the framelet expansion \er{qtfr} is determined by vanishing moments of the framelet generator $\psi$.  By a simple argument (e.g., see \cite[Proposition~3.3.1]{hanbook}), one has
$$
\min\{\vmo(b_1),\ldots,\vmo(b_s)\}=\min\{\vmo(\psi^1),\ldots,\vmo(\psi^s)\}\le \min\left\{\sr(a), \tfrac{1}{2}\vmo\left(\pTh(z)-\pTh(z^2)\pa^\star(z)\pa(z)\right)\right\},
$$
no matter how we choose $\pTh$. Thus we always try to find $\pTh$ such that $\vmo\left(\pTh(z)-\pTh(z^2)\pa^\star(z)\pa(z)\right)$ is as large as possible, and make $\min\{\vmo(\psi^1),\ldots,\vmo(\psi^s)\}$ also large.

From the above discussions, we see that constructing a tight framelet with symmetry and high order of vanishing moments is not easy. We need to impose an extra condition \er{sym:tf} on the low-pass filter $a$ (or stringent conditions on the underlying refinable function $\phi$ having stable integer shifts, we refer the readers to \cite{ch00,chs02,han13,han14acha,hm05} for more details on this issue), and the moment correcting filter $\Theta$ that we choose must make $\cM_{\pa,\pTh}$ positive semidefinite. This motivates us to consider quasi-tight framelets, which behave almost identical to tight framelets, but have much more flexibility and are much less difficult to construct. As shown in \cite{dh18pp,dhacha,hl19pp,hl20pp}, quasi-tight framelets can be constructed from arbitrary compactly supported refinable (vector) functions, with the desired features (e.g., symmetry, a high order of vanishing moments, and high balancing order) being achieved. This makes quasi-tight framelets of interest in both theory and applications.

\subsection{Illustrative examples}

In this subsection, we provide several examples to illustrate our main results in Theorems~\ref{thm:SymSpecFactSym} and~\ref{thm:qtfsym}. To check whether the function $\phi$ defined in \eqref{def:phi} belongs to $\Lp{2}$ or not, we shall employ the technical quantity $\sm(a)$ (known as the \emph{smoothness exponent of $a$}) defined in \cite[(4.3)]{han03jat} and \cite[(5.6.44)]{hanbook} with $p=2$, which can be easily computed by \cite[Corollary~5.8.5]{hanbook}.
For $\tau\in \R$, a function
$\phi$ belongs to the Sobolev space $H^\tau(\R)$ if $\int_\R |\wh{\phi}(\xi)|^2(1+|\xi|^2)^\tau d\xi<\infty$. Define $\sm(\phi):=\sup\{\tau\in \R \setsp \phi\in H^\tau(\R)\}$.
For a refinable function $\phi$ defined in \eqref{def:phi},
the inequality $\sm(\phi)\ge \sm(a)$ holds and we further have $\sm(\phi)=\sm(a)$ if the integer shifts of $\phi$ are stable, e.g., see \cite[Theorem~6.3.3]{hanbook}.
Hence, $\phi\in \Lp{2}$ if $\sm(a)>0$.
Define $\td\in\lp{0}$ to be the filter such that $\td(0)=1$ and $\td(k)=0$ for all $k\ne 0$.

The first example $\{a;b_1,b_2\}_{\td, (1,-1)}$ of quasi-tight framelet filter banks with symmetry was heuristically obtained in
\cite[Example 3.2.2]{hanbook}, where
\begin{align*}
&\pa(z)=-\frac{1}{16}(z^2+z^{-2})+\frac{1}{4}(z+z^{-1})+\frac{5}{8},\\
&\pb_1(z)=\frac{\sqrt{2}}{4}z(-z+2-z^{-1}),\quad \pb_2(z)=\frac{1}{16}(z^2-4z+6-4z^{-1}+z^{-2}),
\end{align*}
such that
$\sm(a)\approx 0.8853>0$,
$\sym\pa(z)=\sym\pb_2(z)=1$,
$\sym\pb_1(z)=z^2$ and
$\vmo(b_1)=2$, $\vmo(b_2)=4$.
\cite[Example 3.2.2]{hanbook} motivates all later research on quasi-tight framelets in
\cite{dh18pp,dhacha,hl19pp,hl20pp} but all such constructed quasi-tight framelets there lack the symmetry property.
As an application of our main results Theorems~\ref{thm:SymSpecFactSym} and \ref{thm:qtfsym} in this paper, we now provide a few more examples here.

\begin{exmp}\label{ex:TwoHighSymEven1}
Choose $\Theta=\td$ and consider a low-pass filter $a\in \lp{0}$ given by
$$
\pa(z)=-\frac{1}{16}(z^2-6z+1)(1+z)^2z^{-2}+\frac{2\sqrt{2}-3}{32}z^{-3}(1+z)^2(1-z)^4.
$$
Note that $\sym\pa(z)=1$, $\sr(a)=2$ and $\vmo(\pTh(z)-\pTh(z^2)\pa^\star(z)\pa(z))=\vmo(1-\pa^\star(z)\pa(z))=4$. Applying Theorem~\ref{thm:qtfsym} with
$n_b=2$, we can construct a quasi-tight framelet filter bank $\{a;b_1, b_2\}_{\Theta,(1,-1)}$ with
\begin{align*}
\pb_1(z)=&\frac{(1-z)^2}{2048}\left[(4-3\sqrt{2})(z^3+z^{-3})-2\sqrt{2}(z^2+z^{-2}) -(2068-1559\sqrt{2})(z+z^{-1})+1084\sqrt{2}\right],\\ \pb_2(z)=&\frac{(1-z)^2}{2048}\left[(3\sqrt{2}-4)(z^3+z^{-3})+2\sqrt{2}(z^2+z^{-2}) -(2028-1513\sqrt{2})(z+z^{-1})+964\sqrt{2}\right].
\end{align*}
Then $\sym\pb_1(z)=\sym\pb_2(z)=z^2$ and $\vmo(b_1)=\vmo(b_2)=2$.
Define $\phi$ through \eqref{def:phi} and $\psi^\ell$ by $\wh{\psi^\ell}(\xi):=\wh{b_\ell}(\xi/2)\wh{\phi}(\xi/2)$ for all $\xi\in\R$ and $\ell=1,2$. Since $\sm(a)\approx 1.0193>0$, we have
$\phi\in \Lp{2}$. Consequently, we conclude from Theorem~\ref{thm:qtfsym} that $\{\phi;\psi^1,\psi^2\}_{(1,-1)}$ is a quasi-tight framelet in $\Lp{2}$ with symmetry and $\vmo(\psi^1)=\vmo(\psi^2)=2$.
	
	\begin{figure}[h!]
		\centering
		 \begin{subfigure}[]{0.24\textwidth}
			 \includegraphics[width=\textwidth, height=0.8\textwidth]{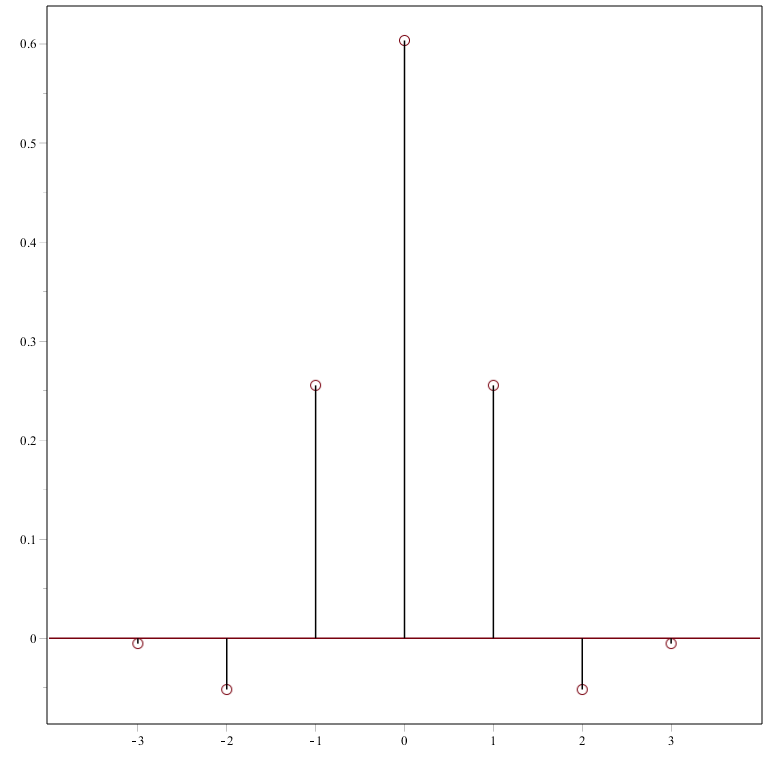}
			\caption{$a$}
		\end{subfigure}
		  \begin{subfigure}[]{0.24\textwidth}
			 \includegraphics[width=\textwidth, height=0.8\textwidth]{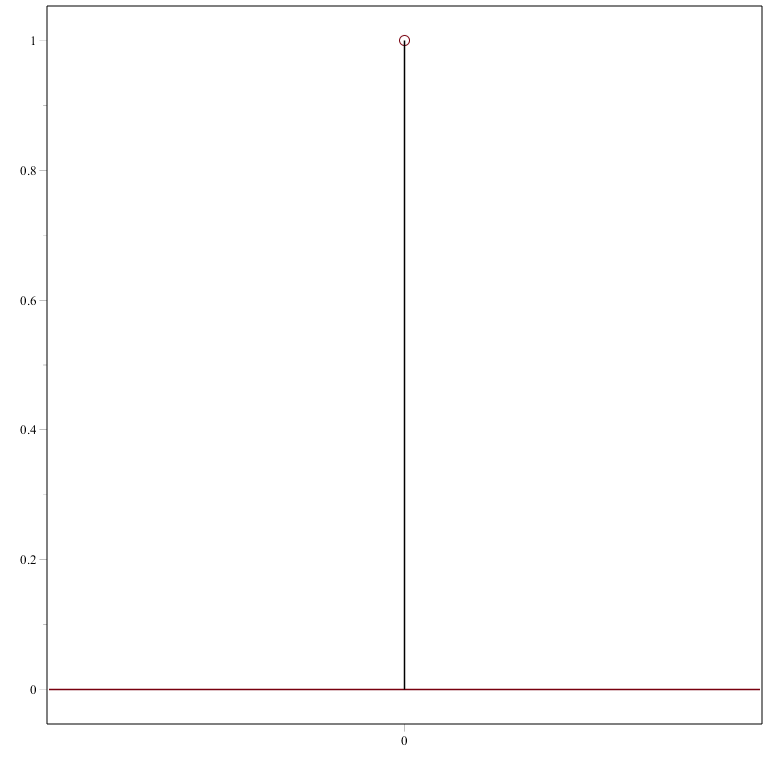}
			\caption{$\td $ }
		\end{subfigure}
\begin{subfigure}[]{0.24\textwidth}
			 \includegraphics[width=\textwidth, height=0.8\textwidth]{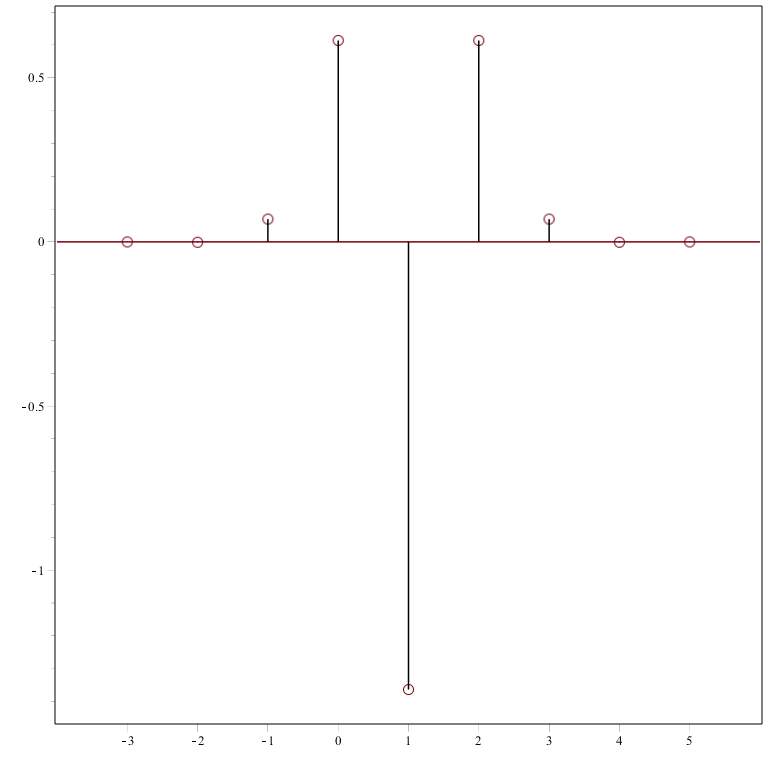}
			\caption{$b_1$}
		\end{subfigure}
		 \begin{subfigure}[]{0.24\textwidth}
			 \includegraphics[width=\textwidth, height=0.8\textwidth]{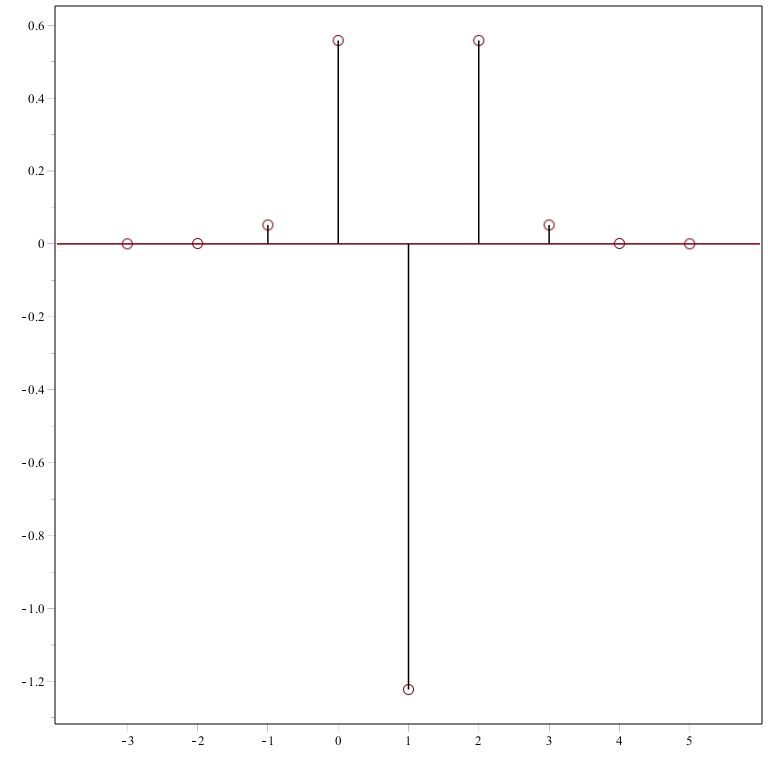}
			\caption{$b_2$}
		\end{subfigure}	
		\\
		 \begin{subfigure}[]{0.24\textwidth}
			 \includegraphics[width=\textwidth, height=0.8\textwidth]{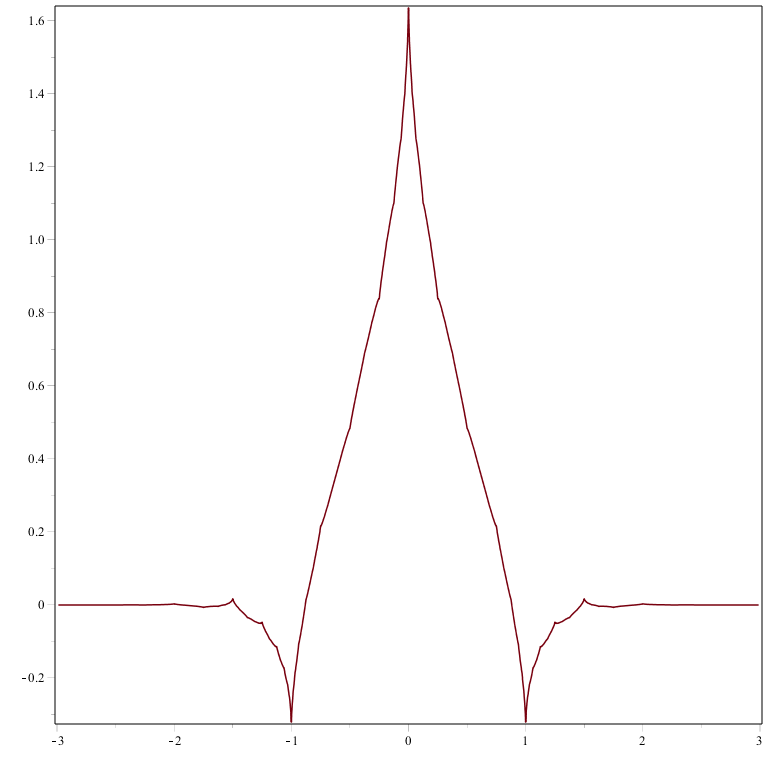}
			\caption{$\phi$}
		\end{subfigure}
		 \begin{subfigure}[]{0.24\textwidth}
			 \includegraphics[width=\textwidth, height=0.8\textwidth]{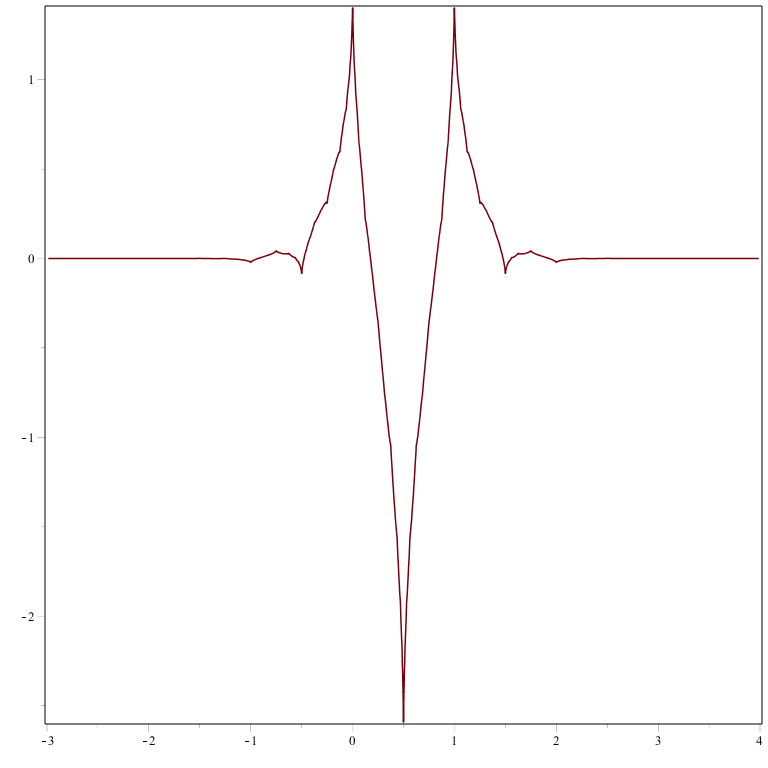}
			\caption{$\psi^1$}
		\end{subfigure}
		 \begin{subfigure}[]{0.24\textwidth}
			 \includegraphics[width=\textwidth, height=0.8\textwidth]{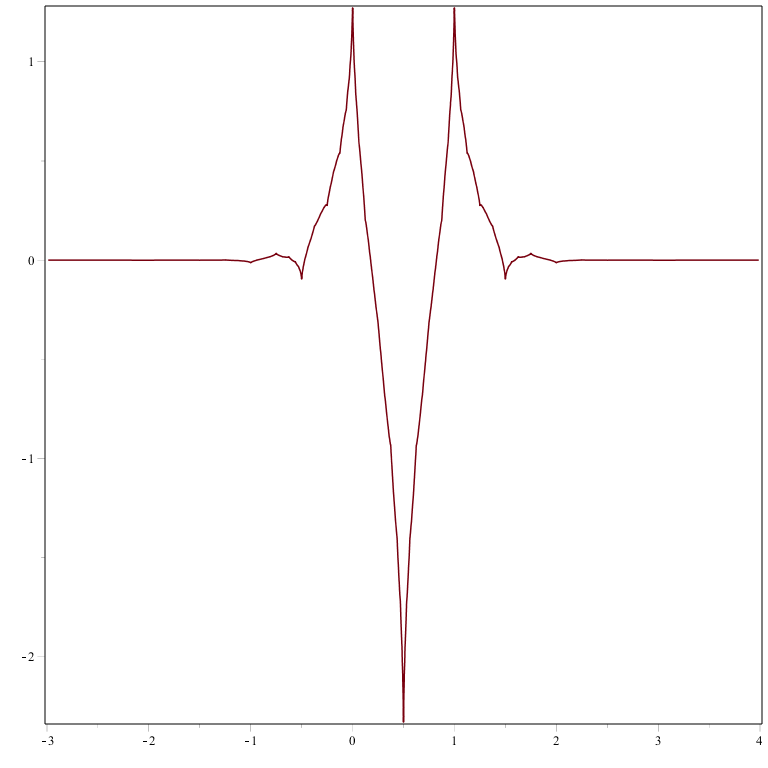}
			\caption{$\psi^2$}
		\end{subfigure}
		 \caption{Example~\ref{ex:TwoHighSymEven1}:
			(A),(B),(C) and (D) are the graphs of the filters $ a, b_1, b_2 $ and $ \Theta $.
			(E), (F) and (G) are the graphs of the refinable function $\phi$ and the framelet generators $\psi^1$ and $\psi^2$.}
	\end{figure}
	
\end{exmp}

\begin{exmp}\label{ex:TwoHighSymEven3}Choose $\Theta=\td$. Let $\phi$ be a compactly supported refinable function with the associated refinement filter $a\in\lp{0}$ such that
$$\pa(z)=\frac{1}{1024}(-z^6+18z^4-32z^3-63z^2+288z+604+288z^{-1}-63z^{-2}-32z^{-3}+18z^{-4}-z^{-6}).$$
	We have
$\sym\pa(z)=1$, $\sr(a)=4$, and $\vmo(1-\pa^\star\pa)=8$. Since $\sm(a)\approx 1.6821>0$, $\phi$ in \eqref{def:phi} belongs to $\Lp{2}$.
For $n_b=4$, we can construct a quasi-tight framelet filter bank $\{a;b_1,b_2\}_{\Theta, (1,-1)}$ with
	 $$\pb_1(z)=\frac{\sqrt{2}}{32}(-z^{4}+9z^2-16z+9-z^{-2}),$$
	 $$\pb_2(z)=\frac{1}{1024}(z^6-18z^4+32z^3+63z^2-288z+420-288z^{-1}+63z^{-2}+32z^{-3}-18z^{-4}+z^{-6}).$$
	We have $\sym\pb_1(z)=z^2$, $\sym\pb_2(z)=1$ and $\vmo(b_1)=4$, $\vmo(b_2)=8$. By letting $\wh{\psi^\ell}(\xi):=\wh{b_\ell}(\xi/2)\wh{\phi}(\xi/2)$ for all $\xi\in\R$ and $\ell=1,2$, we see that $\{\phi;\psi^1,\psi^2\}_{(1,-1)}$ is a quasi-tight framelet in $\Lp{2}$ with symmetry.
	
	\begin{figure}[h!]
		\centering
		 \begin{subfigure}[]{0.24\textwidth}
			 \includegraphics[width=\textwidth, height=0.8\textwidth]{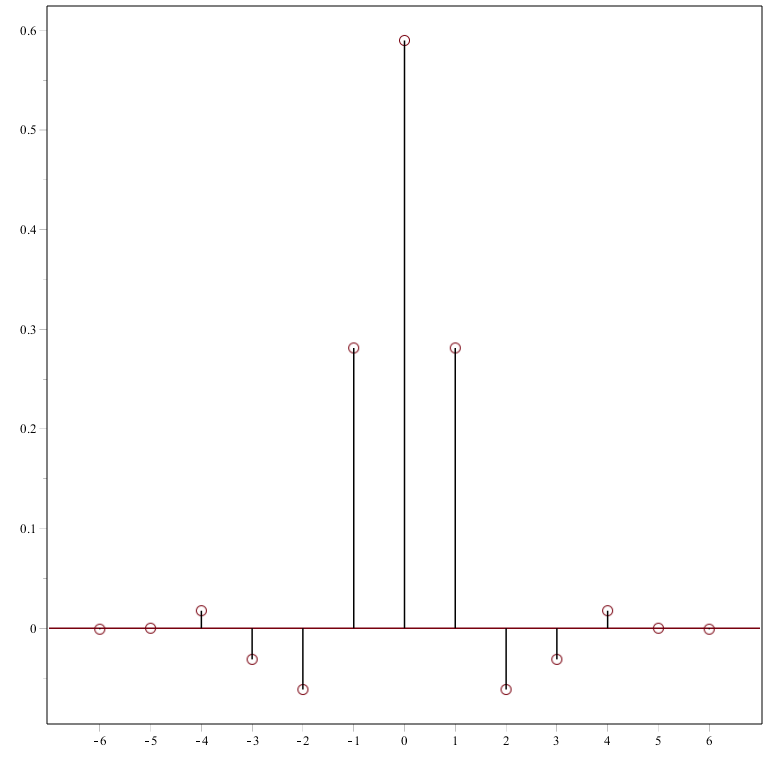}
			\caption{$a$}
		\end{subfigure}	
	 \begin{subfigure}[]{0.24\textwidth}
			 \includegraphics[width=\textwidth, height=0.8\textwidth]{StemDelta}
			\caption{$\td $ }
		\end{subfigure}
		 \begin{subfigure}[]{0.24\textwidth}
			 \includegraphics[width=\textwidth, height=0.8\textwidth]{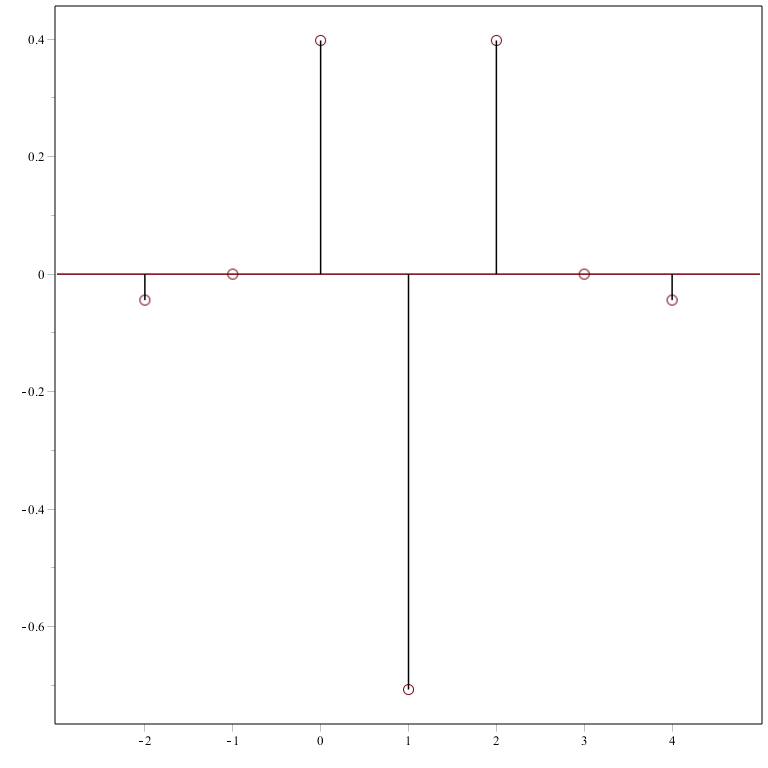}
			\caption{$b_1$}
		\end{subfigure}
		 \begin{subfigure}[]{0.24\textwidth}
			 \includegraphics[width=\textwidth, height=0.8\textwidth]{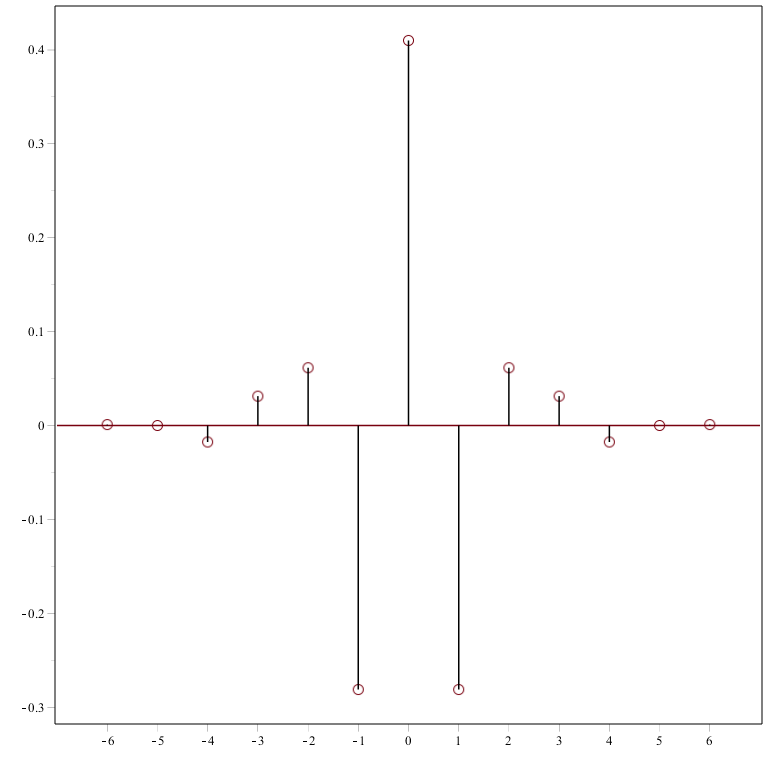}
			\caption{$b_2$}
		\end{subfigure}
		\\
		 \begin{subfigure}[]{0.24\textwidth}
			 \includegraphics[width=\textwidth, height=0.8\textwidth]{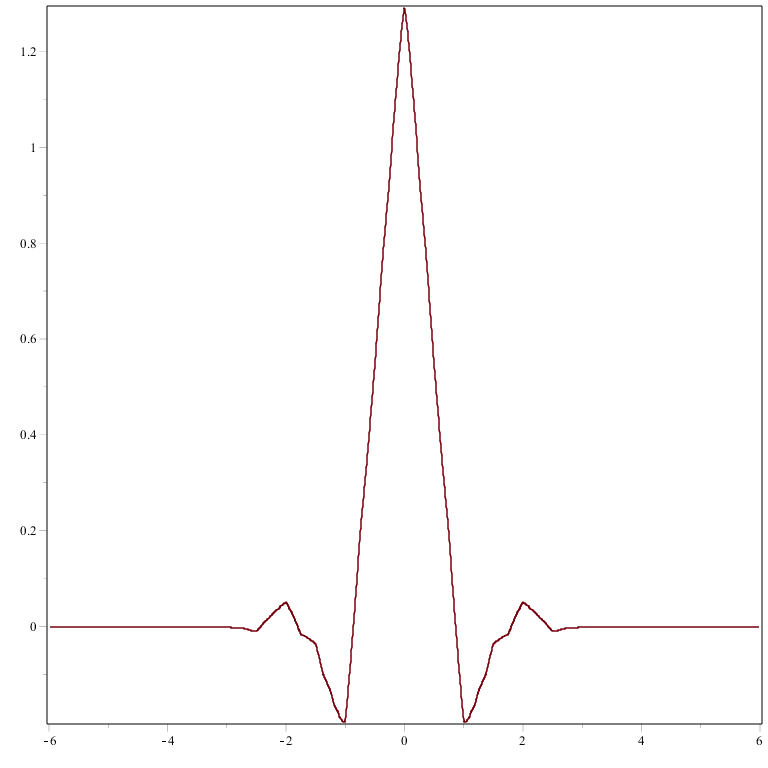}
			\caption{$\phi$}
		\end{subfigure}
		 \begin{subfigure}[]{0.24\textwidth}
			 \includegraphics[width=\textwidth, height=0.8\textwidth]{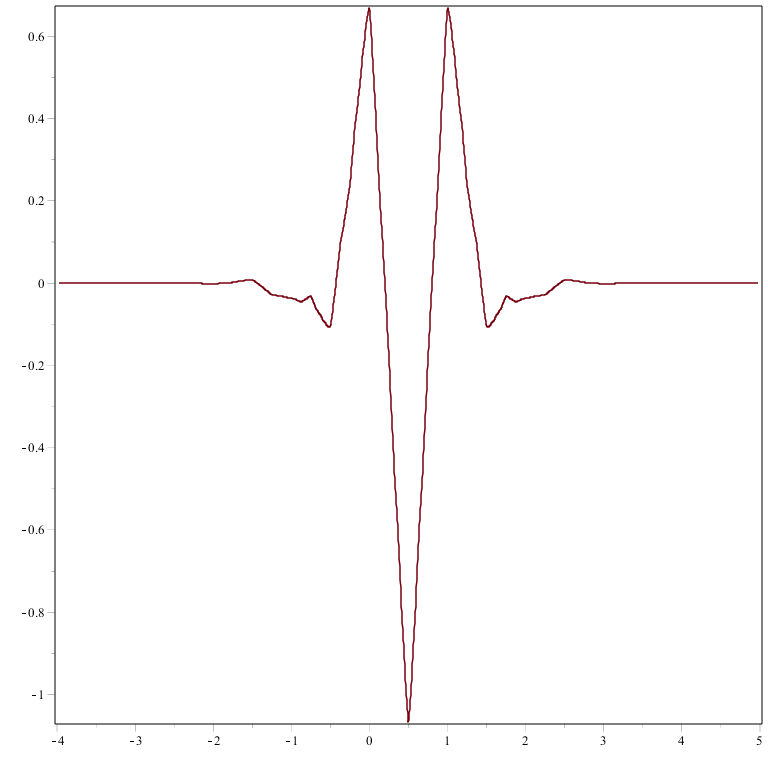}
			\caption{$\psi^1$}
		\end{subfigure}
		 \begin{subfigure}[]{0.24\textwidth}
			 \includegraphics[width=\textwidth, height=0.8\textwidth]{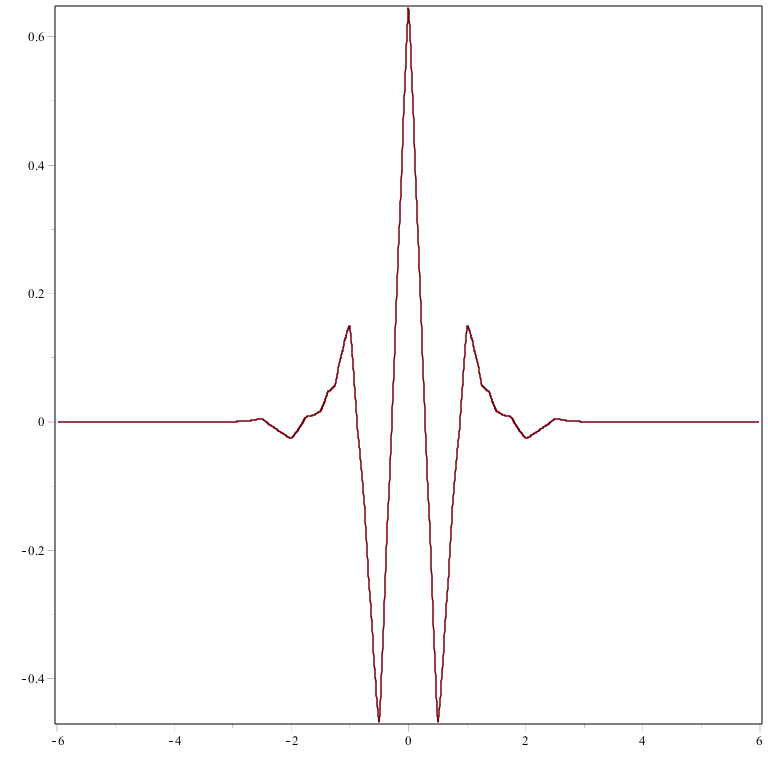}
			\caption{$\psi^2$}
		\end{subfigure}
		\caption{In Example~\ref{ex:TwoHighSymEven3}:
			(A),(B),(C) and (D) are the graphs of the filters $a, \Theta, b_1, b_2 $.
			(E), (F) and (G) are the graphs of the refinable function $\phi$ and the framelet generators $\psi^1$ and $\psi^2$. }
	\end{figure}
	
\end{exmp}

\begin{exmp}\label{ex:TwoHighSymOdd1}Choose $\Theta=\td$. Let $\phi$ be a compactly supported refinable function with the associated refinement filter $a\in\lp{0}$ such that
	 $$\pa(z)=\frac{1}{1024}(15z^4-63z^3+35z^2+525z+525+35z^{-1}-63z^{-2}+15z^{-3}).$$
	We have $\sym\pa(z)=z$, $\sr(a)=3$, and $\vmo(1-\pa^\star\pa)=6$. Since $\sm(a)\approx 1.1543>0$, $\phi$ in \eqref{def:phi} belongs to $\Lp{2}$. For $n_b=3$, we can construct a quasi-tight framelet filter bank $\{a;b_1,b_2\}_{\Theta,(1,-1)}$ with
	 $$\pb_1(z)=\frac{1}{1024}(15z^4-63z^3+385z^2-945z+945-385z^{-1}+63z^{-2}-15z^{-3}),$$
	 $$\pb_2(z)=\frac{\sqrt{105}}{512}(-5z^4+21z^3-38z^2+38z-21+5z^{-1}).$$
	We have $\sym\pb_1(z)=-z$, $\sym\pb_2(z)=-z^3$ and $\vmo(b_1)=\vmo(b_2)=3$. By letting $\wh{\psi^\ell}(\xi):=\wh{b_\ell}(\xi/2)\wh{\phi}(\xi/2)$ for all $\xi\in\R$ and $\ell=1,2$, we see that $\{\phi;\psi^1,\psi^2\}_{(1,-1)}$ is a quasi-tight framelet in $\Lp{2}$.
	
	\begin{figure}[h!]
		\centering
		 \begin{subfigure}[]{0.24\textwidth}
			 \includegraphics[width=\textwidth, height=0.8\textwidth]{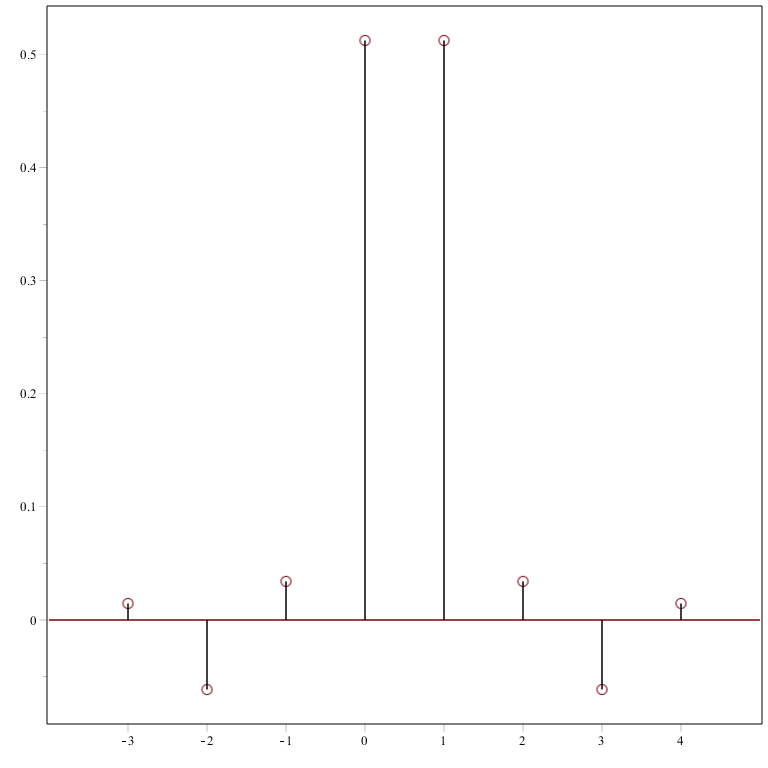}
			\caption{$a$}
		\end{subfigure}	
 \begin{subfigure}[]{0.24\textwidth}
			 \includegraphics[width=\textwidth, height=0.8\textwidth]{StemDelta}
			\caption{$\td $ }
		\end{subfigure}
		 \begin{subfigure}[]{0.24\textwidth}
			 \includegraphics[width=\textwidth, height=0.8\textwidth]{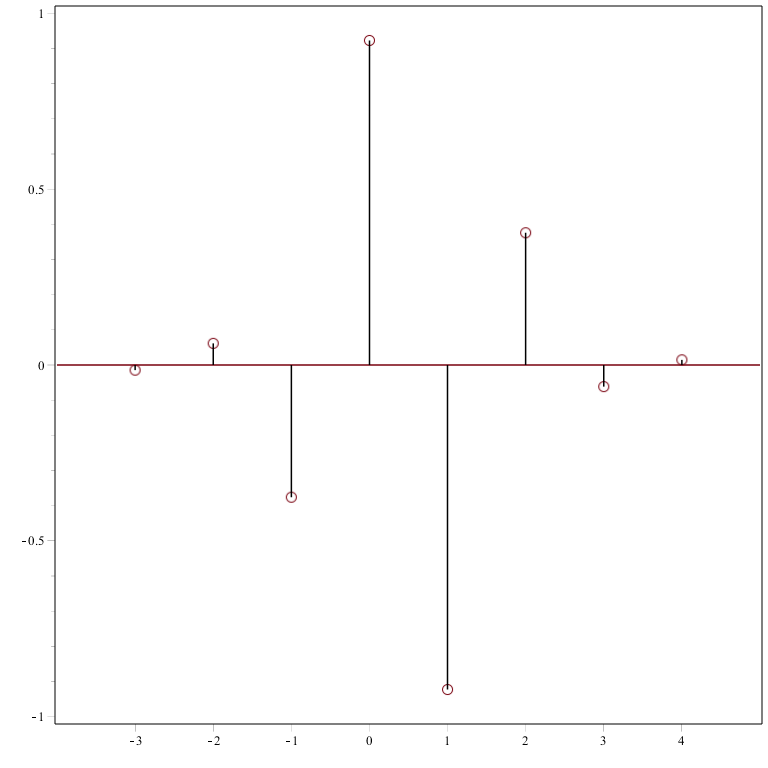}
			\caption{$b_1$}
		\end{subfigure}
		 \begin{subfigure}[]{0.24\textwidth}
			 \includegraphics[width=\textwidth, height=0.8\textwidth]{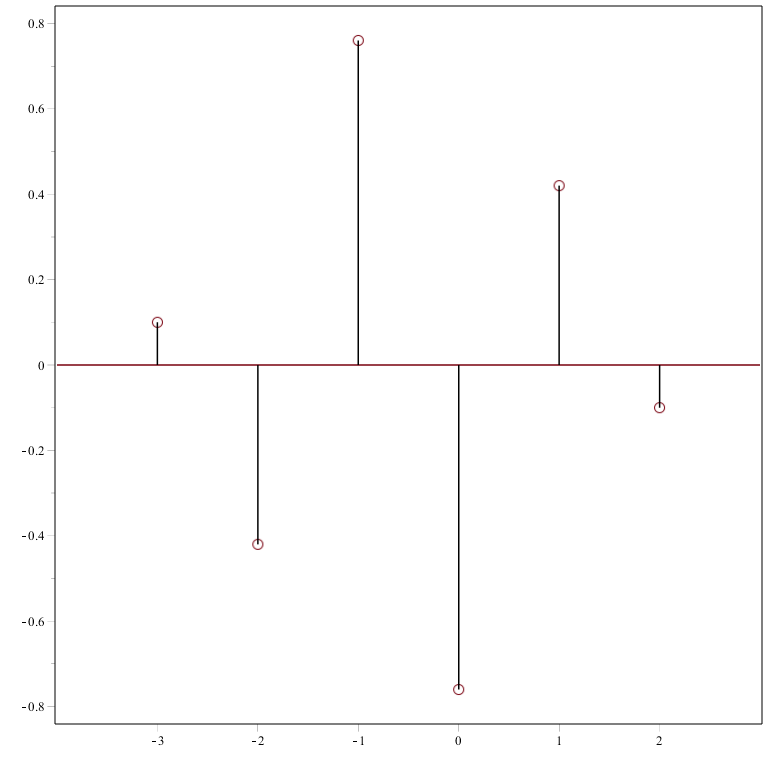}
			\caption{$b_2$}
		\end{subfigure}
		\\
		 \begin{subfigure}[]{0.24\textwidth}
			 \includegraphics[width=\textwidth, height=0.8\textwidth]{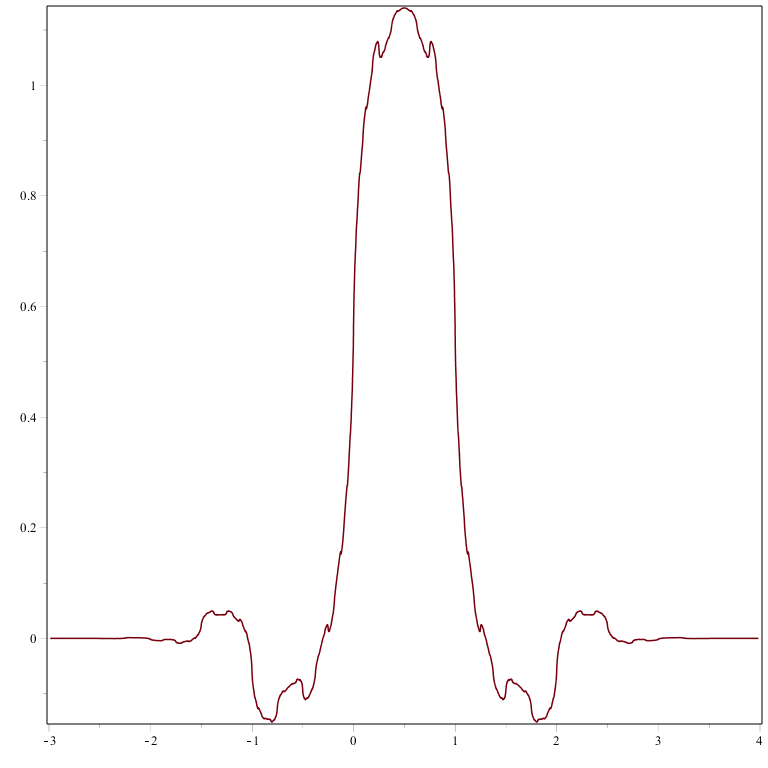}
			\caption{$\phi$}
		\end{subfigure}
		 \begin{subfigure}[]{0.24\textwidth}
			 \includegraphics[width=\textwidth, height=0.8\textwidth]{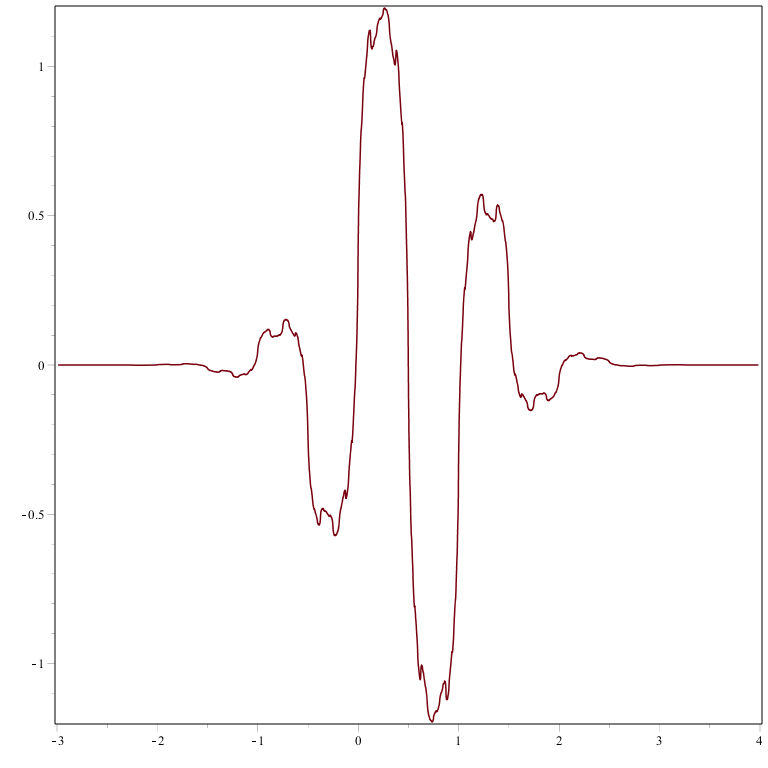}
			\caption{$\psi^1$}
		\end{subfigure}
		 \begin{subfigure}[]{0.24\textwidth}
			 \includegraphics[width=\textwidth, height=0.8\textwidth]{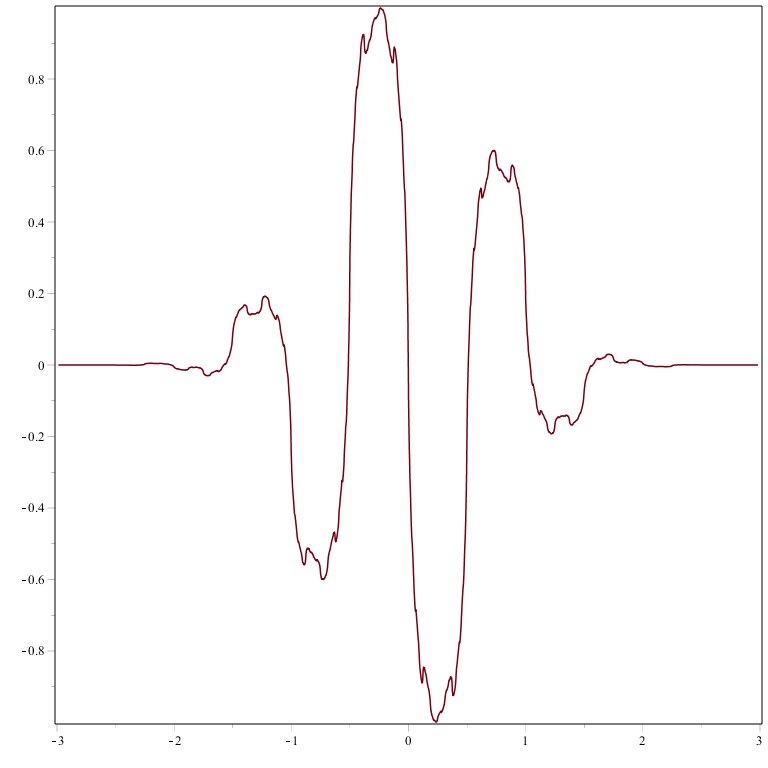}
			\caption{$\psi^2$}
		\end{subfigure}
		\caption{In Example~\ref{ex:TwoHighSymOdd1}:
			(A),(B),(C) and (D) are the graphs of the filters $a, \Theta, b_1, b_2 $.
			(E), (F) and (G) are the graphs of the refinable function $\phi$ and the framelet generators $\psi^1$ and $\psi^2$. }
	\end{figure}
\end{exmp}

\begin{exmp}\label{ex:TwoHighSymOdd2}Choose $\Theta=\td$. Let $\phi$ be a compactly supported refinable function with the associated refinement filter $a\in\lp{0}$ such that
	 $$\pa(z)=\frac{1}{32}(-z^3+z^2+16z+16+z^{-1}-z^{-2}).$$
	We have $\sym\pa(z)=z$, $\sr(a)=1$, and $\vmo(1-\pa^\star\pa)=4$. Since $\sm(a)\approx 0.7184>0$, $\phi$ in \eqref{def:phi} belongs to $\Lp{2}$. For $n_b=1$, we can construct a quasi-tight framelet filter bank $\{a;b_1,b_2\}_{\Theta,(1,-1)}$ with
	 $$\pb_1(z)=\frac{\sqrt{2}}{4096}(z^5-z^4-32z^3-32z^2+2032z-2032+32z^{-1}+32z^{-2}+z^{-3}-z^{-4}),$$
	 $$\pb_2(z)=\frac{\sqrt{2}}{4096}(-z^5+z^4+32z^3+32z^2+1794z-1794-32z^{-1}-32z^{-2}-z^{-3}+z^{-4}).$$
	We have $\sym\pb_1(z)=\sym\pb_2(z)=-z$ and $\vmo(b_1)=\vmo(b_2)=1$. By letting $\wh{\psi^\ell}(\xi):=\wh{b_\ell}(\xi/2)\wh{\phi}(\xi/2)$ for all $\xi\in\R$ and $\ell=1,2$, we see that $\{\phi;\psi^1,\psi^2\}_{(1,-1)}$ is a quasi-tight framelet in $\Lp{2}$ with symmetry.

	\begin{figure}[h]
		\centering
		 \begin{subfigure}[]{0.24\textwidth}
			 \includegraphics[width=\textwidth, height=0.8\textwidth]{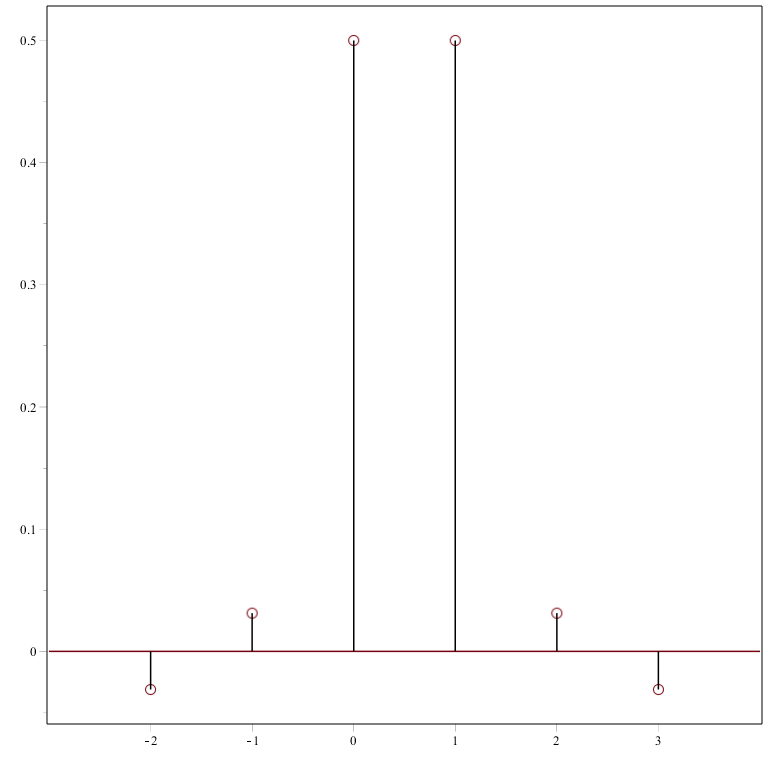}
			\caption{$a$}
		\end{subfigure}	
\begin{subfigure}[]{0.24\textwidth}
			 \includegraphics[width=\textwidth, height=0.8\textwidth]{StemDelta}
			\caption{$\td $ }
		\end{subfigure}
		 \begin{subfigure}[]{0.24\textwidth}
			 \includegraphics[width=\textwidth, height=0.8\textwidth]{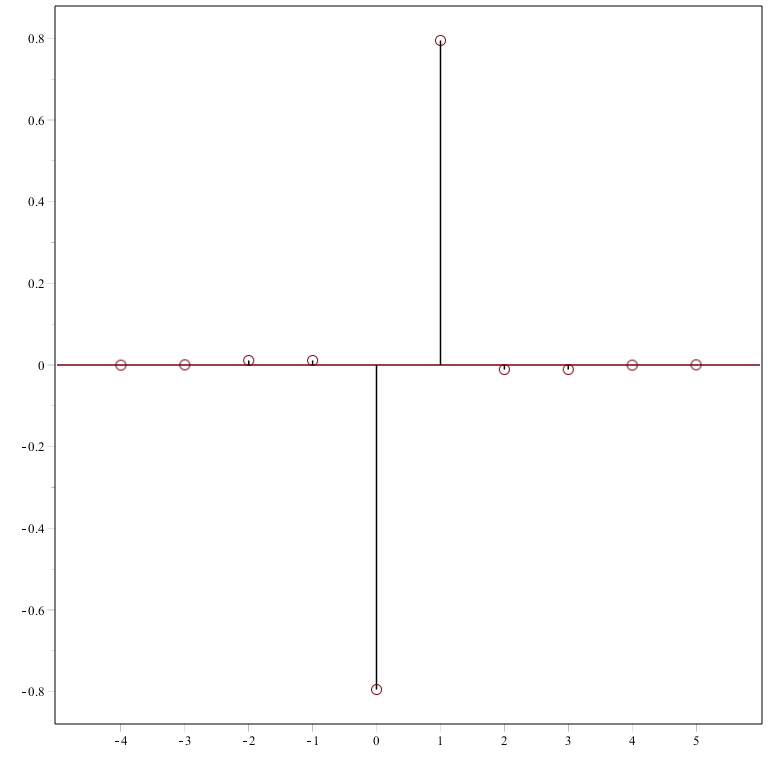}
			\caption{$b_1$}
		\end{subfigure}
		 \begin{subfigure}[]{0.24\textwidth}
			 \includegraphics[width=\textwidth, height=0.8\textwidth]{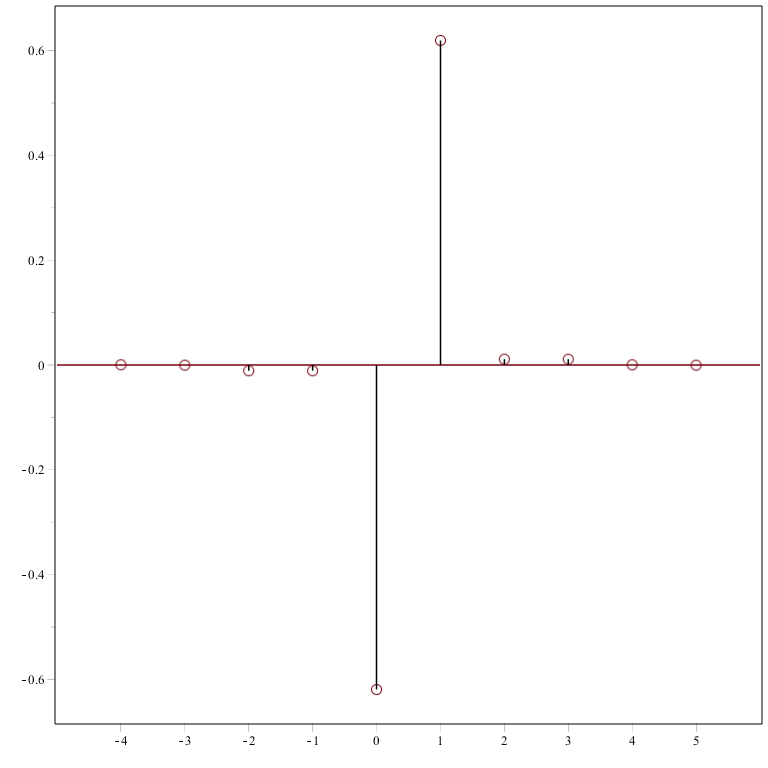}
			\caption{$b_2$}
		\end{subfigure}
		\\
		 \begin{subfigure}[]{0.24\textwidth}
			 \includegraphics[width=\textwidth, height=0.8\textwidth]{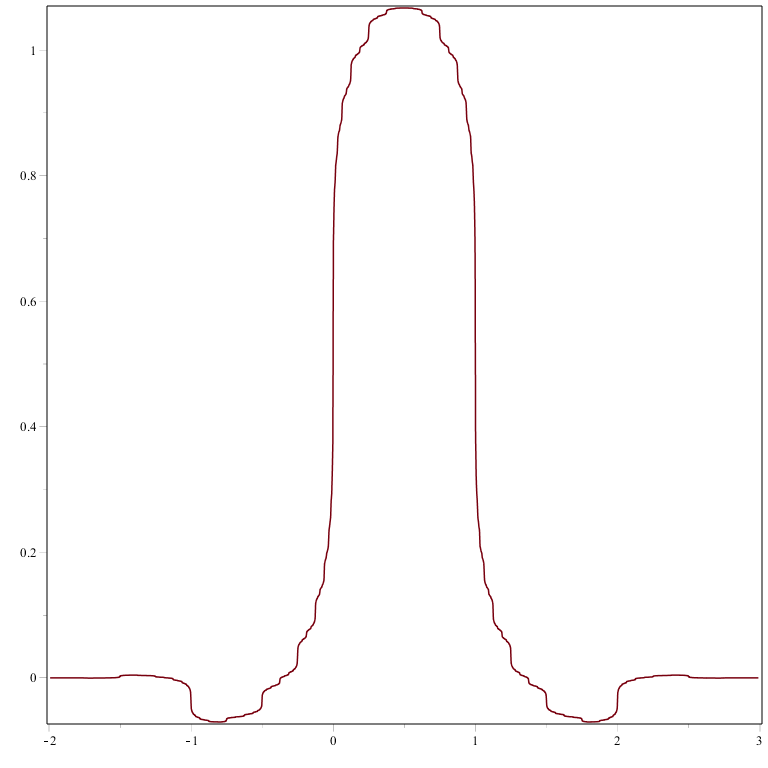}
			\caption{$\phi$}
		\end{subfigure}
		 \begin{subfigure}[]{0.24\textwidth}
			 \includegraphics[width=\textwidth, height=0.8\textwidth]{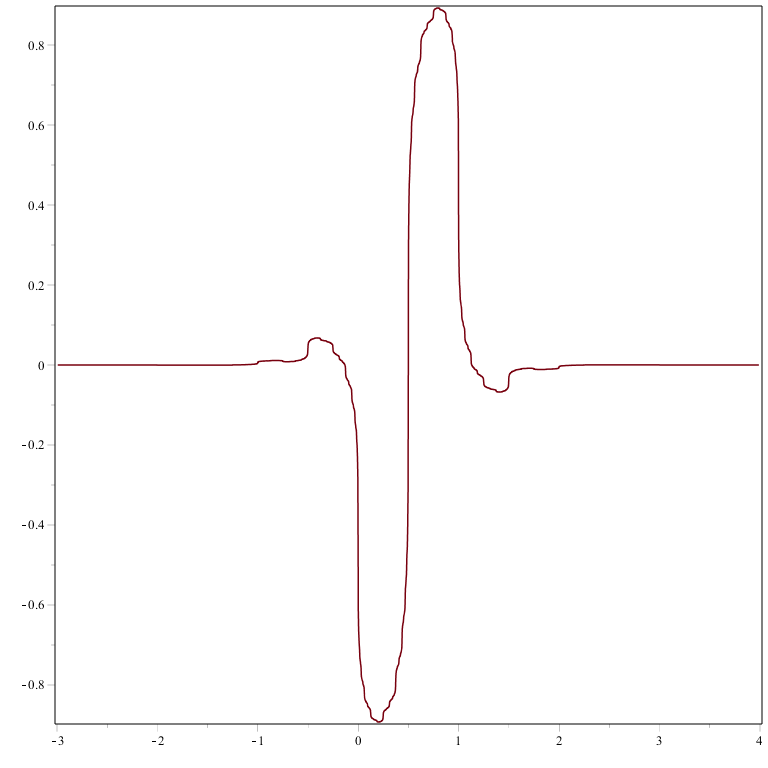}
			\caption{$\psi^1$}
		\end{subfigure}
		 \begin{subfigure}[]{0.24\textwidth}
			 \includegraphics[width=\textwidth, height=0.8\textwidth]{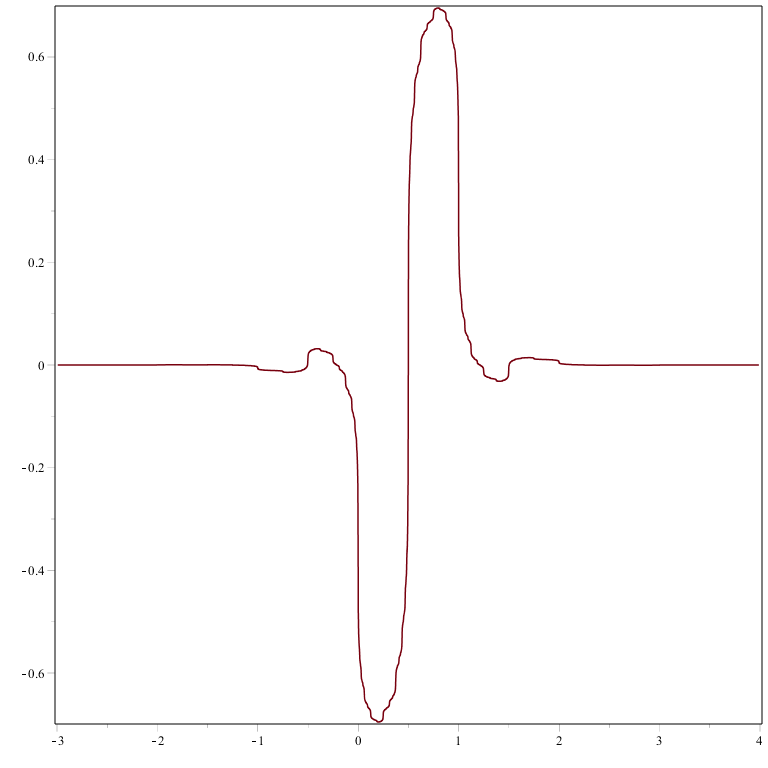}
			\caption{$\psi^2$}
		\end{subfigure}
		\caption{In Example~\ref{ex:TwoHighSymOdd2}:
			(A),(B),(C) and (D) are the graphs of the filters $a, \Theta, b_1, b_2 $.
			(E), (F) and (G) are the graphs of the refinable function $\phi$ and the framelet generators $\psi^1$ and $\psi^2$.}
	\end{figure}	
	
\end{exmp}

\begin{exmp}\label{ex:TwoHighSym6}Choose $\pTh(z)=\tfrac{z+z^{-1}}{2}$. Let $\phi$ be a compactly supported refinable function with the associated refinement filter $a\in\lp{0}$ such that
	 $$\pa(z)=\frac{1}{8}(-z^3+3z+4+3z^{-1}-z^{-3}).$$
	We have $\sym\pa(z)=1$, $\sym\pTh(z)=1$, $\sr(a)=2$, and $\vmo(\pTh-\pTh(\cdot^2)\pa^\star\pa)=4$. Since $\sm(a)\approx 1$, $\phi$ in \eqref{def:phi} belongs to $\Lp{2}$. For $n_b=2$, we can construct a quasi-tight framelet filter bank $\{a;b_1,b_2\}_{\Theta,(1,-1)}$ with
	 $$\pb_1(z)=\frac{1}{64}z^{-1}(z-1)^2(z^4+2z^3-12z^2-30z-46-30z^{-1}-12z^{-2}+2z^{-3}+z^{-4}),$$
	 $$\pb_2(z)=\frac{1}{64}z^{-1}(z-1)^2(z^2-2z+4-2z^{-1}+1)(z^2+4z+8+4z^{-1}+1).$$
	We have $\sym\pb_1(z)=\sym\pb_2(z)=1$ and $\vmo(b_1)=\vmo(b_2)=2$. By letting $\wh{\eta}(\xi):=\wh{\Theta}(\xi)\wh{\phi}(\xi)$ and  $\wh{\psi^\ell}(\xi):=\wh{b_\ell}(\xi/2)\wh{\phi}(\xi/2)$ for all $\xi\in\R$ and $\ell=1,2$, we see that $\{\eta,\phi;\psi^1,\psi^2\}_{(1,-1)}$ is a quasi-tight framelet and thus $\{\psi^1,\psi^2\}_{(1,-1)}$ is a homogeneous quasi-tight framelet in $\Lp{2}$ with symmetry.

	\begin{figure}[h]
	\centering
	\begin{subfigure}[]{0.24\textwidth}
		\includegraphics[width=\textwidth, height=0.8\textwidth]{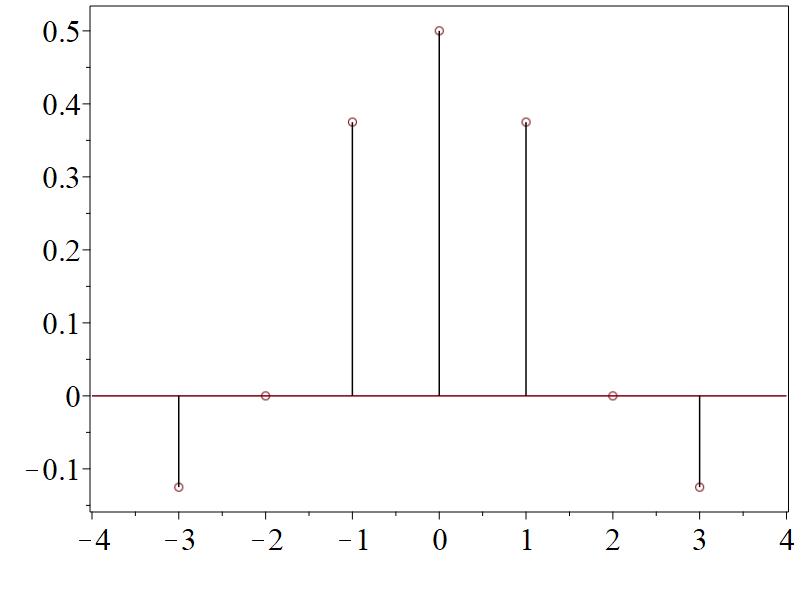}
		\caption{$a$}
	\end{subfigure}
		\begin{subfigure}[]{0.24\textwidth}
		\includegraphics[width=\textwidth, height=0.8\textwidth]{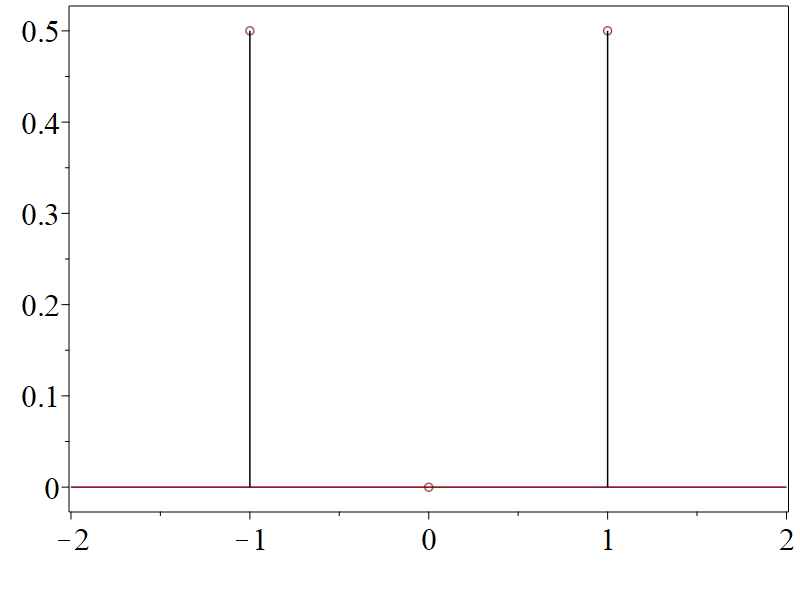}
		\caption{$\Theta $ }
	\end{subfigure}
\begin{subfigure}[]{0.24\textwidth}
		\includegraphics[width=\textwidth, height=0.8\textwidth]{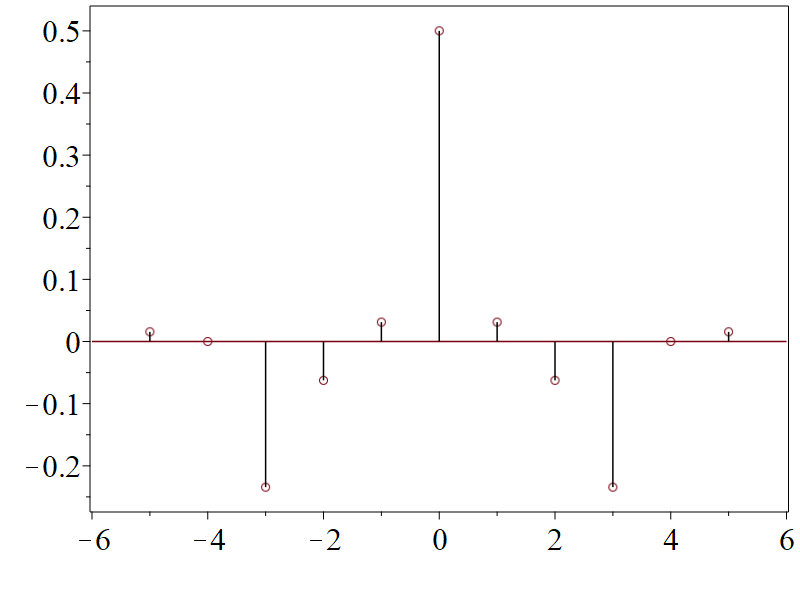}
		\caption{$b_1$}
	\end{subfigure}
	\begin{subfigure}[]{0.24\textwidth}
		\includegraphics[width=\textwidth, height=0.8\textwidth]{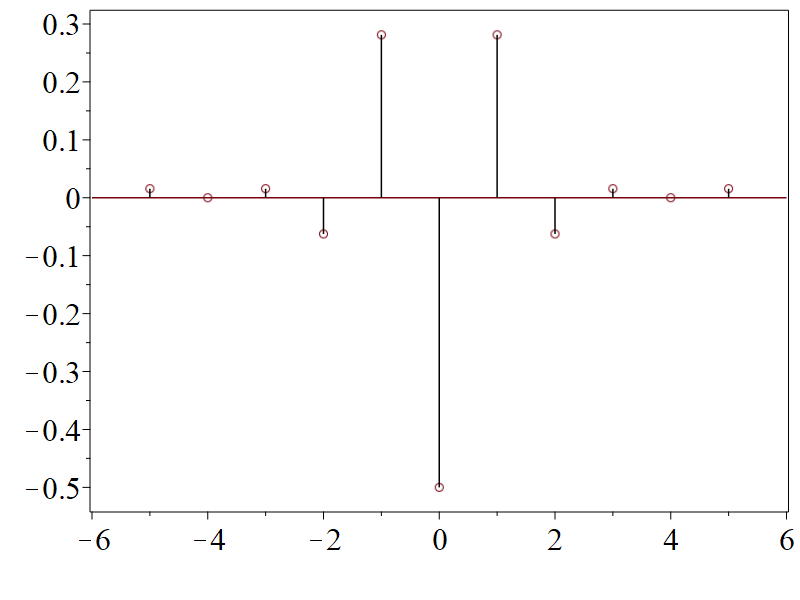}
		\caption{$b_2$}
	\end{subfigure}
\\
	\begin{subfigure}[]{0.24\textwidth}
		\includegraphics[width=\textwidth, height=0.8\textwidth]{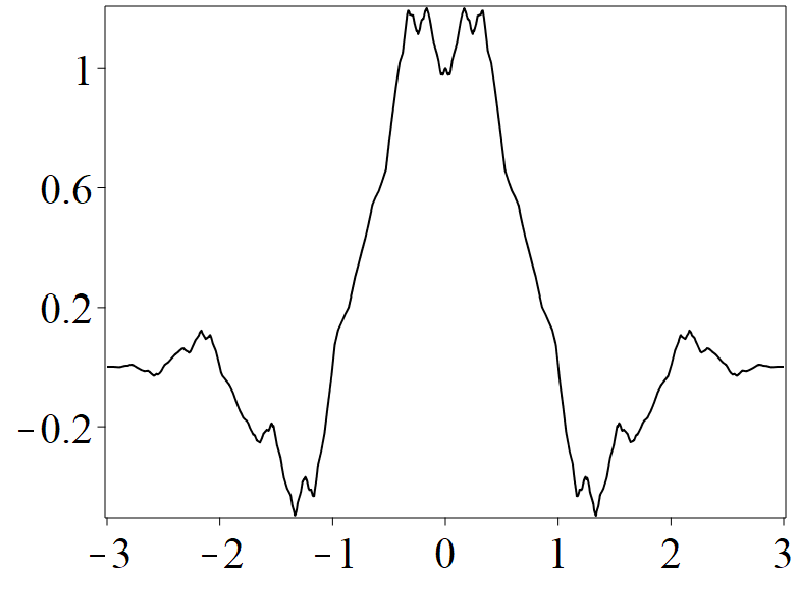}
		\caption{$\phi$}
	\end{subfigure}
\begin{subfigure}[]{0.24\textwidth}
		\includegraphics[width=\textwidth, height=0.8\textwidth]{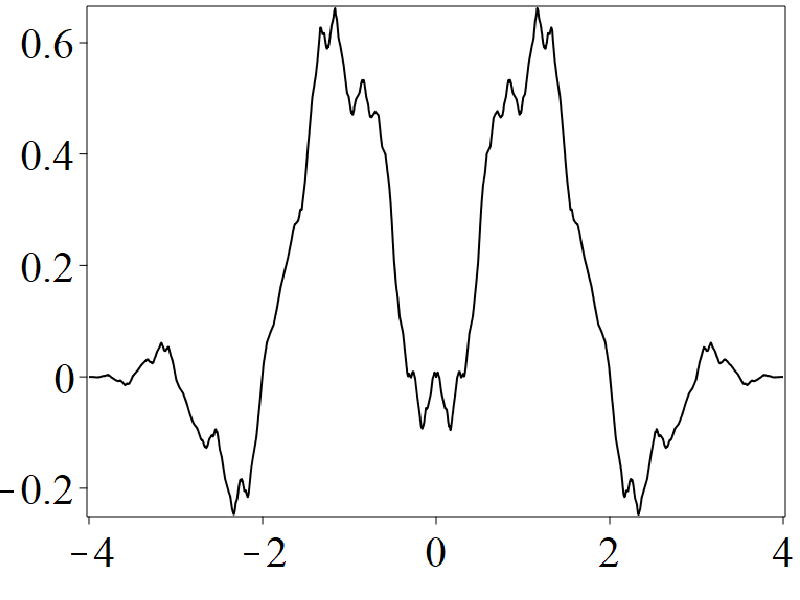}
		\caption{$\eta$}
	\end{subfigure}
	\begin{subfigure}[]{0.24\textwidth}
		\includegraphics[width=\textwidth, height=0.8\textwidth]{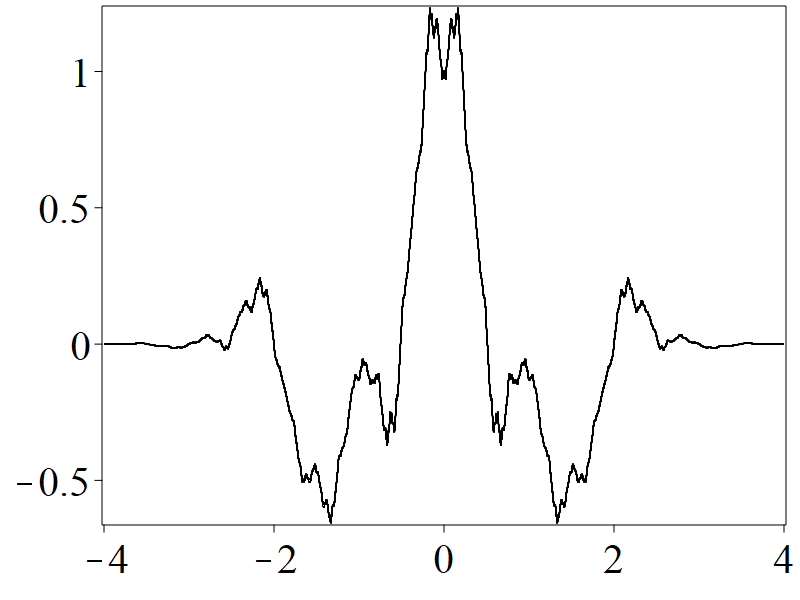}
		\caption{$\psi^1$}
	\end{subfigure}
	\begin{subfigure}[]{0.24\textwidth}
		\includegraphics[width=\textwidth, height=0.8\textwidth]{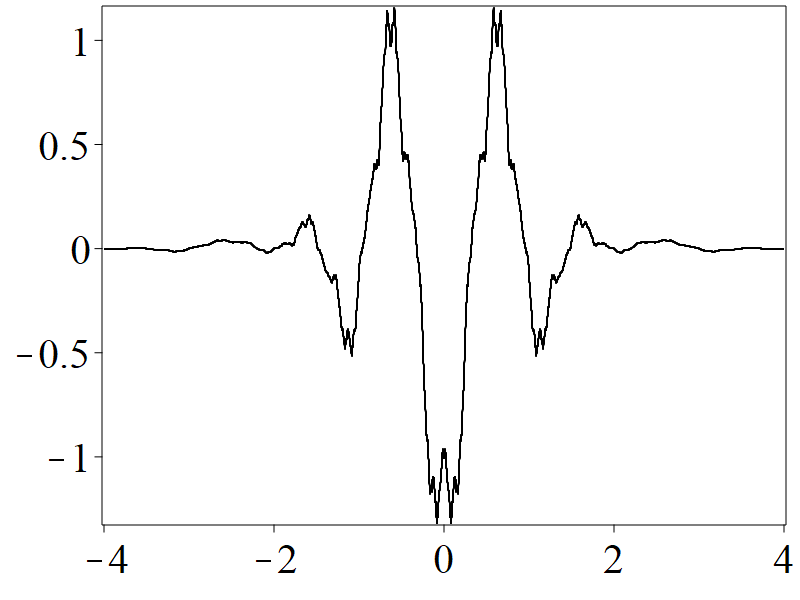}
		\caption{$\psi^2$}
	\end{subfigure}
	\caption{In Example~\ref{ex:TwoHighSym6}:
		(A),(B),(C) and (D) are the graphs of the filters $a, \Theta, b_1, b_2 $.
		(E), (F), (G) and (H) are the graphs of the refinable functions $\phi$, $\eta$ and the framelet generators $\psi^1$ and $\psi^2$.}
\end{figure}

\end{exmp}

\section{Factorization of Matrices of Laurent Polynomials with Symmetry}\label{sec:fac:lau}

In this section, we study the factorization problem proposed at the beginning of the paper. Our goal of this section is to prove the main result Theorem~\ref{thm:SymSpecFactSym} on the existence of a generalized spectral factorization as in \er{eq:spectfact:0} with $\eps_1=1$ and $\eps_2=-1$.

\subsection{Laurent polynomials with symmetry}
To prove Theorem~\ref{thm:SymSpecFactSym}, we need to better understand some properties of Laurent polynomials with symmetry.
We first introduce several notations which will be used throughout this paper.
For a nonzero Laurent polynomial $\pu(z)=\sum_{k\in\Z}u(k)z^k$, define its \emph{lower degree} $\ldeg(\pu)$, \emph{degree} $\deg(\pu)$ and its \emph{length} $\len(\pu)$ by
$$
\ldeg(\pu):=\min\{k\in\Z:u(k)\ne 0\},\quad \deg(\pu):=\max\{k\in\Z:u(k)\ne 0\},
\quad \len(\pu):=\deg(\pu)-\ldeg(\pu).
$$
If $\pu=\mathbf{0}$, then we define $\len(\mathbf{0}):=-\infty$.
Here are some basic facts about the symmetry operator.

\begin{prop}\label{prop:sym:lau}Suppose $\pu$ and $\pv$ are Laurent polynomials with symmetry types $\eps_uz^{c_u}$ and $\eps_vz^{c_v}$ respectively. Then the following hold:

\begin{enumerate}
\item[(1)] $\pu\pv$ has symmetry with $\sym[\pu\pv](z)=\eps_u\eps_vz^{c_u+c_v}$;

\item[(2)] If $\pv$ divides $\pu$, then $\pu/\pv$ has symmetry with $\sym[\pu/\pv](z)=\eps_u\eps_vz^{c_u-c_v}$;

\item[(3)] $\pu^\star$ has symmetry with $\sym[\pu^\star](z)=[\sym\pu]^\star(z)=\eps_uz^{-c_u}.$

\end{enumerate}
If $\pu\ne 0$, then $\sym[\pu^\star]=\sym\pu^{-1}$,
$c_u=\ldeg(\pu)+\deg(\pu)$,
and $\odd(\len(\pu))=\odd(c_u)$, where $\odd$ is defined as in \er{odd}.
\end{prop}

\bp Items (1) -- (3) are obvious. Note that $\ldeg(\pu(\cdot^{-1}))=-\deg(\pu)$. By the definition of $\sym$, we have $\ldeg(\pu)=c_u+\ldeg(\pu(\cdot^{-1}))=c_u-\deg(\pu)$ and $\len(\pu)=\deg(\pu)-\ldeg(\pu)=c_u-2\ldeg(\pu)$.
Now all the rest claims follow easily.
\ep

The symmetry of a Laurent polynomial is related to its multiplicities of roots at $\pm 1$. Recall that $\mz(\pu,z_0)$ stands for the multiplicity of the root of $\pu$ at $z_0$. We have the following lemma.

\begin{lemma}\label{lem:rt:pm1}Let $\pu$ be a nonzero Laurent polynomial with symmetry type $\sym\pu(z)=\eps z^c$ for some $\eps\in\{\pm 1\}$ and $c\in\Z$. Then $\eps=(-1)^{\mz(\pu,1)}$ and $\odd(c)=\odd(\mz(\pu, 1)+\mz(\pu,-1))$.
\end{lemma}

\bp Write $\pu(z)=(z-1)^{\mz(\pu,1)}\pv(z)$ for some Laurent polynomial $\pv$ with $\pv(1)\neq 0$. We have $\eps z^c=\frac{\pu(z)}{\pu(z^{-1})}=(-z)^{\mz(\pu,1)}\frac{\pv(z)}{\pv(z^{-1})}.$ By letting $z=1$, we obtain $\eps=(-1)^{\mz(\pu,1)}$.
We can also write $\pu(z)=(z+1)^{\mz(\pu,-1)} \pw(z)$ for a unique Laurent polynomial $\pw$ with $\pw(-1)\ne 0$. By $\eps=(-1)^{\mz(\pu,1)}$ and the definition of $\sym$ in \eqref{eq:RealSymZ}, we have
$$
(z+1)^{\mz(\pu,-1)}\pw(z)=\pu(z)=\eps z^c u(z^{-1})=
(-1)^{\mz(\pu,1)}z^{c-\mz(\pu,-1)}(z+1)^{\mz(\pu,-1)} \pw(z^{-1}),
$$
from which we have $\pw(z)=(-1)^{\mz(\pu,1)}z^{c-\mz(\pu,-1)} \pw(z^{-1})$. Taking $z=-1$ and noting $\pw(-1)\ne 0$, we conclude that $(-1)^{c-\mz(\pu,-1)+\mz(\pu,1)}=1$, i.e., $\odd(c)=\odd(\mz(\pu, 1)+\mz(\pu,-1))$.\ep

\subsection{Compatibility of matrices of Laurent polynomials with symmetry}
We now consider matrices of Laurent polynomials with symmetry. First, we extend the symmetry operator $\sym$. For an $r\times s$ matrix $\pA:=\pA(z)$ of Laurent polynomials with symmetry, define $\sym\pA$ via $[\sym\pA (z)]_{j,k}:=\sym(\pA_{j,k}(z))$ for all $1\le j\le r$ and $1\le k\le s$. In this case, we say that the matrix $\pA$ \emph{has symmetry type $\sym\pA$}.

Next, we discuss how the symmetry property behaves under matrix summations. If $\pA$ and $\pB$ are $r\times s$ matrices of Laurent polynomials with symmetry and $\sym\pA=\sym\pB$, then it is not difficult to conclude that $\pA\pm\pB$ have symmetry with
$$\sym[\pA+\pB]=\sym[\pA-\pB]=\sym\pA=\sym\pB.$$
In this case, we say that \emph{the operations $\pA+\pB$ and $\pA-\pB$ are compatible} for preserving symmetry.

For an $r\times s$ matrix $\pP:=\pP(z)$ with symmetry, we say that the symmetry type of $\pP$ is \emph{compatible} or $\pP$ has \emph{compatible symmetry} if
\begin{equation} \label{eq:CompSymDef}
\sym \pP(z) = \sym \pth_1^{\star}(z)\sym \pth_2(z)
\end{equation}
holds for some $1\times r$ and $1\times s$ row vectors of Laurent polynomials $\pth_1$ and $\pth_2$ with symmetry. It is easy to observe from the definition \eqref{eq:CompSymDef} that $ \sym\pth_1(z) $ gives the symmetry relationship between the rows of $ \pP(z) $, while $ \sym\pth_2(z) $ gives the symmetry relationship between the columns of $ \pP(z) $:
\begin{align}
\frac{\sym\pP_{j,l}(z)}{\sym\pP_{k,l}(z)} =& \frac{[\sym\pth^\star_1(z)]_j}{[\sym\pth^\star_1(z)]_k}, \quad
\forall l=1,\dots,s, \quad  \forall j,k =1,\ldots,r, \label{eq:CompSymDefRowCol1} \\
\frac{\sym\pP_{l,j}(z)}{\sym\pP_{l,k}(z)} =& \frac{[\sym\pth_2(z)]_j}{[\sym\pth_2(z)]_k}, \quad
\forall l=1,\dots,s, \quad  \forall j,k =1,\ldots,r. \label{eq:CompSymDefRowCol2}
\end{align}
For an $r\times s$ matrix $\pP$ and an $s\times t$ matrix $\pQ$ of Laurent polynomials with symmetry, we say \emph{the multiplication $\pP\pQ$ is compatible} if
$$\sym \pP(z) = \sym \pth_1^{\star}(z)\sym \pth_2(z),\quad \sym \pQ(z) = \sym \pth_2^{\star}(z)\sym \pth_3(z),$$
for some $1\times r$, $1\times s$ and $1\times t$ row vectors $\pth_1$, $\pth_2$ and $\pth_3$ of Laurent polynomials with symmetry. In this case, it is not hard to see that $\pP\pQ$ has compatible symmetry with
$$\sym[\pP\pQ](z)=\sym\pP(z)\sym\pQ(z)=\sym\pth_1^\star(z)\sym\pth_3(z).$$
For an $n\times n$ matrix $\pP$ of Laurent polynomials with compatible symmetry, one can show by induction that $\det(\pP)$ also has symmetry with type
\be\label{sym:det}
\sym[\det(\pP)](z)
=\prod_{j=1}^n\left([\sym\pth_1^\star(z)]_j[\sym\pth_2(z)]_j\right)
=\prod_{j=1}^n\sym\pP_{j,j}(z).
\ee
If we further assume that $\pP$ is \emph{strongly invertible}, i.e., $\det(\pP)$ is a nonzero monomial, then $\pP^{-1}$ is a matrix of Laurent polynomials with compatible symmetry. To calculate its symmetry type, recall that $\pP^{-1}=\det(\pP)^{-1}\adj(\pP)$. By \er{sym:det}, we have
$$[\sym[\adj(\pP)](z)]_{j,k}=\prod_{1\le l\le n, l\ne k}\prod_{1\le m\le n, m\ne j}[\sym\pth_1^\star(z)]_l[\sym\pth_2(z)]_m=\frac{\sym[\det(\pP)](z)}{[\sym\pth_1^\star(z)]_k[\sym\pth_2(z)]_j}=\sym[\det(\pP)](z)[\sym\pth_1(z)]_k[\sym\pth_2^\star(z)]_j.$$
It follows that $\pP^{-1}$ has compatible symmetry with
\be\label{p:inv:sym}\sym[\pP^{-1}](z)=\sym\pth_2^\star(z)\sym\pth_1(z).\ee
Consequently, if $\pP$, $\pQ$ and $\pR$ are $n\times n$, $n\times m$ and $m\times n$ matrices of Laurent polynomials with symmetry, and if $\pP$ is strongly invertible, then
\begin{itemize}
	\item the multiplication $\pP\pQ=:\pA$ is compatible implies that the multiplication $\pP^{-1}\pA$ is compatible;
	
	\item the multiplication $\pR\pP=:\pB$ is compatible implies that the multiplication $\pB\pP^{-1}$ is compatible.
\end{itemize}

Here we discuss some basic properties of compatibility. We summarize them as the following proposition.

\begin{prop}\label{sym:mat:com} Let $\pP$ be an $r\times s$ matrix of Laurent polynomials with compatible symmetry satisfying \er{eq:CompSymDef} for some $1\times r$ and $1\times s$ row vectors of Laurent polynomials $\pth_1$ and $\pth_2$ with symmetry.
\begin{enumerate}
\item[(1)] Let $\pP^{l,j}$ be the $s\times s$ permutation matrix such that $\pP\pP^{l,j}$ switches the $l$-th and $j$-th columns of $\pP$, and let $\ptP^{l,j}$ be the $r\times r$ permutation matrix such that $\ptP^{k,m}\pP$ switches the $k$-th and $m$-th rows of $\pP$. Then $\pP\pP^{l,j}$ and $\ptP^{k,m}\pP$ have compatible symmetry, with
$$\sym[\pP\pP^{l,j}](z)=\sym\pth_1^\star(z)\sym[\pth_2\pP^{l,j}](z),\quad \sym[\ptP^{k,m}\pP](z)=\sym[\pth_1(\ptP^{k,m})^{\tp}]^\star(z)\sym\pth_2(z).$$

\item[(2)] Let $\pD_r$ and $\pD_s$ be $r\times r$ and $s\times s$ diagonal matrices of Laurent polynomials with symmetry. Then $\pD_r\pP$ and $\pP\pD_s$ have compatible symmetry with
$$\sym[\pD_r\pP](z)=\sym[\pth_1\pD_r^\star]^\star(z)\sym\pth_2(z),\quad \sym[\pP\pD_s](z)=\sym\pth_1^\star(z)\sym[\pth_2\pD_s](z).$$

\item[(3)] Let $\pU$ be a $t\times r$ matrix of Laurent polynomials with symmetry. If there exists $k\in\{1,\dots,s\}$ such that
\be\label{su:k}\sym\pU_{j,1}(z)\sym\pP_{1,k}(z)=\sym\pU_{j,2}(z)\sym\pP_{2,k}(z)=\dots=\sym\pU_{j,r}(z)\sym\pP_{r,k}(z),\quad \forall j=1,\dots,t,\ee
then $\pU$ has compatible symmetry with $\sym\pU(z)=\sym\pth_3^\star(z)\sym\pth_1(z)$, where $\pth_3$ is any $1\times t$ row vector of Laurent polynomials with symmetry satisfying $[\sym\pth_3^\star(z)]_{j}=\frac{\sym\pU_{j,1}(z)}{[\sym\pth_1(z)]_1}$ for all $j=1,\dots,t$. Conversely, if the multiplication $\pU\pP$ is compatible, then there must exist $k\in\{1,\dots,s\}$ such that \er{su:k} holds.

\item[(4)] Let $\pU$ be a $s\times t$ matrix of Laurent polynomials with symmetry.
If there exists $l\in\{1,\dots,r\}$ such that
\be\label{su:l}\sym\pP_{l,1}(z)\sym\pU_{1,m}(z)=\sym\pP_{l,2}(z)\sym\pU_{2,m}(z)=\dots=\sym\pP_{l,s}(z)\sym\pU_{s,m}(z),\quad \forall m=1,\dots,t,\ee
then $\pU$ has compatible symmetry with $\sym\pU(z)=\sym\pth_2^\star(z)\sym\pth_3(z)$, where $\pth_3$ is any $1\times t$ row vector of Laurent polynomials with symmetry satisfying $[\sym\pth_3(z)]_{m}=\frac{\sym\pU_{1,m}(z)}{[\sym\pth_2^\star(z)]_1}$ for all $m=1,\dots, t$. Conversely, if the multiplication $\pP\pU$ is compatible, then there must exist $m\in\{1,\dots,s\}$ such that \er{su:l} holds.

\end{enumerate}
\end{prop}

\bp Items (1) and (2) can be easily verified by direct calculation. To prove item (3), by the fact that $\pP$ has compatible symmetry and using \er{su:k} and \er{eq:CompSymDefRowCol1}, we have
\be\label{su:k:2}\frac{\sym\pU_{j,q}(z)}{\sym\pU_{j,1}(z)}=\frac{\sym\pP_{1,k}(z)}{\sym\pP_{q,k}(z)}=\frac{[\sym\pth_1^\star(z)]_{1}}{[\sym\pth_1^\star(z)]_{q}}=\frac{[\sym\pth_1(z)]_{q}}{[\sym\pth_1(z)]_{1}},\quad \forall q=1,\dots,r,\quad j=1,\dots,t.\ee
Let $\pth_3$ be an $1\times t$ row vector of Laurent polynomials with symmetry satisfying $[\sym\pth_3^\star(z)]_{j}=\frac{\sym\pU_{j,1}(z)}{[\sym\pth_1(z)]_1}$ for all $j=1,\dots,t$. It follows from \er{su:k:2} that
$$\sym\pU_{j,q}(z)=\frac{[\sym\pth_1(z)]_{q}}{[\sym\pth_1(z)]_{1}}\sym\pU_{j,1}(z)=\frac{[\sym\pth_1(z)]_{q}}{[\sym\pth_1(z)]_{1}}[\sym\pth_3^\star(z)]_{j}[\sym\pth_1(z)]_1=[\sym\pth_3^\star(z)]_{j}[\sym\pth_1(z)]_q,$$
for all $q=1,\dots,r$ and $j=1,\dots,t$. Hence $\sym\pU(z)=\sym\pth_3^\star(z)\sym\pth_1(z)$, that is, $\pU$ has compatible symmetry and the multiplication $\pU\pP$ is compatible.

Conversely, suppose the multiplication $\pU\pP$ is compatible. This means there exists an $1\times t$ row vector $\pth_3$ of Laurent polynomials with symmetry such that $\sym\pU(z)=\sym\pth_3^\star(z)\sym\pth_1(z)$. Using \er{eq:CompSymDefRowCol1} and \er{eq:CompSymDefRowCol2}, we see that \er{su:k} in fact holds for all $k=1,\dots,s$. This proves item (3).

Finally, item (4) can be proved by using similar arguments as in th proof of item (3).\ep

\subsection{Factorization of matrices of Laurent polynomials with symmetry and negative constant determinants}

In this subsection, we consider a special case of Theorem~\ref{thm:SymSpecFactSym}, which is summarized as the following theorem.

\begin{theorem}\label{thm:SymSpecFactSym:sp}
Let $ \pA $ be a $ 2\times 2 $ Hermitian matrix of Laurent polynomials with compatible symmetry and define its symmetry type by
	$ \sym\pA(z) = \begin{bmatrix}
	1 & \mathsf{\alpha}(z) \\
	\mathsf{\alpha}^\star(z) & 1
	\end{bmatrix}$, where $\alpha(z):=\sym \pA_{1,2}(z)$.
Suppose $\det(\pA(z))=C< 0$ for all $z\in\T$ for some negative constant $C$. Then we can find a matrix $ \pU(z) = \begin{bmatrix}
	\pU_{1,1}(z) & \pU_{1,2}(z) \\
	\pU_{2,1}(z) & \pU_{2,2}(z)
	\end{bmatrix} $ of Laurent polynomials with compatible symmetry such that
	$\pA(z) = \pU(z)\DG(1,-1)\pU^\star(z)$ holds, and
	the symmetry type satisfies
	\begin{equation} \label{eq:Sym1:sp}
	 \frac{\sym\pU_{1,1}(z)}{\sym\pU_{2,1}(z)} = \frac{\sym\pU_{1,2}(z)}{\sym\pU_{2,2}(z)} = \alpha(z).
	\end{equation}
	
\end{theorem}

 To prove Theorem~\ref{thm:SymSpecFactSym:sp}, we need to establish several auxiliary results, which are the foundations for algorithms of constructing a matrix $\pU$ as in the statement of Theorem~\ref{thm:SymSpecFactSym:sp}.

\begin{lemma}\label{lem:ab}
Let $\pa$ and $\pb$ be Laurent polynomials with symmetry types $\sym\pa(z)=\eps_az^{c_a}$ and $\sym\pb(z)=\eps_bz^{c_b}$, respectively, for some $\eps_a,\eps_b\in\{\pm 1\}$ and $c_a,c_b\in\Z$. Suppose $\pb\ne \mathbf{0}$ and $\len(\pa)>\len(\pb)$. Then there exists a Laurent polynomial $\pq_1$ with symmetry such that
\be\label{eq:abq}\pa_1(z):=\pa(z)-\pb(z)\pq_1(z)\ee
satisfies $\len(\pa_1)<\len(\pa)$ and $\sym\pa_1(z)=\sym\pa(z)=\sym\pb(z)\sym\pq_1(z)$.
\end{lemma}
\bp Define $M_a:=\deg(\pa), m_a:=\ldeg(\pa), M_b:=\deg(\pb)$ and $m_b:=\ldeg(\pb)$. We have
$\pa(z)=\sum_{k=m_a}^{M_a}a(k)z^k$ and $\pb(z)=\sum_{k=m_b}^{M_b}b(k)z^k$.
The symmetry properties of $\pa$ and $\pb$ yields
$$a(m_a)=\eps_aa(M_a),\quad b(m_b)=\eps_bb(M_b),\quad M_a+m_a=c_a,\quad m_b+M_b=c_b.$$
Define
\be\label{q1}\pq_1(z):=\frac{a(M_a)}{b(M_b)}z^{M_a-M_b}+\frac{a(m_a)}{b(m_b)}z^{m_a-m_b}.\ee
We have
\begin{align*}\sym\pq_1(z)&=\frac{\pq_1(z)}{\pq_1(z^{-1})}=\frac{\frac{a(M_a)}{b(M_b)}z^{M_a-M_b}+\frac{a(m_a)}{b(m_b)}z^{m_a-m_b}}{\frac{a(M_a)}{b(M_b)}z^{M_b-M_a}+\frac{a(m_a)}{b(m_b)}z^{m_b-m_a}}=\frac{\frac{a(M_a)}{b(M_b)}z^{M_a-M_b}+\eps_a\eps_b\frac{a(M_a)}{b(M_b)}z^{(c_a-c_b)-(M_a-M_b)}}{\frac{a(M_a)}{b(M_b)}z^{M_b-M_a}+\eps_a\eps_b\frac{a(M_a)}{b(M_b)}z^{(M_a-M_b)-(c_a-c_b)}}\\
&=\eps_a\eps_bz^{c_1-c_b}=\frac{\sym\pa(z)}{\sym\pb(z)}.\end{align*}
Define $\pa_1$ as \er{eq:abq}, it follows immediately that $\sym\pa_1(z)=\sym\pa(z)=\sym\pb(z)\sym\pq_1(z)$.

We are left to show that $\len(\pa_1)<\len(\pa)$. As $\len(\pa)>\len(\pb)$, we have $M_a-M_b>m_a-m_b$. Thus by the definition of $\pq_1$ in \er{q1}, we have $\ldeg(\pq_1)=m_a-m_b$ and $\deg(\pq_1)=M_a-M_b$. So $\ldeg(\pb\pq_1)=m_a$ and $\deg(\pb\pq_1)=M_a$. Moreover, we obtain from \er{q1} that the coefficient of the $z^{M_a}$ term in $\pb\pq_1$ is $a(M_a)$, and the coefficient of the $z^{m_a}$ term in $\pb\pq_1$ is $a(m_a)$. Hence we conclude that $\ldeg(\pa_1)>m_a$ and $\deg(\pa_1)<M_a$, which implies $\len(\pa_1)<\len(\pa)$. This completes the proof.\ep

The next lemma provides a long division algorithm for Laurent polynomials with symmetry.

\begin{lemma}\label{lem:LD}
Let $\pa$ and $\pb$ be Laurent polynomials with symmetry types $\sym\pa(z)=\eps_az^{c_a}$ and $\sym\pb(z)=\eps_bz^{c_b}$ respectively, for some $\eps_a,\eps_b\in\{\pm 1\}$ and $c_a,c_b\in\Z$. Suppose $\pb\ne \pzr$. Consider the following four cases:
\begin{enumerate}
	\item[Case (1):] $\eps_a\eps_b=1$, $c_a-c_b\in2\Z+1$;
	\item[Case (2):] $\eps_a\eps_b=-1$, $c_a-c_b\in2\Z+1$;
	\item[Case (3):] $\eps_a\eps_b=1$, $c_a-c_b\in2\Z$;
	\item[Case (4):] $\eps_a\eps_b=-1$, $c_a-c_b\in2\Z$.
\end{enumerate}
Then the following results hold:
\begin{enumerate}
	\item[(i)] For cases (1) -- (3), we can construct a Laurent polynomial $\pq$ with symmetry such that
\be\label{eq:rab}\pr(z):=\pa(z)-\pb(z)\pq(z)\ee
satisfies $\len(\pr)<\len(\pb)$ and $\sym\pr(z)=\sym\pa(z)=\sym\pb(z)\sym\pq(z).$

\item[(ii)] For case (4), we can construct a Laurent polynomial $\pq$ with symmetry such that the Laurent polynomial $\pr$ defined in \er{eq:rab} satisfies $\len(\pr)\le\len(\pb)$ and $\sym\pr(z)=\sym\pa(z)=\sym\pb(z)\sym\pq(z).$
\end{enumerate}
\end{lemma}

\bp If $\len(\pa)<\len(\pb)$, simply take $\pq:=\mathbf{0}$ and $\pr:=\pa$, and all claims hold.

If $\len(\pa)\ge \len(\pb)$, define $\pa_0:=\pa$. Starting from $j=0$, if $\len(\pa_j)>\len(\pb)$, then apply Lemma~\ref{lem:ab} to construct a Laurent polynomial $\pq_j$ with symmetry such that
$$\pa_{j+1}(z):=\pa_j(z)-\pb(z)\pq_j(z)$$
satisfies $\len(\pa_{j+1})<\len(\pa_j)$ and $\sym\pa_{j+1}(z)=\sym\pa_j(z)=\sym\pb(z)\sym\pq_j(z)$. Note that this process cannot iterate forever: If $\len(\pa)=\len(\pb)$, the process terminates at the first step $j=0$ and we set $\pq_0:=\mathbf{0}$. Otherwise, the process stops at some step $j=K-1$ with $K\in\N$ when $\len(\pa_{K-1})>\len(\pb)\ge\len(\pa_K)$. In both cases, we have
\be\label{aK}\pa_K(z)=\pa_{K-1}(z)-\pb(z)\pq_{K-1}(z)=\dots=\pa_0(z)-\pb(z)(\pq_{K-1}(z)+\dots+\pq_0(z)).\ee
Moreover, the symmetry types of $\pa_j$ satisfy
$$\sym\pa_{K}(z)=\sym\pa_{K-1}(z)=\dots=\sym\pa_0(z)=\sym\pa(z),$$
$$\sym\pq_j(z)=\frac{\sym\pa_{j+1}(z)}{\sym\pb(z)}=\frac{\sym\pa(z)}{\sym\pb(z)},\qquad \forall j=0,\dots,K-1.$$
Define
\be\label{ptq}\ptq(z):=\pq_{K-1}(z)+\dots+\pq_0(z).\ee
It follows that
\be\label{ptq:ab}\sym\ptq(z)=\frac{\sym\pa(z)}{\sym\pb(z)}\quad\mbox{and}\quad \pa_K(z)=\pa(z)-\pb(z)\ptq(z).\ee
We now consider cases (1)--(4) separately in the following way:

\begin{itemize}
\item \textbf{Cases (1) and (2):} The condition $c_a-c_b\in2\Z+1$ implies that $\odd(c_a)\ne\odd(c_b)$. Thus by Proposition~\ref{prop:sym:lau}, we have $\len(\pa)\ne\len(\pb)$. Therefore, we only consider the situation when $\len(\pa)>\len(\pb)$. Define $\pq:=\ptq$ and $\pr:=\pa_K$, where $\ptq$ and $\pa_K$ are defined as \er{ptq} and \er{aK} respectively. By \er{ptq:ab}, we have $$\sym\pr(z)=\sym\pa_{K}(z)=\sym\pa(z)=\sym\pb(z)\sym\pq(z).$$
Moreover, we have $\len(\pr)=\len(\pa_K)\le\len(\pb)$. On the other hand, note that $\sym\pr(z)=\eps_az^{c_a}$. It follows from Proposition~\ref{prop:sym:lau} that $\len(\pr)\ne\len(\pb)$, and thus $\len(\pr)<\len(\pb)$. This proves the claims for cases (1) and (2) in item (i).

\item \textbf{Case (3):} Define $\ptq$ and $\pa_K$ as \er{ptq} and \er{aK} respectively. We have $\len(\pa_K)\le\len(\pb)$ and $\sym\pa_K(z)=\sym\pa(z)=\sym\pb(z)\sym\ptq(z)$. If $\len(\pa_K)<\len(\pb)$, simply set $\pa:=\ptq$ and $\pr:=\pa_K$, and all claims in item (1) hold. If $\len(\pa_K)=\len(\pb)$, define $\ldeg(\pa_K):=m_{a,K}, \deg(\pa_K):=M_{a,K}, \ldeg(\pb):=m_b, \deg(\pb):=M_b$, and define
$$\mrpq(z):=\frac{a_K(M_{a,K})}{b(M_b)}z^{M_{a,K}-M_b},$$
\be\label{r:K}\pr(z):=\pa_K(z)-\pb(z)\mrpq(z).\ee
From $\len(\pa_K)=\len(\pb)$, we have $M_{a,k}-M_b=m_{a,K}-m_b$. Thus $\deg(\pb\mrpq)=M_b+\deg(\mrpq)=M_{a,K}$ and $\ldeg(\pb\mrpq)=m_b+(M_{a,K}-M_b)=m_b+(m_{a,K}-m_b)=m_{a,K}$. Hence from the definitions of $\mrpq$ and $\pr$, it is easy to see that $\ldeg(\pr)>m_{a,K}$ and $\deg(\pr)<M_{a,K}$. In particular, $\len(\pr)<\len(\pa_K)=\len(\pb)$.

Next, note that $\sym\mrpq(z)=z^{2M_{a,K}-2M_b}$. As $\eps_a\eps_b=1$, we have $\eps_a=\eps_b$. Using Proposition~\ref{prop:sym:lau}, direct calculation yields
\begin{align*}\sym\pb(z)\sym\mrpq(z)&=\eps_bz^{c_b+2M_{a,K}-2M_b}=\eps_bz^{m_b+M_b+2M_{a,K}-2M_b}=\eps_bz^{2M_{a,K}-\len(\pb)}\\
&=\eps_bz^{2M_{a,K}-\len(\pa_K)}=\eps_az^{m_{a,K}+M_{a,K}}=\sym\pa_K(z).
\end{align*}
Define $\pq:=\ptq+\mrpq$. We have
$$\pr(z)=\pa_K(z)-\pb(z)\mrpq(z)=\pa(z)-\pb(z)\ptq(z)-\pb(z)\mrpq(z)=\pa(z)-\pb(z)\pq(z),$$
and $\sym\pr(z)=\sym\pa_K(z)=\sym\pa(z)=\sym\pb(z)\sym\pq(z).$ This proves all claims for case (3) in item (i).

\item \textbf{Case 4:} Define $\pq:=\ptq$ and $\pr:=\pa_K$, where $\ptq$ and $\pa_K$ are defined as \er{ptq} and \er{aK} respectively. All claims of item (ii) follows immediately.

\end{itemize}

The proof is now complete.\ep

We move on to perform some further analysis on $2\times 2$ Hermitian matrices of Laurent polynomials with symmetry. Let $\pA(z)$ be such a matrix. Then as $\pA(z)=\pA^\star(z)$ for all $z\in\T$, its symmetry type is of the form $\sym\pA=\begin{bmatrix}1 & \alpha(z)\\
\alpha^\star(z) &1\end{bmatrix}$ with four possible cases for $\alpha(z)=\eps z^c$:
\begin{enumerate}
	\item[Case 1:] $\eps=1$, $c\in2\Z+1$;
	\item[Case 2:] $\eps=-1$, $c\in2\Z+1$;
	\item[Case 3:] $\eps=1$, $c\in2\Z$;
	\item[Case 4:] $\eps=-1$, $c\in2\Z$.
\end{enumerate}
To prove the main result Theorem~\ref{thm:SymSpecFactSym:sp}, all the cases 1 -- 3 can be investigated together, while case 4 need to be dealt with individually. For cases 1 -- 3, we have the following technical lemma, which plays a key role for  constructing a matrix $\pU$ as required in Theorem~\ref{thm:SymSpecFactSym:sp}.

\begin{lemma}\label{lem:c1to3}
Let $\pA$ be a $2\times 2$ Hermitian matrix	 of Laurent polynomials with symmetry type $\sym\pA=\begin{bmatrix}1 & \alpha(z)\\
\alpha^\star(z) &1\end{bmatrix}$, where $\alpha$ is one of the following:
\begin{enumerate}
\item[(1)] $\alpha(z)=\eps z^{2k+1}$ for some $\eps\in\{\pm 1\}$ and $k\in\Z$;
\item[(2)] $\alpha(z)=z^{2k}$ for some $k\in\Z$.
\end{enumerate}	
Suppose $\det(\pA)=C$ for some negative constant $C$. Then there exist a positive integer $K$ and $2\times 2$ strongly invertible matrices $\pV^{(0)},\dots,\pV^{(K-1)}$ of Laurent polynomials with compatible symmetry, such that
$$\pA^{(K)}(z):=\pV^{(K-1)}(z)\cdots\pV^{(0)}
\pA(z)\pV^{(0)\star}(z)\cdots\pV^{(K-1)\star}(z)$$
satisfies
\begin{enumerate}
	\item[(i)] all the multiplications are compatible;
	\item[(ii)] at least one entry of $\pA^{(K)}$ is $\pzr$.
\end{enumerate}
\end{lemma}
\bp If $\pA$ has one zero entry, simple take $K=1$, $\pV^{(0)}=\pI_2$ and $\pA^{(1)}=\pA$. Here $\pI_n$ stands for the $n\times n$ identity matrix.
Otherwise, all for entries of $\pA$ are nonzero, and we claim that
\be\label{len:pA}\len(\pA_{1,1})+\len(\pA_{2,2})=\len(\pA_{1,2})+\len(\pA_{2,1})=2\len(\pA_{2,1}).\ee
For any nonzero Laurent polynomial $\pu$, define its \emph{support} as
\be\label{fsupp}\fsupp(\pu):=[\ldeg(\pu),\deg(\pu)].\ee

As $\pA$ is Hermitian, it is easy to see that $\fsupp(\pA_{1,1}\pA_{2,2})$ and $\fsupp(\pA_{1,2}\pA_{2,1})$ are symmetric intervals with centre $0$. Define
$$[-m,m]:=\fsupp(\pA_{1,1}\pA_{2,2}),\qquad [-n,n]:=\fsupp(\pA_{1,2}\pA_{2,1}).$$
Since $\det(\pA)=\pA_{1,1}\pA_{2,2}-\pA_{1,2}\pA_{2,1}=C$ is a constant, we must have $m=n$. It follows that
$$\len(\pA_{1,1})+\len(\pA_{2,2})=\len(\pA_{1,1}\pA_{2,2})=2m=2n=\len(\pA_{1,2}\pA_{2,1})=\len(\pA_{1,2})+\len(\pA_{2,1})=2\len(\pA_{2,1}).$$
This proves \er{len:pA}. As a consequence, we have
\be\label{min:pA}\min\{\len(\pA_{1,1}),\len(\pA_{2,2})\}\le \len(\pA_{2,1}).\ee
We now construct a strongly invertible matrix $\pV^{(0)}$ of Laurent polynomials with compatible symmetry, such that $\pA^{(1)}:=\pV^{(0)}\pA\pV^{(0)\star}$ satisfies $\len(\pA^{(1)}_{2,1})<\len(\pA_{2,1})$ and the multiplications are compatible. We consider the following two cases:
\begin{itemize}
	\item If $\len(\pA_{1,1})\le \len(\pA_{2,1})$, as $\sym\pA_{1,1}(z)=1$ and $\sym\pA_{2,1}(z)=\alpha^\star$ with $\alpha$ satisfying (1) or (2), we can apply item (i) of Lemma~\ref{lem:LD} to construct a Laurent polynomial $\pQ$ with symmetry, such that
$$\pr(z):=\pA_{2,1}(z)-\pA_{1,1}(z)\pq(z)$$
satisfies $\len(\pr)<\len(\pA_{1,1})$ and $\sym\pr(z)=\sym\pq(z)\sym\pA_{1,1}(z)=\sym\pA_{2,1}(z).$ Define
$$\pV^{(0)}(z):=\begin{bmatrix} 1 & 0\\
-\pq(z) & 1\end{bmatrix}.$$
Then
$$\pA^{(1)}(z):=\pV^{(0)}(z)\pA(z)\pV^{(0)\star}(z)=\begin{bmatrix}\pA_{1,1}(z) & \pr^\star(z)\\
\pr(z) &\pA_{2,2}^{(1)}(z)\end{bmatrix}.$$
So
\be\label{len:pA1:1}\len(\pA_{2,1}^{(1)})=\len(\pr)<\len(\pA_{1,1})\le\len(\pA_{2,1}).\ee

\item If $\len(\pA_{2,2})\le \len(\pA_{2,1})$, we apply item (i) of Lemma~\ref{lem:LD} to construct a Laurent polynomial $\pq$ with symmetry, such that
$$\pr(z):=\pA_{1,2}(z)-\pA_{2,2}(z)\pq(z)$$
satisfies $\len(\pr)<\len(\pA_{2,2})$ and $\sym\pr(z)=\sym\pq(z)\sym\pA_{2,2}(z)=\sym\pA_{1,2}(z).$ Define
$$\pV^{(0)}(z):=\begin{bmatrix} 1 & -\pq(z)\\
0 & 1\end{bmatrix}.$$
Then
$$\pA^{(1)}(z):=\pV^{(0)}(z)\pA(z)\pV^{(0)\star}(z)=\begin{bmatrix}\pA_{1,1}^{(1)}(z) & \pr(z)\\
\pr^\star(z) &\pA_{2,2}(z)\end{bmatrix}.$$
So
\be\label{len:pA1:2}\len(\pA_{2,1}^{(1)})=\len(\pr)<\len(\pA_{2,2})\le\len(\pA_{2,1}).\ee

\end{itemize}
In both cases above,  $\sym\pV^{(0)}(z)=\sym\pA(z)=\sym\pth_1^\star(z)\sym\pth_1(z)$ where $\sym\pth_1(z):=[1,\alpha(z)]$. This shows that the multiplications $\pA^{(1)}=\pV^{(0)}\pA\pV^{(0)\star}$ are compatible. Hence we have found a desired $\pV^{(0)}$ as required.

Now if $\pA^{(1)}$ has a zero entry, we stop here. Otherwise, we repeat the above process to further shrink the length of $\pA_{2,1}^{(1)}$. In general, if all entries of $\pA^{(j)}$ are nonzero, we apply the above process to find a strongly invertible matrix $\pV^{(j)}$ such that $\pA^{(j+1)}:=\pV^{(j)}\pA^{(j)}\pV^{(j)\star}$ satisfies $\len(\pA^{(j+1)}_{2,1})<\len(\pA_{2,1}^{(j)})$ and the multiplications are compatible. Note that $\len(\pA_{2,1})>\len(\pA_{2,1}^{(1)})>\len(\pA_{2,1}^{(2)})>\dots$, so the iteration terminates within finite steps, say at step $j=K-1$, so that the matrix
$$\pA^{(K)}:=\pV^{(K-1)}\pA^{(K-1)}\pV^{(K-1)\star}=\dots=\pV^{(K-1)}\dots\pV^{(0)}\pA\pV^{(0)\star}\dots\pV^{(K-1)\star}$$
has one zero entry. Moreover, it can be easily seen from previous discussions that all the multiplications above are compatible, with
\be\label{sym:vj}\sym\pA^{(j)}=\sym\pA=\sym\pV^{(j)}=\sym\pth_1^\star\sym\pth_1,\qquad j=0,\dots,K-1,\ee
where $\pA^{(0)}:=\pA$ and $\sym\pth_1(z):=[1,\alpha(z)]$. This completes the proof.\ep

We now have all the necessary tools to prove Theorem~\ref{thm:SymSpecFactSym:sp}.

\bp[\textbf{Proof of Theorem~\ref{thm:SymSpecFactSym:sp} for cases 1--3}]
We constructively prove Theorem~\ref{thm:SymSpecFactSym:sp} for cases 1--3 by the following construction algorithm for $\pU$ in Theorem~\ref{thm:SymSpecFactSym:sp} for cases 1--3
as follows.

\begin{algorithm}\label{alg:c1to3}Suppose $\alpha(z)=\eps z^{2k+1}$ for some $\eps\in\{\pm 1\}$ and $k\in\Z$ or $\alpha(z)=z^{2k}$ for some $k\in\Z$.
\begin{enumerate}
\item[Step 1.] Apply Lemma~\ref{lem:c1to3} to construct $2\times2$ strongly invertible matrices $\pV^{(0)},\dots,\pV^{(K-1)}$ of Laurent polynomials with compatible symmetry such that
$$\pA^{(K)}(z):=\pV^{(K-1)}(z)\dots\pV^{(0)}\pA(z)\pV^{(0)\star}(z)\dots\pV^{(K-1)\star}(z)$$
has at least one zero entry, and all above multiplications are compatible.

\item[Step 2.] Construct the $2\times 2$ matrix $\ptU$ of Laurent polynomials with compatible symmetry such that $\pA^{(K)}=\ptU \DG(1,-1)\ptU^{\star}$. The explicit expression of $\ptU$ is given by the following:
\begin{enumerate}
	\item[(1)] If $\pA^{(K)}_{1,2}=\pA^{(K)}_{2,1}=\pzr$, then $\pA_{1,1}^{(K)}=c_1$ and $\pA_{2,2}^{(K)}=c_2$ for some constants $c_1$ and $c_2$ satisfying $c_1c_2=\det(\pA^{(K)})=\det(\pA)=C<0$. Set \be\label{ptu:1}\ptU:=\begin{cases}\DG(\sqrt{c_1},\sqrt{-c_2}), &\text{if }c_1>0,c_2<0,\\[0.2cm]
	\begin{bmatrix}0 &\sqrt{-c_1}\\
	\sqrt{c_2} & 0\end{bmatrix},&\text{if }c_1<0,c_2>0;\end{cases}\ee
	
	\item[(2)] If $\pA^{(K)}_{1,1}=0$, set
	 \be\label{ptu:2}\ptU:=\frac{1}{\sqrt{2}}\begin{bmatrix}\pA^{(K)}_{1,2} & \pA^{(K)}_{1,2}\\
	 \frac{1}{2}\pA^{(K)}_{2,2}+1 & \frac{1}{2}\pA^{(K)}_{2,2}-1\end{bmatrix};\ee
	
	\item[(3)] If $\pA^{(K)}_{2,2}=0$, set
\be\label{ptu:3}\ptU:=\frac{1}{\sqrt{2}}\begin{bmatrix}\frac{1}{2}\pA^{(K)}_{1,1}+1 & \frac{1}{2}\pA^{(K)}_{1,1}-1\\
	\pA^{(K)}_{2,1} & \pA^{(K)}_{2,1}\end{bmatrix}.\ee
	
\end{enumerate}

\item[Step 3.] The matrix $\pU:=\pV^{(K-1)}\dots\pV^{(0)}\ptU$ is the one satisfying all claims of Theorem~\ref{thm:SymSpecFactSym:sp}.

	\end{enumerate}
		
\end{algorithm}

\textbf{Justification of Algorithm~\ref{alg:c1to3}:} By the assumptions of the algorithm, the validity of step 1 is trivial.

For step 2: If $\pA_{1,2}^{(K)}=\pA_{2,1}^{(K)}=0$, by $\sym\pA^{(K)}=\sym\pA$, we have $\sym\pA_{1,1}^{(K)}=\sym\pA_{2,2}^{(K)}=1$. Note that $\pA_{1,1}^{(K)}\pA_{2,2}^{(K)}=\det(\pA^{(K)})=\det(\pA)=C<0$. This means $\pA_{1,1}^{(K)}$ and $\pA_{2,2}^{(K)}$ are constants, say $\pA_{1,1}^{(K)}=c_1$ and $\pA_{2,2}^{(K)}=c_2$. Moreover, $c_1$ and $c_2$ have opposite signs. Define $\ptU$ as \er{ptu:1}, direct calculation shows that $\pA^{(K)}=\ptU\DG(1,-1)\ptU^\star$. Furthermore, since $0$ has any types of symmetry, we see that $\ptU$ has compatible symmetry with $\sym\ptU=\sym\pth_1^\star\sym\pth_2$, where $\sym\pth_1:=[1,\alpha]$ and $\sym\pth_2 :=[\alpha,1]$. If $\pA_{1,1}^{(K)}=0$, define $\ptU$ as \er{ptu:2}. Then it is straightforward to verify that $\pA^{(K)}=\ptU\DG(1,-1)\ptU^\star$ and $\sym\ptU=\sym\pth_1^\star\sym\pth_3$ where $\sym\pth_3 :=[\alpha,\alpha]$. Lastly, if $\pA_{2,2}^{(K)}=0$, it is easy to see that $\pA^{(K)}=\ptU\DG(1,-1)\ptU^\star$ with $\ptU$ being defined as \er{ptu:3}, and $\sym\ptU=\sym\pth_1^\star\sym\pth_4$ where $\sym\pth_4 :=[1,1]$.

Finally, define $\pU$ as in step 3, it follows immediately that $\pU$ has compatible symmetry and satisfies $\pA=\pU\DG(1,-1)\pU^{\star}$. Moreover, it follows from \er{sym:vj} that all the multiplication of $\pU=\pV^{(K-1)}\dots\pV^{(0)}\ptU$ are compatible, and $\sym\pU=\sym\pth_1^\star\sym\pth_l$ where $l=2,3$ or $4$. Hence \er{eq:Sym1:sp} holds. This completes the justification of Algorithm~\ref{alg:c1to3}. \ep

\bp[\textbf{Proof of Theorem~\ref{thm:SymSpecFactSym:sp} for case 4}] This is a harder case, and the construction for $\pU$ in Theorem~\ref{thm:SymSpecFactSym:sp} for case 4 is summarized by the following technical algorithm.

\begin{algorithm}\label{alg:c4}
Suppose $\alpha(z)=-z^{2k}$ for some $k\in\Z$.
\begin{enumerate}
\item[(S1)] Define $[-n,n]:=\fsupp(\pA_{1,1})$. Parametrize Laurent polynomials
$$\pV(z)=\sum_{j=0}^nV(j)z^j,\qquad \pW(z)=\sum_{j=0}^nW(j)z^j,$$
where $V(0),\dots,V(n),W(0),\dots,W(n)$ are free parameters. Let $\pQ$ and $\pR$ be the Laurent polynomials which are determined uniquely through long division using $\pA_{1,1}$ by
\be\label{ld:a11}\pA_{2,1}(z)\pV(z)+dz^{n-k}\pW^\star(z)=\pA_{1,1}(z)\pQ(z)+\pR(z)\text{ with }\fsupp(\pR)\subseteq [-n,n-1],\ee
where $d:=\sqrt{-C}=\sqrt{-\det(\pA)}$. Then the linear space
$$E:=\{(V(0),\dots,V(n),W(0),\dots,W(n)):\pR=\mathbf{0}\}\subseteq \C^{2n+2}$$
has dimension at least two. Take any nonzero element $(t_0,\dots,t_n,s_0,\dots,s_n)\in E$ and define
$$\ptU_{1,1}(z):=\sum_{k=0}^nt_kz^k,\qquad \ptU_{1,2}(z):=\sum_{k=0}^n s_kz^k.$$

\item[(S2)] If $\sym\ptU_{1,1}(z)=z^n$ and $\sym\ptU_{1,2}=-z^n$, set $\mrpU_{1,1}:=\ptU_{1,1}$ and $\mrpU_{1,2}:=\ptU_{1,2}$. Otherwise, define
\be\label{up:u11:u12}\mrpU_{1,1}(z):=\frac{\ptU_{1,1}(z)-z^n\ptU_{1,1}^\star(z)}{2},\qquad \mrpU_{1,2}(z):=\frac{\ptU_{1,2}(z)+z^n\ptU_{1,2}^\star(z)}{2}.\ee
Then the updated $\mrpU_{1,1}$ and $\mrpU_{1,2}$ are not both $\mathbf{0}$, and their symmetry types are $\sym\mrpU_{1,1}(z)=-z^n$ and $\sym\mrpU_{1,2}(z)=z^n$.

\item[(S3)] There exists $\lambda\in\R\setminus\{0\}$ such that
\be\label{a11:lambda}\pA_{1,1}(z)=\lambda\left(\mrpU_{1,1}(z)\mrpU_{1,1}^\star(z)-\mrpU_{1,2}(z)\mrpU_{1,2}^\star(z)\right).\ee
If $\lambda>0$, define
\be\label{up:u11:u12:2}\pU_{1,1}(z):=\sqrt{\lambda}\mrpU_{1,1}(z),\qquad \pU_{1,2}(z):=\sqrt{\lambda}\mrpU_{1,2}(z).\ee
Otherwise, define
\be\label{up:u11:u12:3}\pU_{1,1}(z):=\sqrt{-\lambda}\mrpU_{1,2}(z),\qquad \pU_{1,2}(z):=\sqrt{-\lambda}\mrpU_{1,1}(z).\ee
Now we have $\sym\pU_{1,1}(z)=\eps z^n$ and $\sym\pU_{1,2}(z)=-\eps z^n$ for some $\eps\in\{\pm 1\}$.

\item[(S4)] Define
\be\label{u21}\pU_{2,1}(z):=\frac{\pA_{2,1}(z)\pU_{1,1}(z)+dz^{n-k}\pU_{1,2}^\star(z)}{\pA_{1,1}(z)},\ee
\be\label{u22}\pU_{2,2}(z):=\frac{\pA_{2,1}(z)\pU_{1,2}(z)+dz^{n-k}\pU_{1,1}^\star(z)}{\pA_{1,1}(z)}.\ee
Then $\pU_{2,1}$ and $\pU_{2,2}$ are well-defined Laurent polynomials with symmetry types $\sym\pU_{2,1}(z)=-\eps z^{n-2k}$ and $\sym\pU_{2,2}(z)=\eps z^{n-2k}$.  Set $\pU:=\begin{bmatrix}\pU_{1,1} & \pU_{1,2}\\
\pU_{2,1} & \pU_{2,2}\end{bmatrix}$. Then $\pU$ is a desired matrix of Laurent polynomials with compatible symmetry satisfying all claims of Theorem~\ref{thm:SymSpecFactSym:sp}.
\end{enumerate}
\end{algorithm}

\textbf{Justification of Algorithm~\ref{alg:c4}:} For (S1): First note that $\pA_{1,1}\ne\mathbf{0}$. Otherwise, we have $\det(\pA)=-\pA_{1,2}\pA_{2,1}=-d^2$. This in particular means that $\pA_{1,2}$ is a nonzero monomial, which cannot have symmetry type $-z^c$. Hence, $\fsupp(\pA_{1,1})$ is a non-empty interval and is symmetric with centre $0$.

By our choice of $\pR$ and $\fsupp(\pR)\subseteq[-n,n-1]$, the condition $\pR=0$ induces a homogeneous linear system of at most $2n$ equations with $2n+2$ unknowns. Therefore, the space of all solutions, namely the space $E$, has dimension at least two. In particular, $E$ contains non-trivial solutions, and the Laurent polynomials $\ptU_{1,1}$ and $\ptU_{1,2}$ in (S1) can be defined.

For (S2): By our choices of $\ptU_{1,1}$ and $\ptU_{1,2}$, we have
\be\label{ptq:1}\pA_{2,1}(z)\ptU_{1,1}(z)+dz^{n-k}\ptU_{1,2}^\star(z)=\pA_{1,1}(z)\ptU_{2,1}(z),\ee
where $\ptU_{2,1}(z):=\frac{\pA_{2,1}(z)\ptU_{1,1}(z)+dz^{n-k}\ptU_{1,2}^\star(z)}{\pA_{1,1}(z)}$ is a well-defined Laurent polynomial. By the symmetry type of $\pA$, we have $\pA_{2,1}^\star(z)=-z^{2k}\pA_{2,1}(z)$. Taking Hermitian conjugates and multiplying $z^{n-2k}$ on both sides of \er{ptq:1}, we have
\be\label{ptq:2}-z^n\pA_{2,1}(z)\ptU_{1,1}^\star(z)+dz^{-k}\ptU_{1,2}(z)=z^{n-2k}\pA_{1,1}(z)\ptU_{2,1}^\star(z).\ee
It follows from \er{ptq:1} and \er{ptq:2} that
\be\label{ptq:3}\pA_{2,1}(z)\left(\ptU_{1,1}(z)-z^n\ptU_{1,1}^\star(z)\right)+dz^{n-k}\left(\ptU_{1,2}^\star(z)+z^{-n}\ptU_{1,2}(z)\right)=\pA_{1,1}(z)\left(\ptU_{2,1}(z)+z^{n-2k}\ptU_{2,1}^\star(z)\right).\ee
After (S1), if $\sym\ptU_{1,1}(z)\ne z^n$ or $\sym\ptU_{1,2}(z)\ne -z^n$, then either $\sym\ptU_{1,1}(z)\ne z^{2n}\sym\ptU_{1,1}^\star(z)$ or $\sym\ptU_{1,2}(z)\ne -z^{2n}\sym\ptU_{1,2}^\star(z)$. Define $\mrpU_{1,1}$ and $\mrpU_{1,2}$ as in \er{up:u11:u12}, it is trivial that $\mrpU_{1,1}$ and $\mrpU_{1,2}$ cannot be both zero. Moreover, direct calculation yields $\sym\mrpU_{1,1}=-z^n$ and $\sym\mrpU_{1,2}=z^n$.

For (S3): Define $\mrpU_{2,1}(z):=\frac{\ptU_{2,1}(z)+z^{n-2k}\ptU_{2,1}^\star(z)}{2}$. Recall that $\det(\pA)=-d^2$. It follows that
\begin{align*}&\pA_{1,2}(z)\left(\pA_{2,1}(z)\mrpU_{1,2}(z)+dz^{n-k}\mrpU_{1,1}^\star(z)\right)=\left(\pA_{1,1}(z)\pA_{2,2}(z)+d^2\right)\mrpU_{1,2}(z)+dz^{n-k}\mrpU_{1,1}^\star(z) \pA_{1,2}(z)\\
=&\pA_{1,1}(z)\pA_{2,2}(z)\mrpU_{1,2}(z)+dz^{n-k}\left(dz^{n-k}\mrpU_{1,2}^\star(z)+\mrpU_{1,1}(z)\pA_{1,2}^{\star}(z)\right)^\star=\pA_{1,1}(z)\left(\pA_{2,2}(z)\mrpU_{1,2}(z)+dz^{n-k}\mrpU_{2,1}^\star(z)\right).\end{align*}
As $\det(\pA)=\pA_{1,1}\pA_{2,2}-\pA_{1,2}\pA_{2,1}$ is a nonzero constant, we see that $\gcd(\pA_{1,1},\pA_{1,2})=1$. Thus from the above equation, we see that $\pA_{1,1}(z)$ divides $\pA_{2,1}(z)\mrpU_{1,2}(z)+dz^{n-k}\mrpU_{1,1}^\star(z)$.
Define $\mrpU_{2,2}(z):=\frac{\pA_{2,1}(z)\mrpU_{1,2}(z)+dz^{n-k}\mrpU_{1,1}^\star(z)}{\pA_{1,1}(z)}$. We have
\be\label{ptu22}\pA_{2,1}(z)\mrpU_{1,2}(z)+dz^{n-k}\mrpU_{1,1}^\star(z)=\pA_{1,1}(z)\mrpU_{2,2}(z).\ee
It follows from \er{ptq:3} and \er{ptu22} that
\be\label{ptu:mat}\begin{bmatrix}\mrpU_{2,2}(z) & -\mrpU_{1,2}(z)\\
	-\mrpU_{2,1}(z) &\mrpU_{1,1}(z)\end{bmatrix}\begin{bmatrix}\pA_{1,1}(z)\\
	 \pA_{2,1}(z)\end{bmatrix}=dz^{n-k}\begin{bmatrix}\mrpU_{1,1}^\star(z)\\
	 -\mrpU_{1,2}^\star(z)\end{bmatrix}.\ee
Multiplying $\left[\mrpU_{1,2}^\star, \mrpU_{1,1}^\star\right]$ to the left on both sides of \er{ptu:mat} yields
\be\label{ptu:0}\left(\mrpU_{2,1}(z)\mrpU_{1,1}^\star(z)-\mrpU_{2,2}(z)\mrpU_{1,2}^\star(z)\right)\pA_{1,1}(z)=\left(\mrpU_{1,1}(z)\mrpU_{1,1}^\star(z)-\mrpU_{1,2}(z)\mrpU_{1,2}^\star(z)\right)\pA_{2,1}(z).\ee
Note that $\gcd(\pA_{1,1},\pA_{2,1})=1$. So $\pA_{1,1}(z)$ divides $\mrpU_{1,1}(z)\mrpU_{1,1}^\star(z)-\mrpU_{1,2}(z)\mrpU_{1,2}^\star(z)$.

On the other hand, since $\fsupp(\ptU_{1,1})\subseteq [0,n]$ and $\fsupp(\ptU_{1,2})\subseteq[0,n]$, we see that $\fsupp(\mrpU_{1,1})\subseteq [0,n]$ and $\fsupp(\mrpU_{1,2})\subseteq[0,n]$. Thus $\fsupp(\mrpU_{1,1}\mrpU_{1,1}^\star)\subseteq [-n,n]$ and $\fsupp(\mrpU_{1,2}\mrpU_{1,2}^\star)\subseteq [-n,n]$. So $\fsupp(\mrpU_{1,1}\mrpU_{1,1}^\star-\mrpU_{1,2}\mrpU_{1,2}^\star)\subseteq[-n,n]=\fsupp(\pA_{1,1})$.

We claim that $\mrpU_{1,1}\mrpU_{1,1}^\star-\mrpU_{1,2}\mrpU_{1,2}^\star\ne\mathbf{0}$, and thus by what we have got above together with the fact that both $\pA_{1,1}$ and $\mrpU_{1,1}\mrpU_{1,1}^\star-\mrpU_{1,2}\mrpU_{1,2}^\star$ are Hermitian, we see that \er{a11:lambda} must hold for some $\lambda\in\R\setminus\{0\}$. In fact, by the choices of $\mrpU_{1,1}$ and $\mrpU_{1,2}$, either $\sym\mrpU_{1,1}(z)=z^n$ and $\sym\mrpU_{1,2}(z)=-z^n$, or vice versa. If $\sym\mrpU_{1,1}(z)=z^n$ and $\sym\mrpU_{1,2}(z)=-z^n$, then by  item (1) of Lemma~\ref{lem:rt:pm1}, we have $\mz(\mrpU_{1,1},1)=\mz(\mrpU_{1,1}^\star,1)\in2\Z$ and $\mz(\mrpU_{1,2},1)=\mz(\mrpU_{1,2}^\star,1)\in2\Z+1$. This means $\mz(\mrpU_{1,1}\mrpU_{1,1}^\star,1)\in 4\Z$ and $\mz(\mrpU_{1,2}\mrpU_{1,2}^\star,1)\in 4\Z+2$. Therefore $\mrpU_{1,1}\mrpU_{1,1}^\star-\mrpU_{1,2}\mrpU_{1,2}^\star\ne\mathbf{0}$. The case when $\sym\mrpU_{1,1}=-z^n$ and $\sym\mrpU_{1,2}=z^n$ can be proved similarly. This proves the claim.

Now we define $\pU_{1,1}$ and $\pU_{1,2}$ as \er{up:u11:u12:2} if $\lambda>0$ or as \er{up:u11:u12:3} if $\lambda<0$. We have $\sym\pU_{1,1}(z)=\eps z^n$ and $\sym\pU_{1,2}(z)=-\eps z^n$ for some $\eps\in\pm\{1\}$, and
\be\label{pa11}\pA_{1,1}(z)=\pU_{1,1}(z)\pU_{1,1}^\star(z)-\pU_{1,2}(z)\pU_{1,2}^\star(z).\ee

For (S4): By \er{ptq:3}, \er{ptu22} and our choices of $\pU_{1,1}$ and $\pU_{1,2}$, one can see that $\pA_{1,1}(z)$ divides $\pA_{2,1}(z)\pU_{1,1}(z)+dz^{n-k}\pU_{1,2}^\star(z)$ and $\pA_{2,1}(z)\pU_{1,2}(z)+dz^{n-k}\pU_{1,1}^\star(z)$. So $\pU_{2,1}$ and $\pU_{2,2}$ are well-defined Laurent polynomials. Since $\sym\pA_{2,1}(z)=-z^{-2k}$ and $\sym\pA_{1,1}(z)=1$, it is easy to verify that $\sym\pU_{2,1}(z)=-\eps z^{n-2k}$ and $\sym\pU_{2,2}(z)=\eps z^{n-2k}$.

Define $\pU:=\begin{bmatrix}\pU_{1,1} & \pU_{1,2}\\
\pU_{2,1} & \pU_{2,2}\end{bmatrix}$. Observe that \er{u21} and \er{u22} is equivalent to
\be\label{pu:mat}\begin{bmatrix}\pU_{2,2}(z) & -\pU_{1,2}(z)\\
	-\pU_{2,1}(z) &\pU_{1,1}(z)\end{bmatrix}\begin{bmatrix}\pA_{1,1}(z)\\
	 \pA_{2,1}(z)\end{bmatrix}=dz^{n-k}\begin{bmatrix}\pU_{1,1}^\star(z)\\
	 -\pU_{1,2}^\star(z)\end{bmatrix}.\ee
Multiplying $\left[\pU_{1,1}, \pU_{1,2}\right]$ to the left on both sides of \er{pu:mat} and using \er{pa11} yields $\det(\pU(z))\pA_{1,1}(z)=dz^{n-k}\pA_{1,1}(z),$ and thus $\det(\pU(z))=dz^{n-k}.$

Multiplying $\left[\pU_{1,2}^\star, \pU_{1,1}^\star\right]$ to the left on both sides of \er{pu:mat} yields $$\left(\pU_{2,1}(z)\pU_{1,1}^\star(z)-\pU_{2,2}(z)\pU_{1,2}^\star(z)\right)\pA_{1,1}(z)=\pA_{1,1}(z)\pA_{2,1}(z).$$
Hence
\be\label{pa21}\pA_{2,1}(z)=\pU_{2,1}(z)\pU_{1,1}^\star(z)-\pU_{2,2}(z)\pU_{1,2}^\star(z).\ee
Multiplying $\left[\pU_{2,2}^\star, \pU_{2,1}^\star\right]$ to the left on both sides of \er{pu:mat}, together with using \er{pa21}, we have
\begin{align*}&\left(\pU_{2,2}(z)\pU_{2,2}^\star(z)-\pU_{2,1}(z)\pU_{2,1}^\star(z)\right)\pA_{1,1}(z)+\pA_{2,1}^\star(z)\pA_{2,1}(z)\\
=&dz^{n-c}\det(\pU(z))^\star=d^2=-\det(\pA(z))=\pA_{1,2}(z)\pA_{2,1}(z)-\pA_{1,1}(z)\pA_{2,2}(z).\end{align*}
Since $\pA$ is Hermitian, we have $\pA_{2,1}^\star=\pA_{1,2}$. So the above identities imply
Hence
\be\label{pa22}\pA_{2,2}(z)=\pU_{2,1}(z)\pU_{2,1}^\star(z)-\pU_{2,2}(z)\pU_{2,2}^\star(z).\ee

Consequently, it follows from \er{pa11},\er{pa21} and \er{pa22} that $\pA=\pU\DG(1,-1)\pU^\star$. Moreover, we have
$$\sym\pU(z)=\begin{bmatrix}\sym\pU_{1,1}(z) & \sym\pU_{1,2}(z)\\
\sym\pU_{2,1}(z) &\sym\pU_{2,2}(z)\end{bmatrix}=\begin{bmatrix}\eps z^{n} & -\eps z^n\\
-\eps z^{n-2k} &\eps z^{n-2k}\end{bmatrix}=\sym\pth_1^\star(z)\sym\pth_2(z),$$
where
$$\sym\pth_1(z)=[\eps z^{-n},-\eps z^{2k-n}]\quad\mbox{and}\quad \sym\pth_2(z)=[1, -1].$$
This show that $\pU$ has compatible symmetry and clearly \er{eq:Sym1} holds. Now the algorithm is fully justified.
\ep

\subsection{Factorization of matrices of Laurent polynomials with symmetry and coprime entries}

We now study the factorization as in Theorem~\ref{thm:SymSpecFactSym} for a $2\times 2$ Hermitian matrix $\pA$ of Laurent polynomials with compatible symmetry, without assuming that $\pA$ has constant determinant. In this subsection, we study the case when all entries of $\pA$ are \emph{coprime}, that is, $\gcd(\pA_{1,1},\pA_{1,2},\pA_{2,1},\pA_{2,2})=1$. The main result for this subsection is the following theorem.

\begin{theorem}\label{thm:SymSpecFactSym:gcd}Let $\pA$ be a $2\times 2$ Hermitian matrix of Laurent polynomials with compatible symmetry. Suppose that
	\begin{enumerate}
		 \item[(1)]$\gcd(\pA_{1,1},\pA_{1,2},\pA_{2,1},\pA_{2,2})=1$;
		
		 \item[(2)]$\det(\pA(z))=-\pd(z)\pd^\star(z)$ for some nonzero Laurent polynomial $\pd(z)$ with symmetry.
	\end{enumerate}
Then we can find $2\times 2$ matrices $\pV$ and $\pB$ of Laurent polynomials with compatible symmetry such that
\begin{enumerate}
	\item[(i)] $\pB=\pB^\star$ and $\det(\pB)=C$ is a negative constant;
	
	 \item[(ii)]$\pA=\pV\pB\pV^\star$ where all multiplications are compatible.
\end{enumerate}
As a consequence, there exists a $2\times2$ matrix $\pU$ of Laurent polynomials with compatible symmetry such that $\pA=\pU\DG(1,-1)\pU^\star$ and \er{eq:Sym1} holds.
\end{theorem}

The proof of Theorem~\ref{thm:SymSpecFactSym:gcd} relies heavily on long division of Laurent polynomials. Note that the symmetry structures are not preserved under long divisions. That is, if $\pa$ and $\pb$ are Laurent polynomials with symmetry and $\pb\ne \mathbf{0}$, then we can always do the long division
$$\pa(z)=\pb(z)\pq(z)+\pr(z),$$
for some Laurent polynomials $\pq$ and $\pr$. However, $\pq$ and $\pr$ may not necessarily have any symmetry. Fortunately, we can still establish the following Extended Euclidean Algorithm with symmetry.

\begin{theorem}\label{thm:eealg} \textbf{(Extended Euclidean Algorithm for Laurent polynomials with symmetry)} Let $\pa$ and $\pb$ be Laurent polynomials with symmetry, and let $\pr:=\gcd(\pa,\pb)$. Then there exist Laurent polynomials $\pu$ and $\pv$ with symmetry such that
	 \be\label{ld:lau}\pa(z)\pu(z)+\pb(z)\pv(z)=\pr(z).\ee
Moreover, we have $\sym\pa(z)\sym\pu(z)=\sym\pb(z)\sym\pv(z)=\sym\pr(z)$ and $\gcd(\pu,\pv)=1$.
\end{theorem}

\bp If $\pa=\mathbf{0}$, then we simply take $\pu=\mathbf{0}$ and $\pv=1$. If $\pb=\pzr$, then we take $\pu=1$ and $\pv=\mathbf{0}$.

Assume now $\pa\ne \mathbf{0}$ and $\pb\ne\mathbf{0}$. Let $\pr:=\gcd(\pa,\pb)$. By applying the classical Euclidean Algorithm, we can construct Laurent polynomials $\pu_1$ and $\pu_2$ such that
\be\label{ealg:ab}\pa(z)\pu_1(z)+\pb(z)\pu_2(z)=\pr(z).\ee
Thus
$\pa(z^{-1})\pu_1(z^{-1})+\pb(z^{-1})\pu_2(z^{-1})=\pr(z^{-1})$,
which can be rewritten as
\be\label{ealg:ab:1}\frac{\pa(z)\pu_1(z^{-1})}{\sym\pa(z)}+\frac{\pb(z)\pu_2(z^{-1})}{\sym\pb(z)}=\frac{\pr(z)}{\sym\pr(z)}.\ee
Combining \er{ealg:ab} and \er{ealg:ab:1}, we get \er{ld:lau} with
\be\label{ealg:uv}\pu(z)=\frac{\pu_1(z)+\frac{\pu_1(z^{-1})\sym\pr(z)}{\sym\pa(z)}}{2},\qquad \pv(z)=\frac{\pu_2(z)+\frac{\pu_2(z^{-1})\sym\pr(z)}{\sym\pb(z)}}{2}.\ee
We claim that the Laurent polynomials $\pu$ and $\pu$ defined as \er{ealg:uv} are the desired ones as required. In fact, note that $\pr=\gcd(\pa,\pb)$ has symmetry. Direct calculation yields $\sym\pu(z)=\frac{\sym\pr(z)}{\sym\pa(z)}$ and $\sym\pv(z)=\frac{\sym\pr(z)}{\sym\pb(z)}$. This implies that $\pu$ and $\pu$ have symmetry, and $\sym\pa(z)\sym\pu(z)=\sym\pb(z)\sym\pv(z)=\sym\pr(z)$. Let $\pd:=\gcd(\pu,\pv)$. Then $\pd\pr$ divides $\pa\pu+\pb\pv=\pr$. Therefore $\pd=1$. The proof is complete.\ep	

With the Extended Euclidean algorithm developed above, we have the following lemma.

\begin{lemma}\label{lem:u1un2}
Let $\pu_1$ and $\pu_2$ be nonzero Laurent polynomials with symmetry, and let $\pr:=\gcd(\pu_1,\pu_2)$. Then one can construct a strongly invertible $2\times 2$ matrix $\pP$ of Laurent polynomials with compatible symmetry such that
\be\label{pp:u1u2}\pP(z)\begin{bmatrix}\pu_1(z)\\
	 \pu_2(z)\end{bmatrix}=\begin{bmatrix}\pr(z)\\
	0\end{bmatrix},\ee
and the above matrix multiplication is compatible.
\end{lemma}

\bp By Theorem~\ref{thm:eealg}, we can find coprime Laurent polynomials $\pu$ and $\pv$ with symmetry, such that
\be\label{uv:u1u2}\pu_1(z)\pu(z)+\pu_2(z)\pv(z)=\pr(z),\ee
and
\be\label{uv:u1u2:1}\sym\pu_1(z)\sym\pu(z)=\sym\pu_2(z)\sym\pv(z)=\sym\pr(z).\ee
As $\pu$ and $\pv$ are coprime, applying Theorem~\ref{thm:eealg} again yields
\be\label{uv:st}\pu(z)\ps(z)-\pv(z)\pt(z)=1,\ee
for some coprime Laurent polynomials $\ps$ and $\pt$ with symmetry satisfying
\be\label{uv:st:1}\sym\pu(z)\sym\ps(z)=\sym\pv(z)\sym\pt(z)=1.\ee
Define $\pP_1:=\begin{bmatrix}\pu & \pv\\
\pt &\ps\end{bmatrix}$. Then $\det(\pP_1)=1$, which means $\pP_1$ is strongly invertible. Moreover, it follows from \er{uv:u1u2} that
\be\label{pp1}\pP_1(z)\begin{bmatrix}\pu_1(z)\\
\pu_2(z)\end{bmatrix}=\begin{bmatrix}\pr(z)\\
\pt(z)\pu_1(z)+\ps(z)\pu_2(z)\end{bmatrix}.\ee
Using \er{uv:u1u2:1} and \er{uv:st:1}, we see that $\frac{\sym\pu_1}{\sym\pu_2}=\frac{\sym\pv}{\sym\pu}=\frac{\sym\pt}{\sym\ps}=\eps z^c$ for some $\eps\in\{\pm1\}$ and $c\in\Z$. This implies that
$$\sym\pP_1(z)=\begin{bmatrix}\sym\pu(z) & \sym\pv(z)\\
\sym\pt(z) &\sym\ps(z)\end{bmatrix}=\sym\pth_1^\star(z)\sym\pth_2(z),\quad \begin{bmatrix}\sym\pu_1(z)\\
\sym\pu_2(z)\end{bmatrix}=\sym\pth_2^\star(z)\sym\pth_3(z),$$
where
$$\sym\pth_1(z)=[\sym\pu^\star(z),\sym\pt^\star(z)],\quad \sym\pth_2(z)=[1,\eps z^{c}],\quad \sym\pth_3(z)=\sym\pu_1(z).$$
Hence $\pP_1$ has compatible symmetry and all multiplications in \er{pp1} are compatible.

Define $\pP_2:=\begin{bmatrix}1 & 0\\
-\pp & 1\end{bmatrix}$, where $\pp:=\tfrac{\pt\pu_1+\ps\pu_2}{\pr}$ is a well-defined Laurent polynomial as $\pr=\gcd(\pu_1,\pu_2)$ definitely divides $\pt\pu_1+\ps\pu_2$. We see that $\pP_2$ is strongly invertible as $\det(\pP_2)=1$. Furthermore, we have
\be\label{pp2}\pP_2(z)\begin{bmatrix}\pr(z)\\
\pt(z)\pu_1(z)+\ps(z)\pu_2(z)\end{bmatrix}=\begin{bmatrix}\pr(z)\\
0\end{bmatrix}.\ee
Applying a similar argument as proving that $\pP_1$ has compatible symmetry, we see that $\pP_2$ has compatible symmetry and all multiplications in \er{pp2} are compatible. By letting $\pP:=\pP_2\pP_1$, we see that $\pP$ is strongly invertible and \er{pp:u1u2} holds. Moreover, $\pP$ has compatible symmetry and all multiplications in \er{pp:u1u2} are compatible. The proof is now complete.\ep

For a $2\times 2$ matrix $\pA$ of Laurent polynomial, recall its Smith normal form:
\be\label{snf}\pA(z)=\pE(z)\DG(\pd_1(z),\pd_2(z))\pF(z),\ee
where $\pE$ and  $\pF$ are $2\times 2$ strongly invertible matrices of Laurent polynomials, $\pd_1$ and $\pd_2$ are monic polynomials such that $\pd_1$ divides $\pd_2$. From the theory of Smith normal form, $\pd_j$ is uniquely determined up to a multiplication by a monomial for $j=1,2$. Thus without loss of generality, we require the polynomials $\pd_1$ and $\pd_2$ in the factorization \er{snf} to have nonzero constant terms. In this case, we call $\pd_1$ and $\pd_2$ the \emph{invariant polynomials} of $\pA$.

A factorization as in \er{snf} is obtained by iteratively using long divisions, which might not preserve symmetry structures. As a result, even if $\pA$ has symmetry, there is not much to say about the symmetry information on $\pE$, $\pF$, $\pd_1$ and $\pd_2$. To resolve this issue, we now introduce a normal form of a $2\times2$ matrix $\pA$ of Laurent polynomials with compatible symmetry, which will be a backbone of the proof of Theorem~\ref{thm:SymSpecFactSym:gcd}.

\begin{theorem}\label{thm:nf:sym} \textbf{(Normal form with compatible symmetry)}
Let $\pA$ be a $2\times 2$ matrix of Laurent polynomials with compatible symmetry. Then there exist strongly invertible $2\times 2$ matrices $\pP$ and $\pQ$ of Laurent polynomials with compatible symmetry such that $\pP\pA\pQ$ is a diagonal matrix:
\be\label{dmat:nf}\pD(z)=\DG(\pe_1(z),\pe_2(z)):=\pP(z)\pA(z)\pQ(z),\ee
where all multiplications above are compatible, and both $\pe_1,\pe_2$ have symmetry. Furthermore, if $\pd_1$ and $\pd_2$ are the invariant polynomials of $\pA$, then $\{\mz(\pe_1,z_0), \mz(\pe_2,z_0)\}=\{\mz(\pd_1,z_0), \mz(\pd_2,z_0)\}$ for all $z_0\in\C\setminus\{0\}$.
\end{theorem}

\bp If $\pA$ is the zero-matrix, then we simply take $\pP=\pQ=\pI_2$ and $\pD$ the zero matrix will work. If $\pA$ is nonzero and is diagonal, then we take $\pP=\pQ=\pI_2$ and $\pD=\pA$.

For the other non-trivial cases, some off-diagonal elements of $\pA$ is nonzero. Here we assume that $\pA_{2,1}\ne\mathbf{0}$, and the case for $\pA_{1,2}\ne \mathbf{0}$ can be dealt with similarly. Set $\pA^{(0,0)}:=\pA$. For $j=0,1,\dots$, perform the following iterative process:
\begin{enumerate}
\item[Step (i):] If $\pA^{(j,j)}_{2,1}\ne\mathbf{0}$, construct a strongly invertible $2\times 2$ matrix $\pP^{(j)}$ of Laurent polynomials with compatible symmstry such that
\be\label{a:j1:j}\pA^{(j+1,j)}(z):=\pP^{(j)}(z)\pA^{(j,j)}(z)=\begin{bmatrix}\pr^{(j,j)}(z) &\pA^{(j+1,j)}_{1,2}(z) \\
0	& \pA^{(j+1,j)}_{2,2}(z)\end{bmatrix},\ee
where $\pr^{(j,j)}:=\gcd(\pA^{(j,j)}_{1,1},\pA^{(j,j)}_{2,1})\ne \pzr$. In fact, if $\pA^{(j,j)}_{1,1}=\mathbf{0}$, simply take $\pP^{(j)}=\begin{bmatrix}0 &1\\
1 &0\end{bmatrix}$. Otherwise, apply Lemma~\ref{lem:u1un2} to find a strongly invertible matrix $\pP^{(j)}$ of Laurent polynomials with compatible symmetry, such that
$$\pP^{(j)}(z)\begin{bmatrix}\pA^{(j,j)}_{1,1}(z)\\
\pA^{(j,j)}_{2,1}(z)\end{bmatrix}=\begin{bmatrix}\pr^{(j,j)}(z)\\
0\end{bmatrix}.$$
Then $\pP^{(j)}$ is the desired matrix as required. Moreover, the multiplication in \er{a:j1:j} is compatible.

\item[Step (ii):] If $\pA^{(j+1,j)}$ is a diagonal matrix, i.e., $\pA^{(j+1,j)}_{1,2}=\mathbf{0}$, we terminate the process. Otherwise, define $\pr^{(j+1,j)}:=\gcd(\pr^{(j,j)},\pA^{(j+1,j)}_{1,2})\ne \pzr$. Apply Lemma~\ref{lem:u1un2} to find a strongly invertible matrix $\pQ^{(j)}$ of Laurent polynomials with compatible symmetry, such that
$$\pQ^{(j)\star}(z)\begin{bmatrix}\pA^{(j+1,j)\star}_{1,1}(z)\\
\pA^{(j+1,j)\star}_{1,2}(z)\end{bmatrix}=\pQ^{(j)\star}(z)\begin{bmatrix}\pr^{(j,j)\star}(z)\\
\pA^{(j+1,j)\star}_{1,2}(z)\end{bmatrix}=\begin{bmatrix}\pr^{(j+1,j)\star}(z)\\
0\end{bmatrix},$$
 It follows that
\be\label{a:j1:j1}\pA^{(j+1,j+1)}(z):=\pA^{(j+1,j)}(z)\pQ^{(j)}(z)=\begin{bmatrix}\pr^{(j+1,j)}(z) &0 \\
\pA^{(j+1,j+1)}_{2,1}(z)	& \pA^{(j+1,j+1)}_{2,2}(z)\end{bmatrix}.\ee
Moreover, item (3) of Proposition~\ref{sym:mat:com} tells that the multiplication above is compatible.

\item[Step (iii)]: If $\pA^{(j+1,j+1)}$ is a diagonal matrix, i.e., $\pA^{(j+1,j+1)}_{2,1}=\mathbf{0}$, we terminate the process. Otherwise, redefine $j:=j+1$ and go to step (i).

\end{enumerate}

We claim that the above iterative process cannot go on forever, i.e. $\pA^{(j,j)}$ or $\pA^{(j+1,j)}$ becomes diagonal for some $j\in\N$. For every $j\in\N$, we see that $\pA_{1,1}^{(j+1,j)}=\pr^{(j,j)}=\gcd(\pA_{1,1}^{(j,j)},\pA_{2,1}^{(j,j)})$, which divides $\pA_{1,1}^{(j,j)}$, and similarly $\pA_{1,1}^{(j+1,j+1)}$ divides $\pA_{1,1}^{(j+1,j)}$.  This implies that $\len(\pA_{1,1}^{(j+1,j+1)})\le\len(\pA_{1,1}^{(j+1,j)})\le\len(\pA_{1,1}^{(j,j)})$. If $\pA_{1,1}^{(j,j)}$ does not divide $\pA_{2,1}^{(j,j)}$, we have $\len(\pA_{1,1}^{(j+1,j)})<\len(\pA_{1,1}^{(j,j)})$. Similarly if $\pA_{1,1}^{(j+1,j)}$ does not divide $\pA_{1,2}^{(j+1,j)}$, we have $\len(\pA_{1,1}^{(j+1,j+1)})<\len(\pA_{1,1}^{(j+1,j)})$. As the length of the $(1,1)$-th entry cannot decrease strictly forever, at some point we must have $\pA_{1,1}^{(j,j)}$ divides $\pA_{2,1}^{(j,j)}$, or $\pA_{1,1}^{(j+1,j)}$ divides $\pA_{1,2}^{(j+1,j)}$. If $\pA_{1,1}^{(j,j)}$ divides $\pA_{2,1}^{(j,j)}$, then choose
$$\pP^{(j)}(z):=\begin{bmatrix}1 & 0\\
-\frac{\pA^{(j,j)}_{2,1}(z)}{\pA^{(j,j)}_{1,1}(z)} & 1\end{bmatrix}.$$
Then $\pP^{(j)}$ has compatible symmetry, \er{a:j1:j} holds and all multiplications are compatible. Moreover, the matrix $\pA^{(j+1,j)}$ defined as  \er{a:j1:j} is diagonal. If $\pA_{1,1}^{(j+1,j)}$ divides $\pA_{1,2}^{(j+1,j)}$, then choose
$$\pQ^{(j)}(z):=\begin{bmatrix}1 & -\frac{\pA^{(j+1,j)}_{1,2}(z)}{\pA^{(j+1,j)}_{1,1}(z)}\\
0 & 1\end{bmatrix}.$$
Then $\pQ^{(j)}$ has compatible symmetry, \er{a:j1:j1} holds and all multiplications are compatible. Moreover, the matrix $\pA^{(j+1,j+1)}$ defined as  \er{a:j1:j} is diagonal.

Therefore, by applying steps (i)-(iii) for finitely many times, we get
strongly invertible matrices $\pP$ and $\pQ$ of Laurent polynomials with compatible symmetry, such that $\pD:=\pP\pA\pQ$ is a diagonal matrix satisfying \er{dmat:nf}, and all multiplications in \er{dmat:nf} are compatible. The symmetry of the diagonal elements $\pe_1$ and $\pe_2$  of $\pD$ follows trivially from our previous construction steps.

Now let $\pd_1$ and $\pd_2$ be invariant polynomials of $\pA$. Then there exist strongly invertible $2\times 2$ matrices $\pE$ and $\pF$ such that \er{snf} holds. By the definition of invariant polynomials, we have $\mz(\pd_1,z_0)\le\mz(\pd_2,z_0)$. Without loss of generality, we assume $\mz(\pe_1,z_0)\le\mz(\pe_2,z_0)$, and prove that $\mz(\pd_j,z_0)=\mz(\pe_j,z_0)$ for $j=1,2$.

For any $z_0\in\C\setminus\{0\}$, write $\pd_j(z)=(z-z_0)^{\mz(\pd_j,z_0)}\pp_j(z)$ and $\pe_j(z)=(z-z_0)^{\mz(\pe_j,z_0)}\pq_j(z)$ for $j=1,2$. It follows that
\be\label{pd:pv}\pA(z)=\pE(z)\DG\left((z-z_0)^{\mz(\pd_1,z_0)},(z-z_0)^{\mz(\pd_2,z_0)}\right)\DG(\pp_1(z),\pp_2(z))\pF(z),\ee
\be\label{pe:pv}\pA(z)=\pP(z)^{-1}\DG\left((z-z_0)^{\mz(\pe_1,z_0)},(z-z_0)^{\mz(\pe_2,z_0)}\right)\DG(\pq_1(z),\pq_2(z))\pQ^{-1}(z).\ee
It follows that
\be\label{pd:pe}\DG\left((z-z_0)^{\mz(\pd_1,z_0)},(z-z_0)^{\mz(\pd_2,z_0)}\right)=\pW(z)\DG\left((z-z_0)^{\mz(\pe_1,z_0)},(z-z_0)^{\mz(\pe_2,z_0)}\right)\pV(z),\ee
where $\pW:=\pE^{-1}\pP^{-1}$ and $\pV:=\DG(\pq_1,\pq_2)\pQ^{-1}\pF^{-1}\DG(\pp_1,\pp_2)^{-1}$. It is trivial that all entries of $\pW$ and $\pV$, as well as $\det(\pW)$ and $\det(\pV)$, are analytic on a neighborhood of $z_0$. Moreover, $\det(\pW)$ and $\det(\pV)$ and are nonzero at $z_0$. By calculating determinants on both sides of \er{pd:pe}, we have
$$(z-z_0)^{\mz(\pd_1,z_0)+\mz(\pd_2,z_0)}=(z-z_0)^{\mz(\pe_1,z_0)+\mz(\pe_2,z_0)}\det(\pW(z))\det(\pV(z)).$$
Note that $\det(\pW)$ and $\det(\pV)$ are analytic in a neighborhood of $z_0$, and are nonzero at $z_0$. Thus
\be\label{mz:pd:pe}\mz(\pd_1,z_0)+\mz(\pd_2,z_0)=\mz(\pe_1,z_0)+\mz(\pe_2,z_0).\ee
Furthermore, \er{pd:pe} yields
$$(z-z_0)^{\mz(\pd_1,z_0)}=(z-z_0)^{\mz(\pe_1,z_0)}\left[\pW_{1,1}(z)\pV_{1,1}(z)+(z-z_0)^{\mz(\pe_2,z_0)-\mz(\pe_1,z_0)}\pW_{1,2}(z)\pV_{2,1}(z)\right],$$
which implies $\mz(\pe_1,z_0)\le\mz(\pd_1,z_0)$. On the other hand, \er{pd:pe} is equivalent to
$$\pW(z)^{-1}\DG\left((z-z_0)^{\mz(\pd_1,z_0)},(z-z_0)^{\mz(\pd_2,z_0)}\right)\pV^{-1}(z)=\DG\left((z-z_0)^{\mz(\pe_1,z_0)},(z-z_0)^{\mz(\pe_2,z_0)}\right),$$
and by using a similar argument we conclude that $\mz(\pd_1,z_0)\le\mz(\pe_1,z_0)$. Therefore, $\mz(\pd_1,z_0)=\mz(\pe_1,z_0)$, and thus by \er{mz:pd:pe} we have $\mz(\pd_2,z_0)=\mz(\pe_2,z_0)$. This completes the proof.\ep

For every square matrix $\pA$ of Laurent polynomials, we define its \emph{spectrum} as
\be\label{spec}\sigma(\pA):=\{z\in\C\setminus\{0\}: \det(\pA(z))=0\}.\ee
When $\det(\pA)$ is non-constant, it is clear that $\sigma(\pA)$ is non-empty. 
The following lemma demonstrates how we can shrink the length of $\det(\pA)$ as well as the set $\sigma(\pA)$.

\begin{lemma}\label{lem:spec}Let $\pA$ be a $2\times 2$ Hermitian matrix of Laurent polynomials with compatible symmetry. Suppose that
	\begin{enumerate}
		 \item[(1)]$\gcd(\pA_{1,1},\pA_{1,2},\pA_{2,1},\pA_{2,2})=1$;
		
		 \item[(2)]$\det(\pA)=-\pd \pd^\star$ for some nonzero Laurent polynomial $\pd(z)$ with symmetry.
	\end{enumerate}
Then we can construct $2\times 2$ matrices $\pU$ and $\ptA$ of Laurent polynomials with compatible symmetry such that

\begin{enumerate}
	\item[(i)] $\ptA$ is Hermitian and $\pA=\pU\ptA\pU^\star$ with all multiplications being compatible;
	
	 \item[(ii)]$\gcd(\ptA_{1,1},\ptA_{1,2},\ptA_{2,1},\ptA_{2,2})=1$;
		
	 \item[(iii)]$\det(\ptA)=-\ptd\ptd^\star$ for some nonzero Laurent polynomial $\ptd$ with symmetry;
	
	\item[(iv)] $\len(\det(\ptA))<\len(\det(\pA))$ and $\sigma(\ptA)\subsetneq\sigma(\pA)$.
	\end{enumerate}
	
\end{lemma}

\bp By Theorem~\ref{thm:nf:sym}, one can construct strongly invertible matrices $\pP$ and $\pQ$ of Laurent polynomials with compatible symmetry, such that \er{dmat:nf} holds, where all multiplications are compatible, and $\pe_1,\pe_2$ have symmetry. Without loss of generality, we may assume that $\det(\pP)=\det(\pQ)=1$, in which case we have
\be\label{det:pe:pd}\pe_1(z)\pe_2(z)=\det(\pD(z))=\det(\pA(z))=-\pd(z)\pd^\star(z).\ee
Furthermore, item (1) implies that the Smith normal form of $\pA$ is $\DG(1,\det(\pA))$. Thus by Theorem~\ref{thm:nf:sym}, we have $\mz(\pe_j,z_0)=0$ or $ \mz(\pe_j,z_0)=\mz(\det(\pA),z_0)$ for all $z_0\in\C\setminus\{0\}$ and $j=1,2$. Define the Hermitian matrix $\mrpA$ of Laurent polynomials via
\be\label{mrpA}\mrpA(z):=\pP^{-1}\pA(z)\pP^{-\star}(z)=\DG(\pe_1(z),\pe_2(z))\pQ(z)\pP^{-\star}(z).\ee
It is trivial that all multiplications above are compatible. Let $z_0\in\sigma(\pA)$. by \er{det:pe:pd}, we must have $(z-z_0)$ divides $\pe_k$ for some $k\in\{1,2\}$. So $\alpha:=\mz(\pe_k,z_0)>0$ and our previous discussion yields $\alpha=\mz(\det(\pA),z_0)$. Moreover, we have $\mz(\pe_k,z_0^{-1})=\alpha$ by the symmetry property of $\pe_k$. There are three possibilities for the location of $z_0$.

\begin{itemize}
	\item \textbf{If $z_0\in\C\setminus(\T\cup\R)$: }Then $z_0, z_0^{-1}, \ol{z_0}$ and $\ol{z_0}^{-1}$ are pairwise distinct. Define the Laurent polynomial $\pp(z):=(z-z_0)^\alpha(z-z_0^{-1})^\alpha$ with symmetry. Then $\pp$ divides $\pe_k$, and thus divides the $k$-th row of $\mrpA$. As $\mrpA$ is Hermitian, we see that $\pp^\star$ divides the $k$-th column of $\mrpA$. Note that $\pp^\star(z)=z^{-2\alpha}(z-\ol{z_0})^\alpha (z-\ol{z_0}^{-1})^\alpha$, which implies $\gcd(\pp,\pp^\star)=1$. Hence $\pp\pp^\star$ divides $\mrpA_{k,k}$. Define
	 \be\label{pvp}\pV_{1,\pp}(z):=\DG(\pp(z),1),\qquad \pV_{2,\pp}(z):=\DG(1,\pp(z)).\ee
Then we have
\be\label{pvp:2}\mrpA(z)=\pV_{k,\pp}(z)\ptA(z)\pV_{k,\pp}^\star(z)\ee
for some $2\times2$ matrix $\ptA$ of Laurent polynomials.

By letting $\pU:=\pP\pV_{k,\pp}$, it follows immediately that $\pA=\pU\ptA\pU^\star$. We claim that the matrices $\pU$ and $\ptA$ satisfy all claims of the lemma. Since $\pp$ has symmetry and $\mrpA$ has compatible symmetry, we conclude that $\pU$ and $\ptA$ have compatible symmetry, and item (i) clearly holds.

Since $\det(\pP)=\det(\pQ)=1$, so the Smith normal form of $\mrpA$ is the same as the one of $\pA$, which is $\DG(1,\det(\pA))$. This in particular implies $\gcd(\mrpA_{1,1},\mrpA_{1,2},\mrpA_{2,1},\mrpA_{2,2})=1$. Therefore, by letting $l\in\{1,2\}\setminus\{k\}$, we have
$$\gcd(\ptA_{1,1},\ptA_{1,2},\ptA_{2,1},\ptA_{2,2})=\gcd\left(\frac{\mrpA_{k,k}}{\pp\pp^\star},\frac{\mrpA_{k,l}}{\pp},\frac{\mrpA_{l,k}}{\pp^\star},\mrpA_{l,l}\right)=1.$$
This proves item (ii).

By $\det(\pA)=\det(\mrpA)=\det(\pV_{k,\pp})\det(\ptA)\det(\pV_{k,\pp}^\star)=\pp\pp^\star\det(\ptA)$, we have $\det(\ptA)=\frac{\det(\pA)}{\pp\pp^\star}=\frac{-\pd\pd^\star}{\pp\pp^\star}$. As $\pd$ has symmetry, we can define
$$\beta_1:=\mz(\pd,z_0)=\mz(\pd,z_0^{-1}),\quad \beta_2:=\mz(\pd,\ol{z_0})=\mz(\pd,\ol{z_0}^{-1}),$$
$$\pq(z):=(z-z_0)^{\beta_1}(z-z_0^{-1})^{\beta_1}(z^{-1}-\ol{z_0})^{\beta_2}(z^{-1}-\ol{z_0}^{-1})^{\beta_2}.$$
We have $\alpha=\mz(\det(\pA),z_0)=\mz(\pd,z_0)+\mz(\pd^\star,z_0)=\beta_1+\beta_2$. Hence
$\pq\pq^{\star}=\pp\pp^\star$. Moreover, $\pq$ has symmetry and divides $\pd$. Define $\ptd:=\frac{\pd}{\pq}$. It is trivial that $\ptd$ is a Laurent polynomial with symmetry, and we have $\det(\ptA)=-\ptd\ptd^\star$. This proves item (iii).

It is easy to see that $\len(\pq)>0$, thus $\len(\ptd)<\len(\pd)$ and $\len(\det(\ptA))<\len(\det(\pA))$. Furthermore, it is trivial that $\sigma(\ptA)\subseteq \sigma(\pA)$. Note that
$$\mz(\det(\ptA),z_0)=\mz\left(\frac{\det(\pA)}{\pp\pp^\star},z_0\right)=\mz(\det(\pA),z_0)-\mz(\pp,z_0)-\mz(\pp^\star,z_0)=\alpha-\alpha-0=0,$$
which implies $z_0\notin\sigma(\ptA)$. Hence $\sigma(\ptA)\subsetneq \sigma(\pA)$, and item (iv) is verified.

\item \textbf{If $z_0\in(\T\cup\R)\setminus\{0,\pm 1\}$: } Then $\ol{z_0}^{-1}\ne z_0$. When $z_0\in\T\setminus\{\pm 1\}$, we have $\ol{z_0}^{-1}=z_0$, so $\mz(\pd,z_0)=\mz(\pd,\ol{z_0}^{-1})=\mz(\pd^\star,z_0)$. When $z_0\in\R\setminus\{0,\pm 1\}$, we have $\ol{z_0}=z_0$, so $ \mz(\pd,z_0)=\mz(\pd,z_0^{-1})=\mz(\pd,\ol{z_0}^{-1})=\mz(\pd^\star,z_0)$. Thus we conclude that $\mz(\pd,z_0)=\mz(\pd^\star,z_0)$ for this case. It follows that
\be\label{alpha:mz}\alpha=\mz(\det(\pA),z_0)=\mz(\pd,z_0)+\mz(\pd^\star,z_0)=2\mz(\pd,z_0)\in2\Z.\ee
Define
$$\pp(z):=z^{-\alpha/2}(z-z_0)^{\alpha/2}(z-z_0^{-1})^{\alpha/2}=\begin{cases}(z+z^{-1}-2\re(z_0))^{\alpha/2},&\text{ if }z_0\in\T\setminus\{\pm 1\},\\
(z+z^{-1}-(z_0+z_0^{-1}))^{\alpha/2},&\text{ if }z_0\in\R\setminus\{0,\pm 1\}.\end{cases}$$
It is easy to check that $\sym\pp=1$ and $\pp^\star=\pp$. Given that $z_0\ne\pm1$, one can conclude that $(z-z_0)^\alpha(z-z_0^{-1})^\alpha$ divides $\pe_k$, i.e., $\pp\pp^\star$ divides $\pe_k$. Hence by \er{det:pe:pd}, $\pp$ divides the $k$-th row of $\mrpA$, $\pp^\star$ divides the $k$-th column of $\mrpA$, and $\pp\pp^\star$ divides $\mrpA_{k,k}$. Define $\pV_{j,\pp}$ as \er{pvp} for $j=1,2$, and define $\pU:=\pP\pV_{k,\pp}$. We see that \er{pvp:2} holds for some $2\times2$ matrix $\ptA$ of Laurent polynomials. Moreover, $\pU$ and $\ptA$ have compatible symmetry, and items (i) and (ii) hold.

From \er{alpha:mz}, we see that $\mz(\pd,z_0)=\mz(\pd,z_0^{-1})=\alpha/2$. As $\pd$ has symmetry, so $\pp$ divides $\pd$. Define $\ptd:=\pd/\pp$, then $\ptd$ has symmetry and $\det(\ptA)=-\ptd\ptd^\star$. This proves item (iii).

As $\len(\pp)>0$, we have $\len(\ptd)<\len(\pd)$ and thus $\len(\det(\ptA))<\len(\det(\pA))$. Furthermore, it is trivial that $\sigma(\ptA)\subseteq\sigma(\pA)$. Note that
$$\mz(\det(\ptA),z_0)=\mz(\det(\pA),z_0)-\mz(\pp,z_0)-\mz(\pp^\star,z_0)=\alpha-\alpha/2-\alpha/2=0,$$
which means $z_0\notin\sigma(\ptA)$. Thus $\sigma(\ptA)\subsetneq\sigma(\pA)$, and item (iv) is proved.

\item \textbf{If $z_0\in\{\pm 1\}$: } Then $z_0=z_0^{-1}=\ol{z_0}=\ol{z_0}^{-1}$, which implies $\mz(\pd,z_0)=\mz(\pd^\star,\ol{z_0}^{-1})=\mz(\pd^\star,z_0)$ and thus \er{alpha:mz} holds. Define $\pp(z):=(z-z_0)^{\alpha/2}$. Direct calculation yields $\sym\pp(z)=(-z_0)^{\alpha/2} z^{\alpha/2}$ and $\pp^\star(z)=(-z_0z)^{-\alpha/2}(z-z_0)^{\alpha/2}$. Moreover, $\pp(z)\pp^\star(z)=(-z_0z)^{-\alpha/2}(z-z_0)^{\alpha}$ divides $\pe_k(z)$. By \er{det:pe:pd}, $\pp$ divides the $k$-th row of $\mrpA$, $\pp^\star$ divides the $k$-th column of $\mrpA$, and $\pp\pp^\star$ divides $\mrpA_{k,k}$. Define $\pV_{j,\pp}$ as \er{pvp} for $j=1,2$, and define $\pU:=\pP\pV_{k,\pp}$. We see that \er{pvp:2} holds for some $2\times2$ matrix $\ptA$ of Laurent polynomials. Moreover, $\pU$ and $\ptA$ have compatible symmetry, and by using the same arguments as in the previous case, we can prove that items (i) - (iv) hold.
\end{itemize}

The proof is now complete.\ep

\bp[\textbf{Proof of Theorem~\ref{thm:SymSpecFactSym:gcd}.}]If $\det(\pA)=C$ is a negative constant, then taking $\pB:=\pA$ and $\pV:=\pI_2$ yields the result.

Assume $\det(\pA)$ is not a constant. Set $\pA^{(0)}:=\pA$. By Lemma~\ref{lem:spec}, we can find two sequences $(\pU_j)_{j=1}^K$ and $(\pA^{(j)})_{j=1}^K$ of Laurent polynomials with compatible symmetry, and a sequence $(\pd_j)_{j=1}^K$ of nonzero Laurent polynomials with symmetry, such that

\begin{itemize}
\item $\pA^{(j)}$ is Hermitian and $\pA^{(j-1)}=\pU_j\pA^{(j)}\pU_j^{\star}$ with all multiplications being compatible for $j=1,\dots,K$;

\item $\gcd(\pA_{1,1}^{(j)},\pA_{1,2}^{(j)},\pA_{2,1}^{(j)},\pA_{2,2}^{(j)})=1$ for all $j=1,\dots,K$;

\item $\det(\pA^{(j)})=-\pd^{(j)}\pd^{(j)\star}$ for all $j=0,\dots,K$

\item $\len(\det(\pA^{(0)}))>\len(\det(\pA^{(1)}))>\dots>\len(\det(\pA^{(K)}))=0$ and $\sigma(\pA^{(0)})\supsetneq\sigma(\pA^{(1)})\supsetneq\dots\supsetneq\sigma(\pA^{(K-1)})=\emptyset$.

\end{itemize}

As a consequence, we have $\pA=\pA^{(0)}=\pU_1
\cdots\pU_K\pA^{(K)}\pU_K^\star\cdots
\pU_1.$
Define $\pB:=\pA^{(K)}$ and $\pU:=\pU_1\cdots\pU_K$. We see that both $\pB$ and $\pU$ have compatible symmetry and all multiplications in $\pA=\pU\pB\pU^\star$ are compatible. Note that $\det(\pB)=\det(\pA^{(K)})=-\pd\pd^\star=\det(\pB)^\star$. Thus by $\len(\det(\pB))=\len(\det(\pA^{(K)}))=0$, we conclude that $\det(\pB)$ is a negative constant. This proves item (i). Finally, item (ii) is a direct consequence of Theorem~\ref{thm:SymSpecFactSym:sp}. This completes the proof.\ep

\subsection{Difference of (Hermitian) squares property}
\label{subsec:dos}

In this subsection, we prove Theorem~\ref{thm:dos} for the difference of squares (DOS) property of a Laurent polynomial  as stated in Section~\ref{sec:introduction}.
Writing a Laurent polynomial into sums and differences of squares is inevitable for a generalized spectral factorization problem. Therefore it is important to deepen our understandings of the DOS property.

We need two lemmas to prove Theorem~\ref{thm:dos}. The following lemma demonstrates that the DOS property is preserved under multiplication.

\begin{lemma}\label{lem:dos:0}Suppose $\pu_1,\pu_2,\pu_3$ and $\pu_4$ are Laurent polynomials with symmetry, denote the symmetry types by
$$\sym\pu_j(z)=\eps_j z^{c_j},\qquad \eps_j\in\{\pm 1\},\quad c_j\in\Z,\quad j=1,2,3,4.$$
Suppose $\frac{\sym\pu_1}{\sym\pu_2}=\frac{\sym\pu_3}{\sym\pu_4}$. Define
\be\label{pu5:pu6}\pu_5(z):=\pu_1(z)\pu_3(z)+z^{c_1}\pu_2^\star(z)\pu_4(z),\quad \pu_6(z):=\pu_2(z)\pu_3(z)+z^{c_1}\pu_1^\star(z)\pu_4(z).\ee
Then
\be\label{dos:sym:1}\sym\pu_5(z)=\eps_1\eps_3z^{c_1+c_3},\quad \sym\pu_6(z)=\eps_2\eps_3z^{c_2+c_3},\quad \frac{\sym\pu_5}{\sym\pu_6}=\frac{\sym\pu_1}{\sym\pu_2}=\frac{\sym\pu_3}{\sym\pu_4},\ee
and
\be\label{dos:sym:2}\pu_5\pu_5^\star-\pu_6\pu_6^\star=(\pu_1\pu_1^\star-\pu_2\pu_2^\star)(\pu_3\pu_3^\star-\pu_4\pu_4^\star).\ee
\end{lemma}

\bp By $\frac{\sym\pu_1}{\sym\pu_2}=\frac{\sym\pu_3}{\sym\pu_4}$, we have $\eps_1\eps_2=\eps_3\eps_4$ and $c_1-c_2=c_3-c_4$. It follows that
$$\sym[z^{c_1}\pu_2^\star\pu_4](z)=\eps_2\eps_4 z^{2c_1}z^{-c_2+c_4}=\eps_1\eps_3z^{c_1+c_3}=\sym[\pu_1^\star\pu_3](z),$$
$$\sym[z^{c_1}\pu_1^\star\pu_4](z)=\eps_1\eps_4 z^{2c_1}z^{-c_1+c_4}=\eps_2\eps_3z^{c_2+c_3}=\sym[\pu_2^\star\pu_3](z).$$
Thus \er{dos:sym:1} follows immediately.

Note that \er{pu5:pu6} is equivalent to
$$\begin{bmatrix}\pu_5\\
\pu_6\end{bmatrix}=\begin{bmatrix}\pu_1 & \pu_2^\star\\
\pu_2 & \pu_1^\star\end{bmatrix}\begin{bmatrix}\pu_3\\
z^{c_1}\pu_4\end{bmatrix}.$$
It follows that
\begin{align*}
\pu_5\pu_5^\star-\pu_6\pu_6^\star&=\begin{bmatrix}\pu_5^\star&
\pu_6^\star\end{bmatrix}\DG(1,-1)\begin{bmatrix}\pu_5\\
\pu_6\end{bmatrix}\\
&=\begin{bmatrix}\pu_3^\star&
z^{-c_1}\pu_4^\star\end{bmatrix}\begin{bmatrix}\pu_1^\star & \pu_2^\star\\
\pu_2 & \pu_1\end{bmatrix}\DG(1,-1)\begin{bmatrix}\pu_1 & \pu_2^\star\\
\pu_2 & \pu_1^\star\end{bmatrix}\begin{bmatrix}\pu_3\\
z^{c_1}\pu_4\end{bmatrix}\\
&=(\pu_1\pu_1^\star-\pu_2\pu_2^\star)\begin{bmatrix}\pu_3^\star&
z^{-c_1}\pu_4^\star\end{bmatrix}\DG(1,-1)\begin{bmatrix}\pu_3\\
z^{c_1}\pu_4\end{bmatrix}\\
&=(\pu_1\pu_1^\star-\pu_2\pu_2^\star)(\pu_3\pu_3^\star-\pu_4\pu_4^\star).
\end{align*}
This completes the proof.\ep

\begin{lemma}\label{lem:dos:1}Let $\pu$ be a Laurent polynomial such that $\pu^\star=\pu$. Define
\be\label{sigma:pu}\sigma(\pu):=\{z\in\C\setminus\{0\}:\pu(z)=0\}.\ee
Then $\sum_{z_0\in\T\cap \sigma(\pu)}\mz(\pu,z_0)\in2\Z.$
\end{lemma}

\bp Write $\pu(z)=\sum_{k\in\Z}u(k)z^k$. By $\pu^\star=\pu$, we have $u(k)=\ol{u(-k)}$ for all $k\in\Z$. So $\fsupp(\pu)$ is a symmetric interval with centre $0$ and $\len(\pu)\in 2\Z$.  Furthermore, we have $\sum_{z_0\in\sigma(\pu)}\mz(\pu,z_0)\in2\Z$. Define
$$\sigma_{in}:=\{z\in\sigma(\pu):0<|z|<1\},\quad \sigma_{out}:=\{z\in\sigma(\pu):|z|>1\},\quad \sigma_{\T}:=\T\cap \sigma(\pu).$$
Note that
$$\mz(\pu,z_0)=\mz(\pu^\star,\ol{z_0}^{-1})=\mz(\pu,\ol{z_0}^{-1}),\quad \forall z_0\in\C\setminus\{0\}.$$
Moreover, the map $z\mapsto \ol{z}^{-1}$ is a bijection between $\sigma_{in}$ and $\sigma_{out}$. So
$$\sum_{z_0\in\sigma_{in}}\mz(\pu,z_0)=\sum_{\ol{z_0}^{-1}\in\sigma_{out}}\mz(\pu,z_0)=\sum_{z_0\in\sigma_{out}}\mz(\pu,\ol{z_0}^{-1})=\sum_{z_0\in\sigma_{out}}\mz(\pu,z_0).$$
It follows that
$$\sum_{z\in\T\cap \sigma(\pu)}\mz(\pu,z_0)=\sum_{z\in\sigma(\pu)}\mz(\pu,z_0)-\sum_{z_0\in\sigma_{in}}\mz(\pu,z_0)-\sum_{z_0\in\sigma_{out}}\mz(\pu,z_0)=\sum_{z\in\sigma(\pu)}\mz(\pu,z_0)-2\sum_{z_0\in\sigma_{in}}\mz(\pu,z_0)\in2\Z.$$
The proof is now complete.\ep

\bp[\textbf{Proof of Theorem~\ref{thm:dos}}]Necessity: Suppose $\pu$ is a Laurent polynomial with the DOS property with respect to the symmetry type $\sym\pu(z)=\eps z^c$. Then there are Laurent polynomials $\pu_1$ and $\pu_2$ with symmetry such that \er{dos} holds, which in particular implies $\pu^\star=\pu$. Write $\pu(z)=\sum_{k\in\Z}u(k)z^k$, we
have $u(k)=\ol{u(-k)}$. Since $\sym[\pu_1\pu^\star_1]=\sym[\pu_2\pu^\star_2]=1$, we have $\sym\pu=1$, and thus $u(k)=u(-k)$. Thus $\pu$ has real coefficients. This proves item (i).

By the definition of a Hermitian conjugate, we have
$$\pu_1^\star(z)=\frac{\ol{\pu_1(\ol{z})}}{\sym\pu_1(z)},\quad \pu^\star_2(z)=\eps z^c\frac{\ol{\pu_2(\ol{z})}}{\sym\pu_1(\ol{z})}.$$
So for $x\in\R\setminus\{0\}$, it follows from \er{dos} that
$$\pu(x)=\frac{|\pu_1(x)|^2-\eps x^c|\pu_2(x)|^2}{\sym\pu_1(x)}.$$
If $\eps=1$ and $c\in 2\Z+1$, we have $\eps x^c<0$ for all $x\in(-1,0)$, and thus $\mz(\pu,x)=2\min\{\mz(\pu_1,x),\mz(\pu_2,x)\}\in 2\Z$. This proves condition (ii), and conditions (iii) and (iv) can be verified similarly. This proves item (2).

Sufficiency: Suppose $\pu$ is a Laurent polynomial satisfying items (1) and (2). By \er{dos}, for every integer $k$, we have
$$(z^k\pu_1(z))(z^k\pu_1(z))^\star-\pu_2(z)\pu_2(z)^\star=\pu(z)\quad\mbox{and}\quad \frac{\sym[z^k\pu_1](z)}{\sym\pu_2(z)}=\eps z^{c+2k}.$$
That is, if $\pu$ has the DOS property with respect to symmetry type $\eps z^c$, then it has the DOS property with respect to symmetry type $\eps z^{c+2k}$ for all $k\in\Z$. As a consequence, it suffices to prove the sufficiency part for $c\in\{0,1\}$.

Since $\pu$ is Hermitian and has real coefficients, we conclude that $\sym\pu=1$ and $\mz(\pu,\ol{z_0}^{-1})=\mz(\pu^\star,z_0)=\mz(\pu,z_0)$. Also $\mz(\pu,z_0)=\mz(\pu,z_0^{-1})$ as $\pu$ has symmetry. It follows that \be\label{mz:z0}\mz(\pu,z_0)=\mz(\pu,z_0^{-1})=\mz(\pu,\ol{z_0}^{-1})=\mz(\pu,\ol{z_0}),\quad \forall z_0\in\C\setminus\{0\}.\ee
Moreover, we conclude from of Lemma~\ref{lem:rt:pm1} that
\be\label{mz:pm1}\mz(\pu,1),\mz(\pu,-1)\in2\Z.\ee
Define $\sigma(\pu)$ as \er{sigma:pu}, and define
$$\sigma_{\pm1}:=\sigma(\pu)\cap\{\pm 1\},\quad \sigma_{\R,in}:=\sigma(\pu)\cap\left((-1,0)\cup(0,1)\right),$$
$$\sigma_{\T,up}:=\sigma(\pu)\cap\{z\in\C\setminus\{0\}: |z|=1, \im(z)>0\},\quad \sigma_{in,up}:=\sigma(\pu)\cap\{z\in\C\setminus\{0\}: |z|<1, \im(z)>0\}.$$
By \er{mz:z0},\er{mz:pm1} and the fact that $\pu$ is Hermitian (particularly, $\fsupp(\pu)$ is symmetric about $0$), up to multiplication by a real number, $\pu(z)$ can be split into the product of the following four types of factors:
\begin{enumerate}
	\item[Type (a):] $z^{-1}(z-z_0)^2$, with $z_0\in\sigma_{\pm 1}$;
	
	\item[Type (b):] $z^{-2}(z-z_0)(z-z_0^{-1})(z-\ol{z_0})(z-\ol{z_0}^{-1})$, with $z_0\in\sigma_{in,up}$;

	\item[Type (c):] $z^{-1}(z-z_0)(z-\ol{z_0})$, with $z_0\in\sigma_{\T,up}$;
	
	\item[Type (d):] $z^{-1}(z-z_0)(z-z_0^{-1})$, with $z_0\in\sigma_{\R,up}$.
\end{enumerate}
According to Lemma~\ref{lem:dos:0}, we see that taking products preserves the DOS property of Laurent polynomials. Therefore, we only need to prove that all types of factors above have the DOS property with respect to the symmetry type $\eps z^c$. We discuss them one by one.

\begin{enumerate}
	\item[Type (a): ]If $z_0\in\sigma_{\pm 1}$, then
$$z^{-1}(z-1)^2=0^2-(z-1)(z-1)^\star,\quad z^{-1}(z+1)^2=(z+1)(z+1)^2-0^2.$$
Recall that $\pzr$ has symmetry of any type. Hence the type (a) factors have the DOS property with respect to symmetry type $\eps z^c$ with $\eps\in\{\pm 1\}$ and $c\in\{0,1\}$.

\item[Type (b): ]If $z_0\in\sigma_{in,up}$, then
\begin{align*}&z^{-2}(z-z_0)(z-z_0^{-1})(z-\ol{z_0})(z-\ol{z_0}^{-1})\\
=&\left(|z_0|^{-1}z^{-1}(z-z_0)(z-\ol{z_0})\right)\left(|z_0|^{-1}z^{-1}(z-z_0)(z-\ol{z_0})\right)^\star-0^2.\end{align*}
For the same reason as in the type (a) case, we see that the type (b) factors have the DOS property with respect to symmetry type $\eps z^c$ with $\eps\in\{\pm 1\}$ and $c\in\{0,1\}$.

\item[Type (c): ]If $z_0\in\sigma_{\T,up}$, we have $\pv(z):=z^{-1}(z-z_0)(z-\ol{z_0})=z-2\re(z_0)+z^{-1}$ and $\re(z_0)<1$. We now discuss symmetry conditions (i) - (iv) in item (2) one by one:
\begin{itemize}
	\item For (i), the symmetry type is $z^0=1$. Define
	 $$\pv_1(z):=\frac{z-4\re(z_0)+2+z^{-1}}{2\sqrt{2-2\re(z_0)}},\qquad  \pv_2(z):=\frac{z-2+z^{-1}}{2\sqrt{2-2\re(z_0)}}.$$
Direct calculation yields $\pv=\pv_1\pv_1^\star-\pv_2\pv_2^\star$ and $\frac{\sym\pv_1}{\sym\pv_2}=1$.

\item For (ii), the symmetry type is $z^1=z$. Define $\pv_1(z):=z+1$ and $\pv_2(z):=\sqrt{2+2\re(z_0)}.$ It is easy to see that $\pv=\pv_1\pv_1^\star-\pv_2\pv_2^\star$ and $\frac{\sym\pv_1}{\sym\pv_2}=z$.

\item For (iii), the symmetry type is $-z^0=-1$. Define
$$\pv_1(z):=\sqrt{\frac{1-\re(z_0)}{2}}(z+1),\qquad  \pv_2(z):=\sqrt{\frac{1-\re(z_0)}{2}}(z-1).$$
It follows that $\pv=\pv_1\pv_1^\star-\pv_2\pv_2^\star$ and $\frac{\sym\pv_1}{\sym\pv_2}=-1$.

\item For (iv), the symmetry type is $-z^1=-z$. Define $\pv_1(z):=\sqrt{2-2\re(z_0)}$ and $\pv_2(z):=1-z^{-1}.$ It is straightforward to verify that $\pv=\pv_1\pv_1^\star-\pv_2\pv_2^\star$ and $\frac{\sym\pv_1}{\sym\pv_2}=-z$.

\end{itemize}
This finishes the case for type (c).

\item[Type (d): ]If $z_0\in\sigma_{\R,in}$, in particular, $z_0\in(-1,0)\cup(0,1)$. Let
$$\pw(z):=z^{-1}(z-z_0)(z-z_0^{-1})=z-(z_0+z_0^{-1})+z^{-1}.$$
We again consider the following two cases:

\begin{itemize}
	\item $\mz(\pu,z_0)\in 2\Z$: Since $\pw$ has symmetry type $\sym\pw=1$, we see that $\mz(\pu,z_0^{-1})\in 2\Z$, and thus $\pw^2$ must divide $\pu$. So it suffices to show that $\pw^2$ has the DOS property. In fact, since $\pw=\pw^\star$, we have $\pw(z)^2=\pw(z)\pw^\star(z)-0^2.$ Therefore, $\pw^2$ has the DOS property with respect to symmetry type $\eps z^c$ with $\eps\in\{\pm 1\}$ and $c\in\{0,1\}$.

\item $\mz(\pu,z_0)\in2\Z+1$: We discuss conditions (i) - (iv) individually:

\begin{itemize}
	\item For (i): the symmetry type is $1$. Define $\pw_1(z):=\pw(z)+\frac{1}{4}$ and $\pw_2(z):=\pv(z)-\frac{1}{4}.$ Direct calculation yields $\pw=\pw_1\pw_1^\star-\pw_2\pw_2^\star$ and $\frac{\sym\pw_1}{\sym\pw_2}=1$.
	
	\item For (ii), the symmetry type is $z$. Moreover, we have $\mz(\pu,z_0)\in2\Z$ for all $z_0\in(-1,0)$ and thus any such points fail the assumption $\mz(\pu,z_0)\in2\Z+1$. So we must have $z_0\in(0,1)$. Now define $\pw_1(z):=z+1$ and $\pw_2(z):=\sqrt{2+z_0+z_0^{-1}}.$ It is easy to see that $\pw=\pw_1\pw_1^\star-\pw_2\pw_2^\star$ and $\frac{\sym\pw_1}{\sym\pw_2}=z$.
	
	\item For (iii), we have $\mz(\pu,z_0)\in2\Z$ for all $z_0\in(-1,0)\cup(0,1)$, which cannot happen as we assumed that $\mz(\pu,z_0)\in2\Z+1$.

	\item For (iv), the symmetry type is $-z$, and we must have $z_0\in (-1,0)$. Define $\pw_1(z):=\sqrt{2-z_0-z_0^{-1}}$ and $\pw_2(z):=z^{-1}-1.$ It is straightforward to verify that $\pw=\pw_1\pw_1^\star-\pw_2\pw_2^\star$ and $\frac{\sym\pw_1}{\sym\pw_2}=-z$.
\end{itemize}	
\end{itemize}
\end{enumerate}
This finishes the proof of the sufficiency part. The theorem is now proved.\ep

\subsection{Proof of Theorem~\ref{thm:SymSpecFactSym}}

In previous subsections, we have proved the special cases of Theorem~\ref{thm:SymSpecFactSym}, namely Theorems~\ref{thm:SymSpecFactSym:sp} and \ref{thm:SymSpecFactSym:gcd}. In this subsection, we are going to complete the proof of Theorem~\ref{thm:SymSpecFactSym}. Before we do this, we need one last lemma.

\begin{lemma}\label{lem:fac:lau}
	Suppose that $\pu$ and $\pv$ are two Laurent polynomials with symmetry, such that $\pv\pv^\star$ divides $\pu\pu^\star$. Then there exists a Laurent polynomial $\pd$ with symmetry, such that \be\label{pu:pv:pd}\pd\pd^\star=\frac{\pu\pu^\star}{\pv\pv^\star}.\ee
Furthermore, for every Laurent polynomial $\pd$ with symmetry satisfying \er{pu:pv:pd}, there exists $k\in\Z$ such that $\sym\pd(z)=z^{2k}\frac{\sym\pu(z)}{\sym\pv(z)}$.

\end{lemma}


\bp The existence of $\pd$ satisfying \er{pu:pv:pd} is a direct consequence of the equivalence between items (1) and (3) in \cite[Theorem~2.9]{han13}. As both $\pu\pu^\star$ and $\pv\pv^\star$ satisfy item (1) and hence, item (3) holds for both. Now it is trivial that $\frac{\pu\pu^\star}{\pv\pv^\star}$ satisfies item (3) and hence, item (1) holds for $\frac{\pu\pu^\star}{\pv\pv^\star}$, which is exactly \er{pu:pv:pd}.
Next, denote $\sym\pd(z)=\eps_dz^{c_d}$, $\sym\pu(z)=\eps_uz^{c_u}$ and $\sym\pv(z)=\eps_vz^{c_v}$  for some $\eps_d,\eps_u,\eps_v\in\{\pm 1\}$ and $c_d,c_u,c_v\in\Z$. By Lemma~\ref{lem:rt:pm1}, we have $\eps_d=(-1)^{\mz(\pd,1)}=(-1)^{\mz(\pu,1)-\mz(\pv,1)}=\eps_u/\eps_v.$ Moreover, by Proposition~\ref{prop:sym:lau}, we get
\begin{align*}\odd(c_d)&=\odd(\len(\pd))=\odd(\len(\pu)-\len(\pv))=\odd(c_u-c_v).\end{align*}
Therefore, we can find $k\in\Z$ such that $\sym\pd(z)=z^{2k}\frac{\sym\pu(z)}{\sym\pv(z)}$.
\ep

\bp[\textbf{Proof of Theorem~\ref{thm:SymSpecFactSym}}]

We first prove sufficiency. Assume that items (1) and (2) hold. If $\det(\pA)$ is identically zero, then by Theorem~\ref{thm:nf:sym}, there exist strongly invertible $2\times 2$ matrices $\pP$ and $\pQ$ of Laurent polynomials with compatible symmetry such that $\pA=\pP\DG(\pe_1,\mathbf{0})\pQ$, where all multiplications are compatible, and $\pe_1$ has symmetry. Define $\mrpA:=\pP^{-1}\pA\pP^{-\star}$, it is trivial that $\mrpA$ has compatible symmetry, and all multiplications are compatible. Since $\mrpA$ is Hermitian, one can conclude that $\mrpA=\DG(\pe_1,\mathbf{0})\pQ\pP^{-\star}=\begin{bmatrix}\mrpA_{1,1} & \mathbf{0}\\
\mathbf{0} &\mathbf{0}\end{bmatrix}$ for some Laurent polynomial $\mrpA_{1,1}$ having symmetry with $\sym\mrpA_{1,1}=\sym\pA_{1,1}=1$ and satisfying $\mrpA_{1,1}=\mrpA_{1,1}^\star$. By letting $\mrpU:=\begin{bmatrix}\mrpA_{1,1}+\frac{1}{4} & \mrpA_{1,1}-\frac{1}{4}\\
\mathbf{0} & \mathbf{0}\end{bmatrix}$, we see that $\mrpU$ has compatible symmetry. Moreover, we have $\mrpA=\mrpU\DG(1,-1)\mrpU^\star$, where all multiplications are compatible. Finally, by letting $\pU:=\pP\mrpU$, we see that $\pA=\pU\DG(1,-1)\pU^{\star}$ holds. Moreover, $\pU$ has compatible symmetry and satisfies \er{eq:Sym1}. This proves the case when $\det(\pA)\equiv 0$.\\

Now assume that $\det(\pA)$ is not identically zero. We construct $\pU$ in the following steps.
	
\begin{enumerate}
		\item[Step 1.] Denote $\pA_{j,k}$ the $(j,k)$-th entry of $\pA$ for $1\le j,k\le 2$. Since $\pA^\star=\pA$, we have $\pA_{1,1}^\star=\pA_{1,1}$. By Lemma~\ref{lem:dos:1}, we have $\mz(\pA_{1,1},1),\mz(\pA_{1,1},-1)\in2\Z$.
		Let
	 $$\alpha_1:=\min\left\{\frac{1}{2}\mz(\pA_{1,1},1),\mz(\pA_{1,2},1)\right\},\quad \alpha_2:=\min\left\{\frac{1}{2}\mz(\pA_{1,1},-1),\mz(\pA_{1,2},-1)\right\},$$
	and define $\pw(z):=(z-1)^{\alpha_1}(z+1)^{\alpha_2}$. Then $\pw$ divides the first row of $\pA$, and $\pw^\star$ divides the first column of $\pA$.  Moreover, $\pw(z)\pw^\star(z)=(-1)^{\alpha_1}(z-1)^{2\alpha_1}(z+1)^{2\alpha_2}z^{-\alpha_1-\alpha_2}$, which divides $\pA_{1,1}$. Define
	 \be\label{ptA}\pF:=\DG(\pw,1),\quad \ptA:=\pF^{-1}\pA\pF^{-\star}=\begin{bmatrix} \frac{\pA_{1,1}}{\pw\pw^\star} & \frac{\pA_{1,2}}{\pw}\\[0.2cm]
		 \frac{\pA_{2,1}}{\pw^\star} &\pA_{2,2}\end{bmatrix}.\ee
We see that $\ptA$ is a Hermitian matrix of Laurent polynomials with compatible symmetry.

\item[Step 2.] Define $\ptA_{j,k}$ to be the $(j,k)$-th entry of $\ptA$ for $1\le j,k\le 2$. Define
\be\label{ptp}\ptp:=\gcd(\ptA_{1,1},\ptA_{1,2},\ptA_{2,1},\ptA_{2,2}).\ee
By our choice of $\alpha_1$, we have
$$\min\{\mz(\ptA_{1,1},1),\mz(\ptA_{1,2},1)\}=\min\left\{\mz(\pA_{1,1},1)-2\alpha_1,\mz(\pA_{1,2},1)-\alpha_1\right\}=0.$$
Hence $\mz(\ptp,1)=0$. Similarly, we get $\mz(\ptp,-1)=0$. Moreover, for $z_0\in\C\setminus\{0,1,-1\}$, we have $\mz(\ptp,z_0)=\mz(\pp_0,z_0)$, where $\pp_0=\gcd(\pA_{1,1},\pA_{1,2},\pA_{2,1},\pA_{2,2}).$ Therefore, we conclude that $\ptp=\pp$, where $\pp$ is defined as \er{pp}. According to item (2), $\ptp$ has the DOS property with respect to the symmetry type $\sym\pd(z)\alpha(z)$. By Theorem~\ref{thm:dos}, we have $\ptp^\star=\ptp$, thus $\sym\ptp=1$. Define $\mrpA:=\frac{1}{\ptp}\ptA$, then $\mrpA$ is a Hermitian matrix of Laurent polynomials with compatible symmetry, with $\sym\mrpA=\sym\ptA$. Furthermore, we have
\be\label{det:mrpa}\det(\mrpA)=\frac{\det(\ptA)}{\ptp^2}=\frac{\det(\pA)}{\det(\pF)\det(\pF^\star)\ptp^2}=\frac{-\pd\pd^\star}{(\pw\ptp)(\pw\ptp)^\star}.\ee
By Lemma~\ref{lem:fac:lau}, there exists a Laurent polynomial $\mrpd$ with symmetry such that $\det(\mrpA)=\mrpd\mrpd^\star$.

\item[Step 3.] Define $\mrpA_{j,k}$ to be the $(j,k)$-th entry of $\mrpA$ for $1\le j,k\le 2$. Note that $\gcd(\mrpA_{1,1},\mrpA_{1,2},\mrpA_{2,1},\mrpA_{2,2})=1$. Thus $\mrpA$ satisfies all assumptions of Theorem~\ref{thm:SymSpecFactSym:gcd}, by which we can find matrices $\brpA$ and $\brpU$ of real Laurent polynomials with compatible symmetry, such that $\mrpA=\brpU\brpA\brpU^\star$, with all multiplications being compatible and $\det(\brpA)=C$ is a negative constant.

\item[Step 4.] We see that $\brpA$ satisfies all assumptions of Theorem~\ref{thm:SymSpecFactSym:sp}, and thus we can find a matrix $\pV$ of real Laurent polynomials with compatible symmetry such that $\brpA=\pV\DG(1,-1)\pV^\star$. Then by our previous discussions, we have
$$\pA=\pF\ptA\pF^\star =\ptp\pF\mrpA\pF^\star=\ptp\pF\brpU\brpA\brpU^\star\pF^\star=\ptp\pF\brpU\pV\DG(1,-1)\pV^\star\brpU^\star\pF^\star.$$
By letting $\ptU:=\pF\brpU\pV$, we see that $\ptU$ has compatible symmetry and all multiplications are compatible. Moreover, we have
$$\pA=\ptp\ptU\DG(1,-1)\ptU^\star\quad\text{and}\quad \det(\pA)=-\ptp\ptp^\star\det(\ptU)\det(\ptU)^\star.$$
By item (1), we have $\det(\ptU)\det(\ptU)^\star=\frac{\pd\pd^\star}{\ptp\ptp^\star}.$ Since $\det(\ptU)$ has symmetry, we conclude from Lemma~\ref{lem:fac:lau} that there exists $c\in\Z$ such that $\sym[\det(\ptU)](z)=z^{2c}\frac{\sym\pd(z)}{\sym\ptp(z)}=z^{2c}\sym\pd(z)$.\\
\end{enumerate}

From the proof of Theorem~\ref{thm:SymSpecFactSym:sp}, we can choose $\brpU$ such that $\sym\brpU=\sym\pth_1^\star\sym\pth$ with $\pth_1=[1,\alpha]$ and $\sym\pth\in\{[\alpha,1],[1,1],[\alpha,\alpha]\}$. Denote $\brpU_{j,k}$ the $(j,k)$-th entry of $\brpU$ for $1\le j,k\le 2$. It follows that
$$\frac{\sym\ptU_{1,1}}{\sym\ptU_{2,1}}=\frac{\sym\pA_{1,1}}{\sym\pA_{2,1}}=\alpha(z).$$
Define $\sym\ptU_{1,1}(z)=\eps_1z^{k_1}$ for some $\eps_1\in\{\pm 1\}$ and $k_1\in\Z$, we conclude that
\be\label{sym:ptu}\frac{\sym\ptU_{1,2}(z)}{\sym\ptU_{1,1}(z)}=\frac{\sym\ptU_{1,2}(z)\sym\ptU_{2,1}(z)}{\sym\ptU_{1,1}(z)\sym\ptU_{2,1}(z)}=\frac{\sym[\det(\ptU)](z)}{\sym\ptU_{1,1}(z)\sym\ptU_{2,1}(z)}=\frac{z^{2k}\sym\pd(z)}{(\sym\ptU_{1,1})^2}\frac{\sym\ptU_{1,1}(z)}{\sym\ptU_{2,1}(z)}=z^{2c-2k_1}\sym\pd(z)\alpha(z).\ee
Since $\ptp$ has the DOS property with respect to symmetry type $\sym\pd(z)\alpha(z)$, it also has the DOS property with respect to symmetry type $z^{2c-2k_1}\sym\pd(z)\alpha(z)$. Hence there exist real Laurent polynomials $\pp_1$ and $\pp_2$ with symmetry such that
\be\label{dos:pp}\pp_1(z)\pp_1^\star(z)-\pp_2(z)\pp_2^\star(z)=\ptp(z)\quad\mbox{and}\quad \frac{\sym\pp_1(z)}{\sym\pp_2(z)}=z^{2c-2k_1}\sym\pd(z)\alpha(z),\ee
which further implies that
\be\label{dos:pp:1}\DG(\ptp,-\ptp)=\begin{bmatrix}\pp_1	 &\pp_2^\star\\
\pp_2 &\pp_1^\star\end{bmatrix}\DG(1,-1)\begin{bmatrix}\pp_1&\pp_2^\star\\
\pp_2 &\pp_1^\star\end{bmatrix}^\star,\ee
and $\begin{bmatrix}\pp_1	 &\pp_2^\star\\
\pp_2 &\pp_1^\star\end{bmatrix}$ has compatible symmetry. Define
\be\label{pU}\pU:=\ptU\begin{bmatrix}\pp_1	 &\pp_2^\star\\
	\pp_2 &\pp_1^\star\end{bmatrix}.\ee
It is obvious that $\pA=\pU\DG(1,-1)\pU^\star$. From \er{sym:ptu} and \er{dos:pp}, we have $\sym\ptU_{1,1}\sym\pp_1=\sym\ptU_{1,2}\sym\pp_2$. Thus by item (4) of Proposition~\ref{sym:mat:com}, we see that the multiplication in \er{pU} is compatible, and $\pU$ has compatible symmetry satisfying \er{eq:Sym1}. This completes the proof of the sufficiency part.

We now prove necessity. Suppose $\pA=\pU\DG(1,-1)\pU^\star$, where $\pU$ is a matrix of real Laurent polynomials with compatible symmetry and \er{eq:Sym1} holds. By letting $\pd:=\det(\pU)$, we have $\det(\pA)=-\pd\pd^\star$ and $\pd$ has symmetry. This proves item (1).

Next, we prove item (2). Note that Step 1 in the proof of the sufficiency part does not require any assumptions in item (2), so we can apply that step here. Define $\ptA$ and $\pF$ as \er{ptA}. Then $\ptA$ is a matrix of real Laurent polynomials with compatible symmetry, and $\pA=\pF\ptA\pF^\star$ holds. Define $\ptp$ as \er{ptp}. Applying the same argument as in Step 2, one can conclude that $\ptp=\pp$, with $\pp$ being the same as in \er{pp}. So it only remains to show that $\ptp$ has the DOS property with respect to symmetry type $\sym\pd(z)\alpha(z)$.

As all entries of $\ptA$ have symmetry, by \cite[Lemma 2.4]{han13}, we see that $\ptp$ also has symmetry. Denote $\sym\ptp(z)=\eps_0z^{c_0}$ for some $\eps_0\in\{\pm 1\}$ and $c_0\in\Z$. From the discussions in Step 2, we see that $\mz(\ptp,1)=\mz(\ptp,-1)=0$, so Lemma~\ref{lem:rt:pm1} yields $\eps_0=1$ and $c_0\in 2\Z$. Since $\ptp$ is a greatest common divisor of all entries of $\ptA$ have symmetry, so is $z^{c_0/2}\ptp$. Thus without loss of generality, we can redefine $\ptp(z):=z^{c_0/2}\ptp(z)$, and get $\sym\ptp=1$.

For $z_0\in\C\setminus\{0\}$, by the definition of $\ptp$ and the fact that $\ptA$ is Hermitian, we conclude that $\mz(\ptp,\ol{z_0}^{-1})=\mz(\ptp,z_0).$ Therefore, there exist $\tilde{c}\in\Z$ and $\tilde{\eps}\in\C\setminus\{0\}$ such that $\ptp(z)=\tilde{\eps}z^{\tilde{c}}\ptp^\star(z)$. As $\sym\ptp=1$, we know that $\fsupp(\ptp)$ is a symmetric interval with center $0$, and thus $\tilde{c}=0$.  Note that $\ol{z}^{-1}=z$ for all $z\in\T$, in which case we have $|\ptp(z)|=|\tilde{\eps}||\ptp(z)|.$ Hence $\tilde{\eps}\in\T$. If $\tilde{\eps}\ne 1$, without loss of generality, we simply redefine $\ptp:=\tilde{\eps}^{-1/2}\ptp$, and thus $\ptp^\star=\ptp$. Note that multiplying $\tilde{\eps}^{-1/2}$ does not change the symmetry type of $\ptp$, so we still have $\sym\ptp=1$. As a result, we see that $\ptp$ must have real coefficients. Therefore, $\ptp$ satisfies item (1) of Theorem~\ref{thm:dos}.

Now define $\mrpA:=\frac{1}{\ptp}\ptA$, which is a Hermitian matrix of Laurent polynomials. By $\sym\ptp=1$, we see that $\mrpA$ has compatible symmetry with $\sym\mrpA=\sym\ptA$, and \er{det:mrpa} holds. By Lemma~\ref{lem:fac:lau}, there exists a Laurent polynomial $\mrpd$ with symmetry such that $\det(\mrpA)=-\mrpd\mrpd^\star$, that is item (1) of Theorem~\ref{thm:SymSpecFactSym} holds with $\pA$ and $\pd$ being replaced by $\mrpA$ and $\mrpd$ respectively. On the other hand, note that $\gcd(\mrpA_{1,1},\mrpA_{1,2},\mrpA_{2,1},\mrpA_{2,2})=1$. Using the fact that $1$ has the DOS property with respect to any symmetry type, we see that item (2) of Theorem~\ref{thm:SymSpecFactSym} holds with $\pA$ being replaced by $\mrpA$. Hence, by the sufficiency part of the theorem that has already been proved, we can find a $2\times 2$ matrix $\mrpU$ of Laurent polynomials with compatible symmetry, such that $\mrpA=\mrpU\DG(1,-1)\mrpU^\star$. It follows that
\begin{align*}
\pU\DG(1,-1)\pU^\star&=\pA=\pF\ptA\pF^\star=\ptp\pF\mrpA\pF^\star=\ptp\ptU\DG(1,-1)\ptU^\star,
\end{align*}
where $\ptU:=\pF\mrpU$ has compatible symmetry. Define $\pQ:=\adj(\ptU)\pU$ and $\ptd:=\det(\ptU)$. Then $\pQ$ has compatible symmetry and $\ptd$ has symmetry. Direct calculation yields
\begin{align*}
&\pQ\DG(1,-1)\pQ^\star=\ptp\adj(\ptU)\ptU\DG(1,-1)\ptU^\star\adj(\ptU)^\star=\ptp\ptd\ptd^\star\DG(1,-1).
\end{align*}
From the above equation, we get $\ptp\ptd\ptd^\star=\pQ_{1,1}\pQ_{1,1}^\star-\pQ_{1,2}\pQ_{1,2}^\star,$ where $\pQ_{j,k}$ denotes the $(j,k)$-th entry of $\pQ$ for $1\le j,k\le 2$. Denote $\sym\pU_{1,2}=\eps_{1,2}z^{c_{1,2}}$ with $\eps_{1,2}\in\{\pm1\}$ and $c_{1,2}\in\Z$. Then
\begin{align*}\frac{\sym\pQ_{1,1}(z)}{\sym\pQ_{1,2}(z)}&=\frac{\sym\pU_{1,1}(z)}{\sym\pU_{1,2}(z)}=\frac{\sym\pd(z)}{\sym\pU_{1,2}(z)\sym\pU_{2,2}(z)}=\frac{\sym\pd(z)}{(\sym\pU_{1,2}(z))^2}\frac{\sym\pU_{1,2}(z)}{\sym\pU_{2,2}(z)}=z^{-2c_{1,2}}\sym\pd(z)\alpha(z),\end{align*}
where we used $\sym\pd=\sym\pU_{1,1}\sym\pU_{2,2}$ (see \er{sym:det}) and $\frac{\sym\pU_{1,2}(z)}{\sym\pU_{2,2}(z)}=\alpha(z)$. Therefore, $\ptp\ptd\ptd^\star$ has the DOS property with respect to symmetry type $z^{-2c_{1,2}}\sym\pd(z)\alpha(z)$. In particular, $\ptp\ptd\ptd^\star$ satisfies item (2) of Theorem~\ref{thm:dos}. For $z_0\in\R\setminus\{0\}$, one can show that $\mz(\ptd^\star,z_0)=\mz(\ptd,\ol{z_0}^{-1})=\mz(\ptd,z_0^{-1})=\mz(\ptd,z_0)$, so $\mz(\ptd\ptd^\star,z_0)\in2\Z$, and $\odd(\mz(\ptp,z_0))=\odd(\mz(\ptp\ptd\ptd^\star,z_0))$ for all $z_0\in\R\setminus\{0\}$. This means that $\ptp$ also satisfies item (2) of Theorem~\ref{thm:dos}.

Consequently, by Theorem~\ref{thm:dos}, $\ptp$ has the DOS property with respect to symmetry type $z^{-2c_{1,2}}\sym\pd(z)\alpha(z)$, and hence also has the DOS property with respect to symmetry type $\sym\pd(z)\alpha(z)$. This proves item (2) of Theorem~\ref{thm:SymSpecFactSym}. This finishes the proof of the necessity part.

The proof of the theorem is now complete.
\ep

\section{Symmetric Quasi-tight Framelets with Two Generators}\label{sec:qtf}

In this section, we study quasi-tight framelet filter banks $\{a;b_1,b_2\}_{\Theta,(1,-1)}$ with symmetry property. Our goal is to prove the main result Theorem~\ref{thm:qtfsym}.

\subsection{A general guideline for constructing quasi-tight framelets}

Let $a,\Theta\in\lp{0}$ be such that $\Theta^\star=\Theta$ and let $b:=(b_1,b_2)^\tp\in\lrs{0}{2}{1}$. Define $\cM_{\pa,\pTh}$ as \er{cond:oep:tf}. Note that $\{a;b_1,b_2\}_{\Theta,(1,-1)}$ is a quasi-tight framelet if and only if \er{fac:1} holds with $s=2$, $\eps_1=1$ and $\eps_2=-1$, that is
\be\label{fac:2}\cM_{\pa,\pTh}(z)=\begin{bmatrix}\pb(z)& \pb(-z)\end{bmatrix}^{\star}\DG(1,-1)\begin{bmatrix}\pb(z),\pb(-z)\end{bmatrix}.\ee
Let $n_b$ be a positive integer satisfying \er{eq:nbRange}. Then we can define the $2\times 2$ matrix $\cM_{\pa,\pTh|n_b}$ of Laurent polynomials as \er{eq:DefcMnb}. More precisely, we have
\be\label{cm:nb}\cM_{\pa,\pTh|n_b}(z)=\begin{bmatrix}\pA(z)& \pB(z)\\
\pB(-z)	&\pA(-z) \end{bmatrix},\ee
where
\be\label{cm:nb:ab}\pA(z):=\dfrac{\pTh(z)-\pTh(z^2)\pa^\star(z)\pa(z)}{(1-z)^{n_b}(1-z^{-1})^{n_b}},\quad \pB(z):=\dfrac{-\pTh(z^2)\pa^\star(z)\pa(-z)}{(1-z^{-1})^{n_b}(1+z)^{n_b}}.\ee


Next, we are going to factorize $\cM_{\pa,\pTh|n_b}$ as the following:
\be\label{cm:ptb}\cM_{\pa,\pTh|n_b}(z)=\begin{bmatrix}\mrpb(z) & \mrpb(-z)\end{bmatrix}^\star\DG(1,-1)\begin{bmatrix}\mrpb(z) & \mrpb(-z)\end{bmatrix},\ee
for some $\mrb\in\lrs{0}{2}{1}$. To properly factorize $\cM_{\pa,\pTh|n_b}$, we need to introduce the notion of a coset sequence of a filter: For a matrix-valued filter $u\in\lrs{0}{t}{r}$ and $\gamma\in\Z$, define the $\gamma$-coset sequence of $u$ by
\be\label{coset:seq}u^{[\gamma]}(k)=u(\gamma+2k ),\qquad k\in\Z.\ee
Denote $\pu$ and $\pu^{[\gamma]}$ the Laurent polynomials associated with $\pu$ and $u^{[\gamma]}$ respectively, we have
\be\label{coset:seq:1}\pu(z)=\pu^{[0]}(z^2)+z\pu^{[1]}(z^2),\qquad z\in\C\setminus\{0\}.\ee
Furthermore, \er{coset:seq:1} is equivalent to
\be\label{coset:seq:0}\pu^{[0]}(z^2)=\frac{u(z)+u(-z)}{2},\quad \pu^{[1]}(z^2)=\frac{u(z)-u(-z)}{2z},\qquad z\in\C\setminus\{0\}.\ee
It follows from \er{coset:seq:0} that
\be\label{cm:cn}\begin{bmatrix}
	1& z^{-1}\\
	1 &-z^{-1}\end{bmatrix}^\star\cM_{\pa,\pTh|n_b}(z)\begin{bmatrix}
	1& z^{-1}\\
	1 &-z^{-1}\end{bmatrix}=\begin{bmatrix}
	\pA^{[0]}(z^2)+\pB^{[0]}(z^2)& \pA^{[1]}(z^2)-\pB^{[1]}(z^2)\\
	 z^2\left(\pA^{[1]}(z^2)+\pB^{[1]}(z^2)\right) &\pA^{[0]}(z^2)-\pB^{[0]}(z^2)\end{bmatrix},\ee
and we define
\be\label{N:a:Th:nb}\cN_{\pa,\pTh|n_n}(z):=\frac{1}{2}\begin{bmatrix}
	\pA^{[0]}(z)+\pB^{[0]}(z)& \pA^{[1]}(z)-\pB^{[1]}(z)\\
	 z\left(\pA^{[1]}(z)+\pB^{[1]}(z)\right) &\pA^{[0]}(z)-\pB^{[0]}(z)\end{bmatrix}.\ee
Now it suffices to find a matrix of Laurent polynomials $\pV=\begin{bmatrix}\pV_{1,1} &\pV_{1,2}\\
\pV_{2,1} &\pV_{2,2}\end{bmatrix}$ such that $\cN_{\pa,\pTh|n_b}=\pV^{\star}\DG(1,-1)\pV$. In this case, by letting
\be\label{mrpb:V}\mrpb^{[0]}(z):=\begin{bmatrix}\pV_{1,1}\\
\pV_{2,1} \end{bmatrix},\quad \mrpb^{[1]}(z):=\begin{bmatrix}\pV_{1,2}\\
\pV_{2,2} \end{bmatrix},\ee
we have
\be\label{cn:mrpb}\cN_{\pa,\pTh|n_b}(z)=\begin{bmatrix}\mrpb^{[0]}(z)&\mrpb^{[1]}(z)\end{bmatrix}^\star\DG(1,-1)\begin{bmatrix}\mrpb^{[0]}(z)&\mrpb^{[1]}(z)\end{bmatrix},\ee
which by \er{coset:seq:1} is equivalent to \er{cm:ptb}. By letting $b:=(b_1,b_2)^\tp\in\lrs{0}{2}{1}$ be such that $\pb(z)=(1-z)^{n_b}\mrpb(z)$, we see that \er{fac:2} holds, and thus $\{a;b_1,b_2\}_{\Theta,(1,-1)}$ is a quasi-tight framelet filter bank. Moreover, we have $\min\{\vmo(b_1),\vmo(b_2)\}\ge n_b$.

By the discussions above, we summarize the guideline for construction as the following theorem.

\begin{theorem}\label{thm:cons:qtf}Let $a,\Theta\in\lp{0}$ be such that $\Theta^\star=\Theta$ and let $b=(b_1,b_2)^\tp\in\lrs{0}{2}{1}$. Suppose that $\{a;b_1,b_2\}_{\Theta,(1,-1)}$ is a quasi-tight framelet filter bank such that $\min\{\vmo(b_1),\vmo(b_2)\}\ge n_b$. Then \er{cn:mrpb} holds, where $\cN_{\pa,\pTh|n_b}$ is defined as \er{N:a:Th:nb}, and $\mrpb\in \lrs{0}{2}{1}$ satisfies $\pb(z)=(1-z)^{n_b}\mrpb(z)$.
Conversely, suppose that there exists a matrix of Laurent polynomials $\pV=\begin{bmatrix}\pV_{1,1} &\pV_{1,2}\\
\pV_{2,1} &\pV_{2,2}\end{bmatrix}$ such that $\cN_{\pa,\pTh|n_b}=\pV^{\star}\DG(1,-1)\pV$. Define $\mrb,b\in\lrs{0}{2}{1}$ via \er{mrpb:V} and $\pb(z):=(1-z)^{n_b}\mrpb(z)$ respectively. Then $\{a;b_1,b_2\}_{\Theta,(1,-1)}$ is a quasi-tight framelet filter bank such that $\min\{\vmo(b_1),\vmo(b_2)\}\ge n_b$.

\end{theorem}

\subsection{Symmetry property of quasi-tight framelet filter banks}

We need to investigate a little further on symmetry property of quasi-tight framelet filter banks.

\begin{theorem}\label{thm:qtf:fb}Let $a, \Theta\in \lp{0}$ be such that $\Theta^\star=\Theta$ and let $b=(b_1,b_2)^\tp\in\lrs{0}{2}{1}$. Suppose that $\{a;b_1, b_2\}_{\Theta,(1,-1)}$ is a quasi-tight framelet filter bank such that
	 \be\label{qtf:sym}\sym\pTh(z)=1,\quad \sym\pa(z)=\eps z^c,\quad \sym\pb_1(z)=\eps_1z^{c_1},\quad \sym\pb_2(z)=\eps_2z^{c_2}, \ee
for some $\eps,\eps_1,\eps_2\in\{\pm 1\}$ and $c,c_1,c_2\in\Z$. Here $\pTh, \pa,\pb_1$ and $\pb_2$ are Laurent polynomials associated with $\Theta, a, b_1$ and $b_2$ respectively. Then
\be\label{sym:cent}c_1+c,\quad c_2+c\in 2\Z.\ee
\end{theorem}

\bp By the definition of a quasi-tight framelet filter bank, we have
\be\label{def:qtf}\pTh(z^2)\pa(z)\pa^\star(-z)+\pb_1(z)\pb_1^\star(-z)-\pb_2(z)\pb_2^\star(-z)=0,\quad z\in \C\setminus\{0\}.\ee
Note that
$$\sym[\pTh(\cdot^2)\pa\pa^\star(-\cdot)](z)=\sym\pTh(z^2)\sym\pa(z)\sym\pa^\star(-z)=\eps z^c\eps(-z)^{-c}=(-1)^c,$$
$$\sym[\pb_1\pb_1^\star(-\cdot)](z)=\sym\pb_1(z)\sym\pb_1^\star(-z)=\eps_1z^{c_1}\eps_1(-z)^{-c_1}=(-1)^{c_1},$$
$$\sym[\pb_2\pb_2^\star(-\cdot)](z)=\sym\pb_2(z)\sym\pb_2^\star(-z)=\eps_2z^{c_2}\eps_2(-z)^{-c_2}=(-1)^{c_2}.$$
We claim that $(-1)^c=(-1)^{c_1}=(-1)^{c_2}$, i.e., $c,c_1,c_2$ have the same parity, thus \er{sym:cent} follows immediately. Assume otherwise. If $(-1)^c=(-1)^{c_1}\ne(-1)^{c_2}$, then we have
	 $$\sym[\pTh(\cdot^2)\pa\pa^\star(-\cdot)+\pb_1\pb_1^\star(-\cdot)]=(-1)^c\ne \sym[\pb_2\pb_2^\star(-\cdot)],$$
which contradicts \er{def:qtf}. Similarly $(-1)^c\ne(-1)^{c_1}=(-1)^{c_2}$ or $(-1)^c=(-1)^{c_2}\ne(-1)^{c_1}$ cannot happen either. This proves the claim, and finishes the proof of the theorem.\ep

Our construction of quasi-tight framelet filter banks is based on factorizing the matrix $\cN_{\pa,\pTh|n_b}$, whose entries are Laurent polynomials associated with coset sequences. So we need the following result on symmetry property of coset sequences.

\begin{lemma}\label{lem:sym:coset}Let $\pu\in \lp{0}$ be such that $\sym\pu(z)=\eps z^c$ for some $\eps\in\{\pm1\}$ and $c\in\Z$. Then the following statements hold:
\begin{enumerate}
	\item[(1)] If $c\in 2\Z$, then $\pu^{[0]}$ and $\pu^{[1]}$ have symmetry, with
$$\sym\pu^{[0]}(z)=\eps z^{c/2},\quad \sym\pu^{[1]}(z)=\eps z^{c/2-1};$$

\item[(2)] If $c\in 2\Z+1$, for all $l\in\Z$, define
\be\label{pp:k}\pP_l(z):=\frac{1}{\sqrt{2}}\begin{bmatrix}1 & z^l\\
1 &-z^l\end{bmatrix},\quad l\in\Z,\ee
and
$$\begin{bmatrix}\pv_1(z)\\
\pv_2(z)\end{bmatrix}:=\pP_l(z)\begin{bmatrix}\pu^{[0]}(z)\\
\pu^{[1]}(z)\end{bmatrix}=\frac{1}{\sqrt{2}}\begin{bmatrix}\pu^{[0]}(z)+z^l\pu^{[1]}(z)\\
\pu^{[0]}(z)-z^l\pu^{[1]}(z)\end{bmatrix}.$$
Then $\pv_1$ and $\pv_2$ have symmetry, with
$$\sym\pv_1(z)=\eps z^{\frac{c-1}{2}+l},\quad \sym\pv_2(z)=-\eps z^{\frac{c-1}{2}+l}.$$
\end{enumerate}
On the other hand, if $\pv_1$ and $\pv_2$ are Laurent polynomials with symmetry, then the following statements hold:
\begin{enumerate}
	\item[(3)] If $\sym\pv_1(z)=\eps z^{c+1}$ and $\sym\pv_2(z)=\eps z^{c}$ for some $\eps\in\{\pm1\}$ and $c\in\Z$, then $\pu(z):=\pv_1(z^2)+z\pv_2(z^2)$ also has symmetry, with $\sym\pu(z)=\eps z^{2c+2}$;
	
	\item[(4)] If $\sym\pv_1(z)=\eps z^{c}$ and $\sym\pv_2(z)=-\eps z^{c}$ for some $\eps\in\{\pm1\}$ and $c\in\Z$, define $\pP_l$ for all $l\in\Z$ as \er{pp:k} and
$$\begin{bmatrix}\pw_1(z)\\
\pw_2(z)\end{bmatrix}:=\pP_l(z)^{-1}\begin{bmatrix}\pv_1(z)\\
\pv_2(z)\end{bmatrix}=\frac{1}{\sqrt{2}}\begin{bmatrix}\pv_1(z)+\pv_2(z)\\
z^{-l}\left(\pv_1(z)-\pv_2(z)\right)\end{bmatrix}.$$
Then $\pu(z):=\pw_1(z^2)+z\pw_2(z^2)$ also has symmetry, with $\sym\pu(z)=\eps z^{2(c-l)+1}$.
\end{enumerate}
\end{lemma}

\bp All claims can be verified directly by using \er{coset:seq:0}.\ep

\subsection{Proof of Theorem~\ref{thm:qtfsym}}
Based on Theorem~\ref{thm:SymSpecFactSym} and the above discussions, we can give a complete characterization of quasi-tight framelets with symmetry.

\bp[\textbf{Proof of Theorem~\ref{thm:qtfsym}}]
Define $\pA$ and $\pB$ as \er{cm:nb:ab}. By the choice of $n_b$, we see that $\pA$ and $\pB$ are well-defined Laurent polynomials. Furthermore, from the symmetry assumptions of $\pa$ and $\pTh$, we have $\sym\pA(z)=1$ and $\sym\pB(z)=(-1)^{c+n_b}$. Hence by item (1) of Lemma~\ref{lem:sym:coset}, we get
\be\label{sym:pa:pb}\sym\pA^{[0]}(z)=1,\quad \sym\pA^{[1]}(z)=z^{-1},\quad \sym\pB^{[0]}(z)=(-1)^{c+n_b},\quad \sym\pB^{[1]}(z)=(-1)^{c+n_b}z^{-1}.\ee

We first prove necessity.
Suppose there exist a quasi-tight framelet filter bank $\{a;b_1,b_2\}_{\Theta,(1,-1)}$ with symmetry with $n_b$ order of vanishing moments for some $b=(b_1,b_2)^\tp\in\lrs{0}{2}{1}$. Denote $\sym\pb_1(z)=\eps_1z^{c_1}$ and $\sym\pb_2=\eps_2z^{c_2}$. Since $b_1$ and $b_2$ have at least $n_b$ order of vanishing moments, there exists $\mrpb=(\mrpb_1,\mrpb_2)^\tp\in\lrs{0}{2}{1}$ such that
	 \be\label{mrpb}\pb_l(z)=(1-z)^{n_b}\mrpb_l(z),\quad l=1,2.\ee
From the symmetry types of $\pb_1$ and $\pb_2$, we have $\sym\mrpb_1(z)=\eps_1z^{c_1+n_b}$ and $\sym\mrpb_2(z)=\eps_2z^{c_2+n_b}$.

If $c+n_b\in 2\Z$, it follows from \er{sym:pa:pb} that $\sym\cN_{\pa,\pTh|n_b}(z)=\begin{bmatrix}1 &z^{-1}\\
z &1\end{bmatrix}$.	By Theorem~\ref{thm:qtf:fb}, we see that $c_1+n_b, c_2+n_b\in 2\Z$. Hence by item (1) of Lemma~\ref{lem:sym:coset}, it is not hard to verify that the matrix
\be\label{pU:mrpb}\pU(z)=\begin{bmatrix}\mrpb^{[0]}(z) & \mrpb^{[1]}(z)\end{bmatrix}^\star\ee
has compatible symmetry, with $\sym\pU(z)=\begin{bmatrix}\eps_1z^{-(c_1+n_b)/2} & \eps_2z^{-(c_2+n_b)/2}\\
\eps_1z^{-(c_1+n_b)/2+1} & \eps_2z^{-(c_2+n_b)/2+1}\end{bmatrix}$. In particular, we have $\frac{\sym\pU_{1,1}(z)}{\sym\pU_{2,1}(z)}=\frac{\sym\pU_{1,2}(z)}{\sym\pU_{2,2}(z)}=z^{-1}$. According to Theorem~\ref{thm:cons:qtf}, $\{a;b_1,b_2\}_{\Theta,(1,-1)}$ is a quasi-tight framelet filter bank with $n_b$ order of vanishing moments implies that $\cN_{\pa,\pTh|n_b}=\pU\DG(1,-1)\pU^\star$. Thus by Theorem~\ref{thm:SymSpecFactSym}, items (1) and (2) follows immediately.

If  $c+n_b\in 2\Z+1$, for any integer $l$, define $\pP_l$ as \er{pp:k}, and define $\pN:=\pP_l\cN_{\pa,\pTh|n_b}\pP_l^\star$. By the definition of $\cN_{\pa,\pTh|n_b}$ and \er{sym:pa:pb}, we see that $\pN$ has compatible symmetry with $\sym\pN=\begin{bmatrix}1 &-1\\
-1 & 1\end{bmatrix}$. Define $\pU$ as \er{pU:mrpb}. By Theorem~\ref{thm:qtf:fb}, the condition $c+n_b\in 2\Z+1$ implies that $c_1+n_b, c_2+n_b\in 2\Z+1$. Hence by item (2) of Lemma~\ref{lem:sym:coset}, $\ptU:=\pP_l\pU$ has compatible symmetry, with
$$\sym\ptU(z)=\begin{bmatrix}\eps_1z^{-(c_1+n_b-1)/2+l} &\eps_2z^{-(c_2+n_b-1)/2+l}\\
-\eps_1z^{-(c_1+n_b-1)/2+l} &-\eps_2z^{-(c_2+n_b-1)/2+l}\end{bmatrix}.$$
In particular, $\frac{\sym\ptU_{1,1}(z)}{\sym\ptU_{2,1}(z)}=\frac{\sym\ptU_{1,2}(z)}{\sym\ptU_{2,2}(z)}=-1$. Using the assumption that $\{a;b_1,b_2\}_{\Theta,(1,-1)}$ is a quasi-tight framelet filter bank, we have
$$\pN=\pP_l\cN_{\pa,\pTh|n_b}\pP_l^\star=\pP_l\pU\DG(1,-1)\pU^\star\pP_l^\star=\ptU\DG(1,-1)\ptU^\star.$$
Since $\pP_l$ is strongly invertible, the Smith normal forms of $\pN$ and $\cN_{\pa,\pTh|n_b}$ are the same. Thus the gcd of the 4 entries of $\cN$, as well as $\det(\cN)$, must by the same as those of $\cN_{\pa,\pTh|n_b}$ up to multiplications by monomials. Therefore, by Theorem~\ref{thm:SymSpecFactSym}, items (1) and (2) must hold. This finishes the proof of the necessity part.

We now prove sufficiency. Suppose items (1) and (2) hold. If $c+n_b\in 2\Z$, by the definition of $\cN_{\pa,\pTh|n_b}$ and \er{sym:pa:pb}, we have $\sym\cN_{\pa,\pTh|n_b}(z)=\begin{bmatrix}1 &z^{-1}\\
z &1\end{bmatrix}$. By Theorem~\ref{thm:SymSpecFactSym}, there exists a $2\times 2$ matrix $\pU$ of Laurent polynomials with compatible symmetry, such that $\frac{\sym\pU_{1,1}(z)}{\sym\pU_{2,1}(z)}=\frac{\sym\pU_{1,2}(z)}{\sym\pU_{2,2}(z)}=z^{-1}$ and $\cN_{\pa,\pTh|n_b}=\pU\DG(1,-1)\pU^\star$. Define $\mrb=(\mrb_1,\mrb_2)^\tp\in\lrs{0}{2}{1}$ via $\begin{bmatrix}\mrpb^{[0]} & \mrpb^{[1]}\end{bmatrix}:=\pU^\star$. According to Lemma~\ref{lem:sym:coset}, both $\mrb_1$ and $\mrb_2$ have symmetry. Define $b=(b_1,b_2)^\tp\in\lrs{0}{2}{1}$ via $\pb(z):=(1-z)^{n_b}\mrpb(z)$, we can easily verify that $b_1$ and $b_2$ have symmetry. Moreover, Theorem~\ref{thm:cons:qtf} tells that $\{a;b_1,b_2\}_{\Theta,(1,-1)}$ is a quasi-tight framelet filter bank with $\min\{\vmo(b_1),\vmo(b_2)\}\ge n_b$.

If $c+n_b\in 2\Z+1$, for any integer $l$, define $\pP_l$ as \er{pp:k}, and define $\pN:=\pP_l\cN_{\pa,\pTh|n_b}\pP_l^\star$. By the definition of $\cN_{\pa,\pTh|n_b}$ and \er{sym:pa:pb}, it is easy to see that $\pN$ has compatible symmetry with $\sym\pN=\begin{bmatrix}1 &-1\\
-1 & 1\end{bmatrix}$. We know that $\pN$ and $\cN_{\pa,\pTh|n_b}$ have the same Smith normal form (up to multiplication by monomials). Thus items (1) and (2) of Theorem~\ref{thm:SymSpecFactSym} holds with $\pA=\pN$. Hence by Theorem~\ref{thm:SymSpecFactSym}, there exists a matrix $\ptU$ of Laurent polynomials with compatible symmetry, such that $\pN=\ptU\DG(1,-1)\ptU^\star$ and $\frac{\sym\ptU_{1,1}(z)}{\sym\ptU_{2,1}(z)}=\frac{\sym\ptU_{1,2}(z)}{\sym\ptU_{2,2}(z)}=-1$. Define $\mrb=(\mrb_1,\mrb_2)^\tp\in\lrs{0}{2}{1}$ via $\begin{bmatrix}\mrpb^{[0]} & \mrpb^{[1]}\end{bmatrix}:=\ptU^\star\pP_l^{-\star}$. By Lemma~\ref{lem:sym:coset}, both $\mrb_1$ and $\mrb_2$ have symmetry. Moreover, we see that \er{cn:mrpb} holds. Define $b=(b_1,b_2)^\tp\in\lrs{0}{2}{1}$ via $\pb(z):=(1-z)^{n_b}\mrpb(z)$, it follows from Theorem~\ref{thm:cons:qtf} that $\{a;b_1,b_2\}_{\Theta,(1,-1)}$ is a quasi-tight framelet filter bank with $\min\{\vmo(b_1),\vmo(b_2)\}\ge n_b$. This finishes the proof of the sufficiency part.

If in addition $\pa(1)=\pTh(1)=1$ and $\phi\in \Lp{2}$, where
	$\wh{\phi}(\xi):=\prod_{j=1}^\infty \pa(e^{-i 2^{-j}\xi})$ for $\xi\in \R$, define $\eta, \psi^1,\psi^2$ as \er{eq:psi}. It follows immediately that $\phi,\eta,\psi^1,\psi^s$ have symmetry, and $\{\eta,\phi;\psi^1,\psi^2\}_{(1,-1)}$ is a quasi-tight framelet in $\Lp{2}$ by the Oblique Extension Principle. Moreover, $\psi^1$ and $\psi^2$ have at least order $n_b$ vanishing moments.

The proof is now complete.\ep

\begin{rem}
From  the remark after \cite[Theorem~4.2]{han13}, we observe that $\det(\cN_{\pa, \pTh|n_b}(z))=(1-z)^{-n_b}(1-z^{-1})^{-n_b} \det(\cN_{\pa, \pTh|0}(z))$. Consequently, the Laurent polynomial $\pp_0$ in \eqref{pp0:z} of Theorem~\ref{thm:qtfsym} is independent of the choice of $n_b$ and its symmetry type $(-1)^{c+n_b} z^{\odd(c+n_b)-1}\sym \pd_{n_b}(z)$ only differs by a factor of $z^{2l}$ for some $l\in \Z$.
Hence, items (1) and (2) of Theorem~\ref{thm:qtfsym} hold for some $n_b$ satisfying \er{eq:nbRange} if and only if item (1) and (2) hold with $n_b=0$.
\end{rem}


\begin{thebibliography}{10}
	
	
	
	
	
	
	
\bibitem{cs08}
M.~Charina and J.~St\"ockler, Tight wavelet frames for irregular multiresolution analysis. \emph{Appl. Comput. Harmon. Anal.} \textbf{25} (2008), 98--113.

	\bibitem{ch00}
	C.~K.~Chui and W.~He, Compactly supported tight frames associated with refinable functions. \emph{Appl. Comput. Harmon. Anal.} \textbf{8} (2000), 293--319.
	
	
	
	\bibitem{chs02}
	C.~K.~Chui, W.~He, and J.~St\"ockler, Compactly supported tight and sibling frames with maximum vanishing moments. \emph{Appl. Comput. Harmon. Anal.} \textbf{13} (2002), 224--262.
	
	
	
	
	\bibitem{daubook}
	I.~Daubechies, Ten lectures on wavelets. CBMS-NSF Series in Applied Mathematics, \textbf{61}. SIAM, 1992.
	
	
	\bibitem{dh04}
	I.~Daubechies and B.~Han, Pairs of dual wavelet frames from any two refinable functions. \emph{Constr. Approx.} \textbf{20} (2004), 325--352.
	
	\bibitem{dhrs03}
	I.~Daubechies, B.~Han, A.~Ron, Z.~Shen, Framelets: MRA-based constructions of wavelet frames. \emph{Appl. Comput. Harmon. Anal.} \textbf{14} (2003), 1--46.
	
	
		\bibitem{dhacha}
	C.~Diao and B.~Han, Quasi-tight framelets with high vanishing moments derived from arbitrary refinable functions, \emph{Appl. Comput. Harmon. Anal.} \textbf{49} (2020), 123--151.
	
	\bibitem{dh18pp}
	C.~Diao and B.~Han,
	Generalized matrix spectral factorization and quasi-tight framelets with minimum number of generators, \emph{Math. Comp.} \textbf{89} (2020), 2867--2911.
	
\bibitem{ds07}
B. Dong and Z. Shen, Pseudo-splines, wavelets and framelets. \emph{Appl. Comput. Harmon. Anal.}
\textbf{22} (2007), 78--104,

	
	
	\bibitem{eh08}
	M.~Ehler and B.~Han, Wavelet bi-frames with few generators from multivariate refinable functions, \emph{Appl. Comput. Harmon. Anal.}, \textbf{25} (2008), 407--414.
	
	
	
	
	\bibitem{han97}
	B.~Han, On dual wavelet tight frames. \emph{Appl. Comput. Harmon. Anal.} \textbf{4} (1997), 380--413.
	
	
	
	
	\bibitem{han03jat}
	B.~Han, Vector cascade algorithms and refinable function vectors in Sobolev spaces. \emph{J. Approx. Theory} \textbf{124} (2003), 44--88.
	
	
	\bibitem{han09}
	B.~Han, Dual multiwavelet frames with high balancing order and compact fast frame transform. \emph{Appl. Comput. Harmon. Anal.} \textbf{26} (2009), 14--42.


	
	
	
	\bibitem{han12}
	B.~Han, Nonhomogeneous wavelet systems in high dimensions. \emph{Appl. Comput. Harmon. Anal.} \textbf{32} (2012), 169--196.
	
	\bibitem{han13}
	B.~Han, Matrix splitting with symmetry and symmetric tight framelet filter banks with two high-pass filters. \emph{Appl. Comput. Harmon. Anal.} \textbf{35} (2013), 200--227.
	
	
	
	\bibitem{han14acha}
	B.~Han, Symmetric tight framelet filter banks with three high-pass filters. \emph{Appl. Comput. Harmon. Anal.} \textbf{37} (2014), 140--161.
	

\bibitem{han17}
B.~Han, Homogeneous wavelets and framelets with the refinable structure. \emph{Sci. China Math.} \textbf{60} (2017), 2173--2198.
	
	\bibitem{hanbook}
	B.~Han, Framelets and wavelets: algorithms, analysis and applications, \emph{Applied and Numerical Harmonic Analysis}, Birkh\"auser/Springer, Cham, (2017), 724 pages.



	
	
	
	
	
		\bibitem{hl19pp}
	B.~Han and R.~Lu, Compactly supported quasi-tight multiframelets with high balancing orders and compact framelet transform. \emph{Appl. Comput. Harmon. Anal.}, \textbf{51} (2021), 295-332.
	
	\bibitem{hl20pp}
	B.~Han and R.~Lu, Multivariate quasi-tight framelets with high balancing orders derived from any compactly supported refinable vector functions. \emph{Sci. China Math.} (2021), https://doi.org/10.1007/s11425-020-1786-9.
	
	
	
	\bibitem{hm04}
	B.~Han and Q.~Mo, Splitting a matrix of Laurent polynomials with symmetry and its application to symmetric framelet filter banks. \emph{SIAM J. Matrix Anal. Appl.} \textbf{26} (2004), 97--124.
	
	\bibitem{hm05}
	B.~Han and Q.~Mo,
	Symmetric MRA tight wavelet frames with three generators and high vanishing moments.	 \emph{Appl. Comput. Harmon. Anal.} \textbf{18} (2005), 67--93.
	
	
	
	\bibitem{hz10}
	B.~Han and X.~Zhuang,  Matrix extension with symmetry and its application to symmetric orthonormal multiwavelets. \emph{SIAM J. Math. Anal.} \textbf{42} (2010), 2297--2317.
	

\bibitem{hjs04}	
D. P. Hardin, T. A. Hogan, and Q. Sun, The matrix-valued Riesz lemma and local orthonormal
bases in shift-invariant spaces. \emph{Adv. Comput. Math.} \textbf{20} (2004), 367--384,

	
	\bibitem{hr08}
	Y.~Hur and A.~Ron, L-CAMP: extremely local high-performance wavelet representations in high spatial dimension. \emph{IEEE Trans. Inform. Theory} \textbf{54} (2008), 2196--2209.
	
	
	
	\bibitem{jiang03}
	Q.~T.~Jiang, Parameterizations of masks for tight affine frames with two symmetric/antisymmetric generators. Frames. \emph{Adv. Comput. Math.} \textbf{18} (2003), 247--268.
	
	\bibitem{js15}
	Q.~T.~Jiang and Z.~Shen, Tight wavelet frames in low dimensions with canonical filters. \emph{J. Approx. Theory} \textbf{196} (2015), 55--78.
	
	
	\bibitem{kps16book}
	A.~Krivoshein, V.~Protasov, and M.~Skopina, Multivariate wavelet frames. \emph{Industrial and Applied Mathematics}. Springer, Singapore, 2016.
	
	
	
	
	
	
	\bibitem{ml11}
	Q.~Mo and S.~Li, Symmetric tight wavelet frames with rational coefficients. \emph{Appl. Comput. Harmon. Anal.} \textbf{31} (2011), 249--263.
	
	
	\bibitem{mz12}
	Q.~Mo and X.~Zhuang, Matrix splitting with symmetry and dyadic framelet filter banks over algebraic number fields. \emph{Linear Algebra Appl.} \textbf{437} (2012), 2650--2679.
	
	
	
	\bibitem{rs97}
	A.~Ron and Z.~Shen, Affine systems in $L_2(\R^d)$: the analysis of the analysis operator. \emph{J. Funct. Anal.} \textbf{148} (1997), 408--447.
	
	
	
	

\bibitem{sz18}	
A.~San Antol\'in and R.~A.~Zalik, Compactly supported Parseval framelets with symmetry associated to $E^{(2)}_d(\Z)$ matrices. \emph{Appl. Math. Comput.} \textbf{325} (2018), 179--190.

	\bibitem{sel01}
	I.~W.~Selesnick, Smooth wavelet tight frames with zero moments. \emph{Appl. Comput. Harmon. Anal.} \textbf{10} (2001), 163--181.
	
	
	\bibitem{sa04}
	I.~W.~Selesnick and A.~F.~Abdelnour, Symmetric wavelet tight frames with two generators. \emph{Appl. Comput. Harmon. Anal.} \textbf{17} (2004),  211--225.

\bibitem{sl11}
Y.~Shen and S.~Li, Wavelets and framelets from dual pseudo splines. \emph{Sci. China Math.} \textbf{54} (2011), 1233--1242.
	
	
	
	

	\bibitem{zhuang12}
	X.~Zhuang, Matrix extension with symmetry and construction of biorthogonal multiwavelets with any integer dilation. \emph{Appl. Comput. Harmon. Anal.} \textbf{33} (2012), 159--181.
	
\end{thebibliography}
\end{document}